\documentclass{imamat}

\usepackage[english]{babel}
\usepackage{amsmath}
\usepackage{amssymb}
\usepackage{mathtools}
\usepackage{enumerate}
\usepackage{graphicx}
\usepackage{array} 
\usepackage{url}
\usepackage{color}
\usepackage{nicefrac}
\usepackage[it,IT]{subfigure}
\usepackage{mathrsfs}
\usepackage{textgreek}

\newcommand{\R}{\ensuremath{\mathbb{R}}}
\newcommand{\N}{\ensuremath{\mathbb{N}}}
\newcommand{\Z}{\ensuremath{\mathbb{Z}}}
\newcommand{\Hil}{\ensuremath{\mathcal{H}}}
\newcommand{\V}{\ensuremath{\mathcal{V}}}
\newcommand{\A}{\ensuremath{\mathcal{A}}}
\newcommand{\eps}{\ensuremath{\varepsilon}}
\renewcommand{\phi}{\ensuremath{\varphi}}
\providecommand{\abs}[1]{\lvert#1\rvert}
\providecommand{\norm}[1]{\lVert#1\rVert}
\providecommand{\Norm}[1]{\left\lVert#1\right\rVert}
\newcommand{\Nf}{N^\mathrm{f}}

\newcommand{\HE}{E}    

\newcommand{\degK}{{}^\circ\hspace*{-0.1em}K}
\renewcommand{\degK}{K}

\newcommand{\tend}{{t_\mathrm{m}}}

\newcommand{\Omeps}{\Omega_{\eps}}
\newcommand{\Gameps}{\Gamma_{\eps}}
\newcommand{\Taeps}{T_{\alpha,\eps}}
\newcommand{\Toneeps}{T_{1,\eps}}
\newcommand{\Ttwoeps}{T_{2,\eps}}
\newcommand{\HEaeps}{\HE_{\alpha,\eps}}
\newcommand{\HEoneeps}{\HE_{1,\eps}}
\newcommand{\HEtwoeps}{\HE_{2,\eps}}
\newcommand{\Omoneeps}{\Omega_{\eps}^1}
\newcommand{\Omtwoeps}{\Omega_{\eps}^2}
\newcommand{\half}{\frac{1}{2}}
\newcommand{\bigC}{C}
\newcommand{\Tout}{T_\mathrm{a}}
\newcommand{\Tcrit}{T_\mathrm{c}}
\newcommand{\defeq}{:=}
\newcommand{\leavethisout}[1]{}
\newcommand{\normvec}{\mathbf{n}}
\newcommand{\ds}{\displaystyle}
\newcommand{\Rgas}{{\cal R}}
\newcommand{\Mmicro}{m} 
\newcommand{\Mmacro}{M} 
\newcommand{\Rtree}{R_\mathrm{tree}}
\newcommand{\Rf}{R^\mathrm{f}}
\newcommand{\Lf}{L^\mathrm{f}}
\newcommand{\Vf}{V^\mathrm{f}}
\newcommand{\Rv}{R^\mathrm{v}}
\newcommand{\Lv}{L^\mathrm{v}}
\newcommand{\Vv}{V^\mathrm{v}}

\newcommand{\Lhyd}{\mathscr{L}}
\newcommand{\myunit}[1]{\mbox{$\mathrm{#1}$}}

\newcommand{\lemmatwo}{%
  There exists a constant $\bigC_1$, independent of $\eps$, such that
  the solution $\Theta_\eps$ of \eqref{problem_weak_together2}
  (equivalently, $E_{1,\eps}$ and $E_{2,\eps}$ of
  \eqref{problem_weak_together1}) satisfies
  \begin{equation*}
    \norm{\Theta_\eps}^2_\Omega + \norm{\kappa_\eps\nabla\Theta_\eps}^2_{\Omega,\,t}
    =\Norm{E_{1,\eps}}^2_{\Omoneeps} + \Norm{\nabla
      \HE_{1,\eps}}^2_{\Omoneeps,\,t} + \Norm{E_{2,\eps}}^2_{\Omtwoeps}
    + \eps^2\Norm{\nabla \HE_{2,\eps}}^2_{\Omtwoeps,\,t} \leq \bigC_1  .
  \end{equation*}
}

\newcommand{\theoremone}{%
  Consider equation \eqref{ext_problem1} satisfying the conditions
  \eqref{ext_cond:abcde}.  Then there exists at least one solution of
  equation \eqref{ext_problem2}.
}

\newcommand{\theoremtwo}{%
  Equations~\eqref{s_problem_limit_nonlinear} have at most one solution
  given by 
  \begin{alignat*}{3}
    T_1 &\in \V^1(\Omega) + \Tout &=& 
    \left( L^2\left([0,\tend], \Hil^1_0(\Omega)\right) + \Tout\right)
    \cap \Hil^1 \left([0,\tend], L^2(\Omega)\right), \\    
    T_2 &\in \V^2(\Omega\times Y^2) + T_1 &=&
    \left(L^2 \left([0,\tend], L^2(\Omega,\Hil^1_0(Y^2))\right) + T_1\right) 
    \cap
    \Hil^1 \left([0,\tend], L^2(\Omega\times Y^2)\right), 
  \end{alignat*}
  where $T_1 = \omega(\HE_1)$ and $T_2 = \omega(\HE_2)$.
}

\begin{document}

\title{A two-scale Stefan problem arising in a model for tree sap exudation}
\shorttitle{A two-scale Stefan problem}
\shortauthorlist{I.\ Konrad, M.\ A.\ Peter \&\ J.\ M.\ Stockie}

\author{%
  \name{Isabell Konrad}
  \address{Department of Mathematics, Simon Fraser University,
    8888 University Drive, Burnaby, BC, V5A 1S6,
    Canada\email{isabella.konrad@gmail.com}} 
  \name{Malte A. Peter}
  \address{Institut f\"ur Mathematik,\ 
    University of Augsburg, Universit\"atsstra\ss e 14, 86159 Augsburg,
    Germany and Augsburg Centre for Innovative Technologies, University
    of Augsburg, 86135 Augsburg, Germany\email{malte.peter@math.uni-augsburg.de}} 
  \and
  \name{John M. Stockie}
  \address{Department of Mathematics, Simon Fraser University,
    8888 University Drive, Burnaby, BC, V5A 1S6, Canada\email{Corresponding author: stockie@math.sfu.ca}}}

\maketitle

\begin{abstract}%
  {The study of tree sap exudation, in which a (leafless) tree generates
    elevated stem pressure in response to repeated daily freeze--thaw
    cycles, gives rise to an interesting multiscale problem involving
    heat and multiphase liquid/gas transport.  The pressure generation
    mechanism is a cellular-level process that is governed by
    differential equations for sap transport through porous cell
    membranes, phase change, heat transport, and generation of osmotic
    pressure.  By assuming a periodic cellular structure based on an
    appropriate reference cell, we derive an homogenized heat equation
    governing the global temperature on the scale of the tree stem, with
    all the remaining physics relegated to equations defined on the
    reference cell.  We derive a corresponding strong formulation of the
    limit problem and use it to design an efficient numerical solution
    algorithm.  Numerical simulations are then performed to validate the
    results and draw conclusions regarding the phenomenon of sap
    exudation, which is of great importance in trees such as sugar maple
    and a few other related species.  The particular form of our
    homogenized temperature equation is obtained using periodic
    homogenization techniques with two-scale convergence, which we
    investigate theoretically in the context of a simpler two-phase
    Stefan-type problem corresponding to a periodic array of melting
    cylindrical ice bars with a constant thermal diffusion coefficient.
    For this reduced model, we prove results on existence, uniqueness
    and convergence of the two-scale limit solution in the weak form,
    clearly identifying the missing pieces required to extend the proofs
    to the fully nonlinear sap exudation model.  Numerical simulations
    of the reduced equations are then compared with results from the
    complete sap exudation model.}%
  {periodic homogenization; two-scale convergence; Stefan problem;
    multiphase flow; phase change.}\\
  MSC(2010):
  35B27, 
  35R37, 
  76T30, 
  80A22, 
  92C80. 
\end{abstract}


\section{Introduction}
\label{sec:intro}

This paper is motivated by the study of sap flow in sugar maple trees
that are subject to repeated cycles of thawing and freezing during the
sap harvest season in late winter~\citep{AHORN5}.  We seek insight into
the phenomenon of \emph{sap exudation}, which refers to the generation
of elevated sap pressure within the maple stem when the tree is in a
leafless state and no transpiration occurs to drive the sap flow.  Our
work is based on the model derived in~\citet{AHORN5} that
captures the physical processes at the microscale (i.e., at the level of
individual wood cells) and includes multiphase flow of ice/water/gas,
heat transport, porous flow through cell walls, and osmosis.  There is
an inherent repeating structure in sapwood at the cellular scale that
lends itself naturally to the use of homogenization ideas that we
exploited in~\citet{graf-ceseri-stockie-2015} to obtain a
multiscale model for the macroscale temperature that is coupled to a
corresponding system of equations governing the microscale cellular
processes.  Our main objective in this paper is to provide a more
rigorous theoretical justification for this multiscale model by working
through the details of the homogenization process and proving results
regarding existence, uniqueness and two-scale convergence.

Multiscale problems such as the one just described are characterized
by geometric, material or other features that exhibit variations on
widely differing spatial scales.  Many mathematical and numerical
methods have been developed to capture such scale separation as well
as the interactions between physical phenomena operating on disparate
scales \citep{engquist-etal-2005,TWOSCALE26}.  For problems having a
periodic microstructure, a mathematical technique that has proven to
be very effective is known as periodic
homogenization \citep{cioranescu}, and more specifically the method of
two-scale convergence \citep{TWOSCALE1, TWOSCALE2}, which has also been
extended to capture non-periodically evolving
microstructures \citep{CRM07,CRM07b,TWOSCALE19}. We are interested here
in applying two-scale convergence to analyze solutions of a
Stefan-type problem that governs the dynamics of the ice/water
interface within individual tree cells.  Locally, temperature obeys
the heat equation and is coupled with a Stefan condition that governs
solid--liquid phase transitions at the interface.  Many different
approaches have been developed to analyze such phase transitions,
which are well-described in~\citet{STEFAN5}.  With the
exception of a few studies of (single-phase) water and solute
transport in plant
tissues \citep{chavarriakrauser-ptashnyk-2010,chavarriakrauser-ptashnyk-2013},
periodic homogenization techniques have not been applied in the
context of heat or sap flow in trees.

The approach we employ in this paper has the advantage that it applies
homogenization techniques in a straightforward manner in order to obtain
an uncomplicated limit model, the simplicity of which ensures that
numerical simulations are relatively easy to perform.  In particular, we
define a reference cell $Y$ that is divided into two sub-regions: $Y^1$,
where the temperature diffuses rapidly; and $Y^2$, on which we define a
second temperature field that diffuses slowly.  Refer
to \citet{PU5} and \citet{TWOSCALE16} for similar homogenization approaches
involving slow and fast transport. One particular challenge arising in
the study of Stefan problems is that the diffusion coefficient depends
on the underlying phases, so that heat diffuses differently in water or
ice.  Consequently, the diffusion coefficient depends on temperature (or
equivalently on enthalpy) so that the governing differential equation is
only quasi-linear.

Rather than attempting to analyze the sap exudation problem in its full
complexity, we find it more convenient to develop our homogenization
results in the context of a simpler ``reduced model'' defined on a
similarly fine-structured domain wherein the cell-level processes are
governed by a Stefan problem that involves only heat transport and
ice/water phase change.  In particular, we consider a domain consisting
of a periodic array of cylindrical ice inclusions immersed in water.  To
handle the multiplicity of the ice bars, we apply the technique of
periodic homogenization with two-scale convergence established in
\citet{TWOSCALE1} and \citet{TWOSCALE2}.  Several authors have
previously applied homogenization to Stefan problems, such as
\citet{bossavit-damlamian-1981}, \citet{STEFAN1} and \citet{STEFAN2},
where the phase change boundary is handled by separately homogenizing an
auxiliary problem. In \citet{STEFAN3} on the other hand, an additional
function $\theta$ is introduced for an aggregate state that diffuses on
a slow time scale and with which all microscopic phase changes are
properly captured.  When we show existence for the heat equation with
phase transitions, we deduce a general existence result for quasi-linear
parabolic differential equations having a non-monotone nonlinearity in
the diffusion operator, which is of general interest in the context of
heat transport and Stefan problems (even in a single-scale setting).

This paper is organized as follows.  We begin in
Section~\ref{sec:sap-model} by providing background material on the
physics of maple sap exudation, along with a description of the
governing equations at the cellular level.  A reduced model involving
only melting of ice is introduced in Section~\ref{sec:problem1} for the
purposes of more easily deriving the two-scale convergence results.  The
main analytical results on existence, a~priori estimates, two-scale
convergence and uniqueness are presented in Section~\ref{sec:analysis},
and detailed proofs of the key results are relegated to the Appendix.
Following that, we state in Section~\ref{sec:strong-form} the strong form
of the limit problem for the reduced model, which in turn suggests a
corresponding strong form of the original sap exudation model in
Section~\ref{sec:sap-limit-problem}.  These limit problems lead
naturally to a multiscale numerical algorithm that is described in
Section~\ref{sec:simulations}, after which numerical simulations of both
problems are presented and compared.

\section{Mathematical model for sap exudation}
\label{sec:sap-model}

Before presenting the details of the mathematical model, it is necessary
to introduce some background material on the phenomenon of sap
exudation.  Sugar maple trees (along with a few related species such as
red or black maple, black walnut, and birch) have a unique ability
compared to other deciduous tree species in that they exude large
quantities of sap during the winter when they are in a leafless state.
Sap exudation originates from an elevated pressure in the tree stem that
is generated over a period of several days during which the air
temperature oscillates above and below the freezing point.  The ability
of maple to exude sap has intrigued tree physiologists for over a
century, and various physical and biological processes have been
proposed to explain this behaviour \citep{johnson-tyree-dixon-1987,
  milburn-kallarackal-1991, tyree-1995}.  Until recently, a significant
degree of controversy existed over the root causes of sap exudation, and
the most plausible and widely-accepted explanation has been a
freeze--thaw hypothesis proposed by \citet{AHORN2}.  This hypothesis
forms the basis of the mathematical model for the cellular processes
underlying exudation during a thawing event that was developed by 
\citet{AHORN5}, which was subsequently extended to capture a
complete freeze--thaw cycle by \citet{graf-ceseri-stockie-2015}.

\subsection{Background: Tree physiology and the Milburn--O'Malley process}
\label{sec:milburn-omalley}

\begin{figure}[tbhp]
  \centering\footnotesize
  \begin{tabular}{ccc}
    (a) Sapwood microstructure &&
    (b) A single fiber/vessel pair\\
    {\setlength{\unitlength}{0.33\textwidth}
       \begin{picture}(1,1.40)
         \put(0,0){\includegraphics[width=0.28\textwidth]{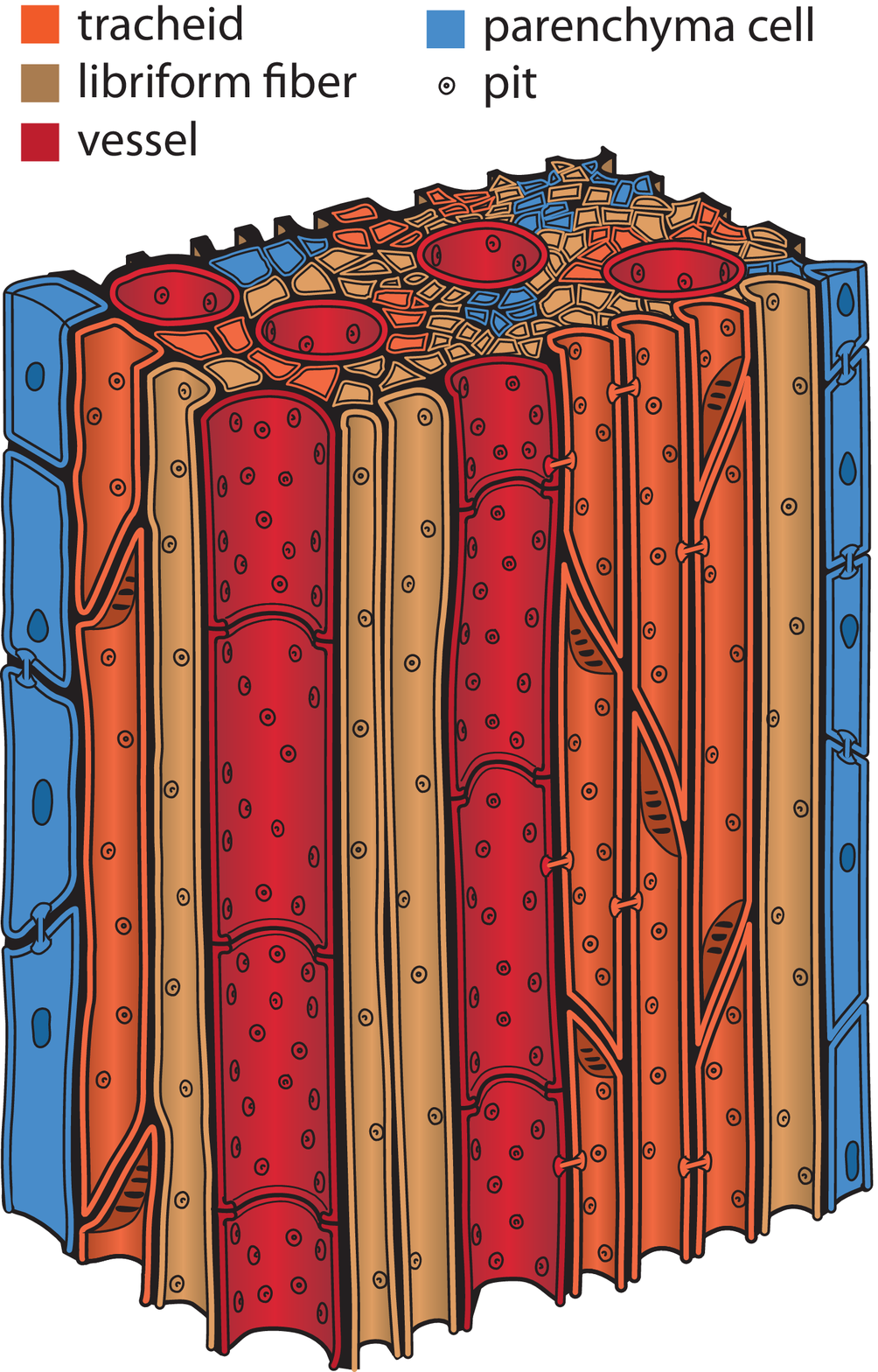}}
         \thicklines
         \linethickness{0.5mm}
         \put(0.65,0.02){\line(1,0){0.25}}
         \put(0.65,0.0){\line(0,1){0.04}}
         \put(0.90,0.0){\line(0,1){0.04}}
         \put(0.70,0.05){\footnotesize\sffamily 100\,$\mu$m} 
       \end{picture}
    }
    & &
    \includegraphics[width=0.25\textwidth]{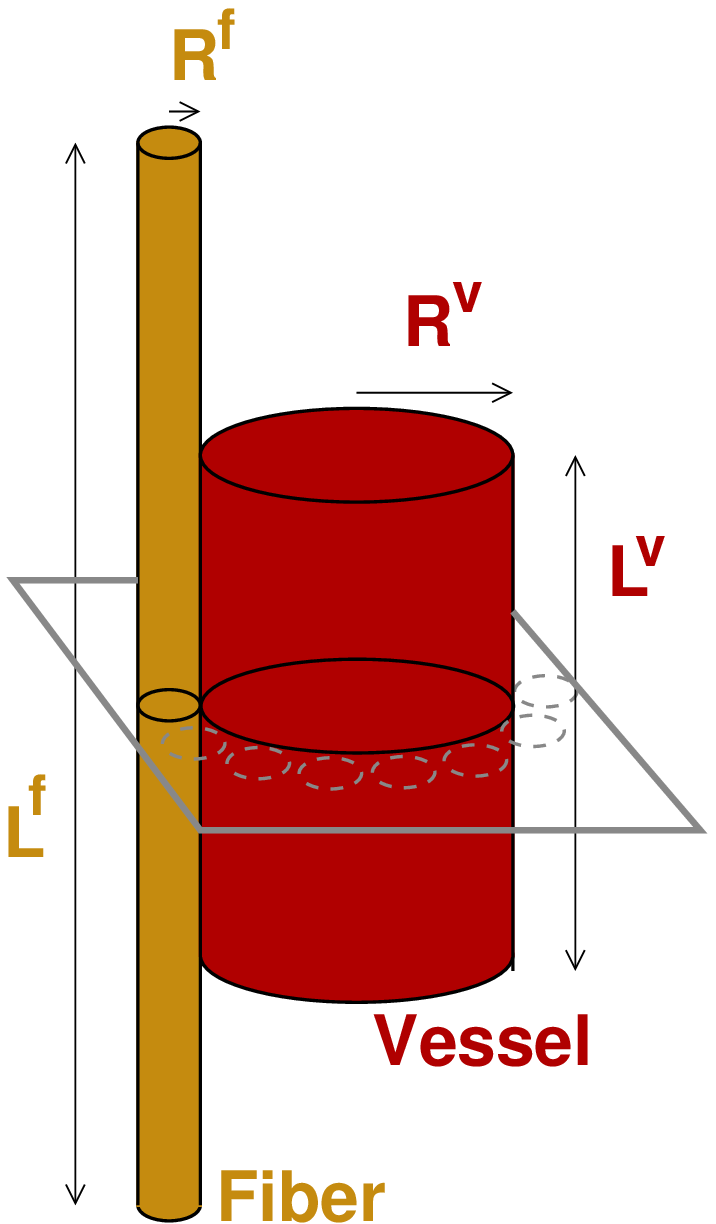}
  \end{tabular}
  \caption{(a) A cut-away view of the sapwood (or xylem) in hardwood
    trees such as sugar maple, depicting the repeating microstructure of
    vessels surrounded by fibers (the other cell types indicated here
    are ignored in our model).  (b) A vessel surrounded by $\Nf$ fibers,
    all depicted as circular cylinders (for simplicity, only one fiber
    is shown).  Typical dimensions of the fiber are length
    $L^\mathrm{f}=1.0\times 10^{-3}\;\mathrm{m}$ and radius
    $R^\mathrm{f}=3.5\times 10^{-6}\;\mathrm{m}$, whereas the vessel has
    $L^\mathrm{v}=5.0\times 10^{-4}\;\mathrm{m}$ and
    $R^\mathrm{v}=2.0\times 10^{-5}~\mathrm{m}$.  The 2D model reference
    cell introduced in what follows is based on a horizontal
    cross-section through the middle of the fiber and vessel.}
  \label{fig:xylem}
\end{figure}

The Milburn--O'Malley process depends crucially on the distinctive
microstructure of sapwood (or xylem) in sugar maple trees (\emph{Acer
  saccharum}).  Wood in most deciduous tree species consists of roughly
cylindrical cells that are on the order of $1\;\mathrm{mm}$ in length.
These cells can be classified into two main types: \emph{vessels} having
an average radius of $20\;$\textmu$\mathrm{m}$, which are surrounded by
the much more numerous \emph{(libriform) fibers} with a radius of
approximately 3--4\;\textmu$\mathrm{ m}$.  The repeating structure of
vessels and fibers is illustrated in Figure~\ref{fig:xylem}a.  The
vessels have a significantly larger diameter and therefore comprise the
main route for sap transport between roots to leaves during the growing
season, whereas the fibers are understood to play a largely passive and
more structural role.  Under normal conditions the vessels are filled
with sap, which is composed primarily of water but also contains as much
as 2--5\%\ sugar by weight in species like \emph{Acer}.  On the other
hand, the fibers are thought to be primarily filled with gas (i.e.,
air).  We note that experiments exhibit small but measurable amounts of
gas also being present within the vessel sap, either as bubbles or in
dissolved form.

\begin{figure}[tbhp]
  \centering\footnotesize
  \includegraphics[width=0.75\textwidth]{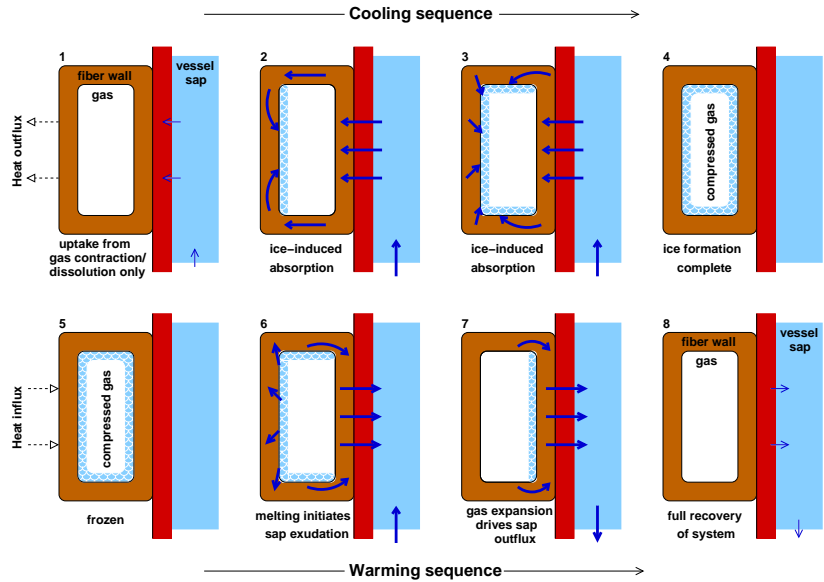}
  \caption{Illustration of the Milburn--O'Malley process in a single
    fiber/vessel pair (adapted from~\citet[Fig.~7]{AHORN2}).
    We focus on the ``warming sequence'' in the bottom row, numbered
    5--8.  The fiber is the large rectangular structure on the left of
    each sub-image, and the vessel is the vertical channel on the right
    (not drawn to scale).}
  \label{fig:milburn-omalley}
\end{figure}

Milburn and O'Malley hypothesized that during late winter when daily
high temperatures peak above the freezing point, and just as evening
temperatures begin to drop below zero, sap is drawn through tiny pores
in the fiber/vessel walls by capillary and adsorption forces into the
gas-filled fibers where it forms ice crystals on the inner surface of
the fiber wall (top ``cooling sequence'' in
Figure~\ref{fig:milburn-omalley}).  As temperatures drop further, the
ice layer grows and the gas trapped inside the fiber is compressed,
forming a pressure reservoir of sorts.  When temperatures rise above
freezing again the next day, the process reverses, with the ice layer
melting and the pressurized gas driving liquid melt-water back into the
vessel where it then (re-)pressurizes the vessel compartment (bottom
``warming sequence'' in Figure~\ref{fig:milburn-omalley}).  Milburn and
O'Malley also stressed the importance of osmotic pressure in terms of
maintaining the high stem pressures actually observed in sugar maple
trees. This essential role of osmosis has since been verified
experimentally by \citet{AHORN4} who confirmed the
existence of osmotic pressure arising from a selectively permeable
membrane within the fiber/vessel wall.  They showed that the cell wall
permits water to pass but prevents larger sugar molecules contained in
the vessel sap from entering the fiber, thereby introducing a
significant osmotic pressure difference between the sugar-rich vessel
sap and the pure water contained in the fiber.

There are two additional physical effects not explicitly addressed by
\citet{AHORN2} that are essential in order to
obtain physically-consistent results for the sap thawing process.  First
of all, \citet{AHORN5} demonstrated the necessity for
including gas bubbles suspended within the vessel sap that permit an
exchange of pressure between the vessel and fiber compartments, which
would otherwise not be possible owing to the incompressibility of water.
Secondly, despite the pervading belief that there is no significant root
pressure in maple during winter \citep{kramer-boyer-1995,wilmot-2011},
we found to the contrary that including uptake of root water during the
freezing process is absolutely essential in order that pressure can
accumulate over multiple freeze--thaw
cycles \citep{graf-ceseri-stockie-2015}.  Indeed, the need to include
root pressure is confirmed by recent experiments \citep{brown-2015} that
demonstrate the existence of root pressure in maple trees during the sap
harvest season.

\subsection{Microscale model for cell-level processes}
\label{sec:cell-level}

The modified Milburn--O'Malley description just presented (with the
exception of root pressure) was employed by \citet{AHORN5} and
\citet{graf-ceseri-stockie-2015} to derive a mathematical model for
cell-level processes governing sap exudation during a thawing cycle.  In
this study, we study the same problem, including the effect of the gas
phase in both cell chambers (fiber and vessel), but we will assume for
the sake of simplicity that the effects of gas dissolution and
nucleation are negligible.  This is the primary difference between our
microscale model and that in \citet{AHORN5} and
\citet{graf-ceseri-stockie-2015}, on which it is based.  Neglecting root
pressure is a reasonable simplification because we are only interested
here in studying a single thawing event and not capturing repeated
freeze--thaw cycles.

\begin{figure}[tbhp]
  \centering\footnotesize
  \begin{tabular}{ccc}
    (a) Reference cell, with ice layer\mbox{\hspace*{2cm}}
    &&
    (b) Reference cell, completely melted\\
    \includegraphics[height=0.25\textheight]{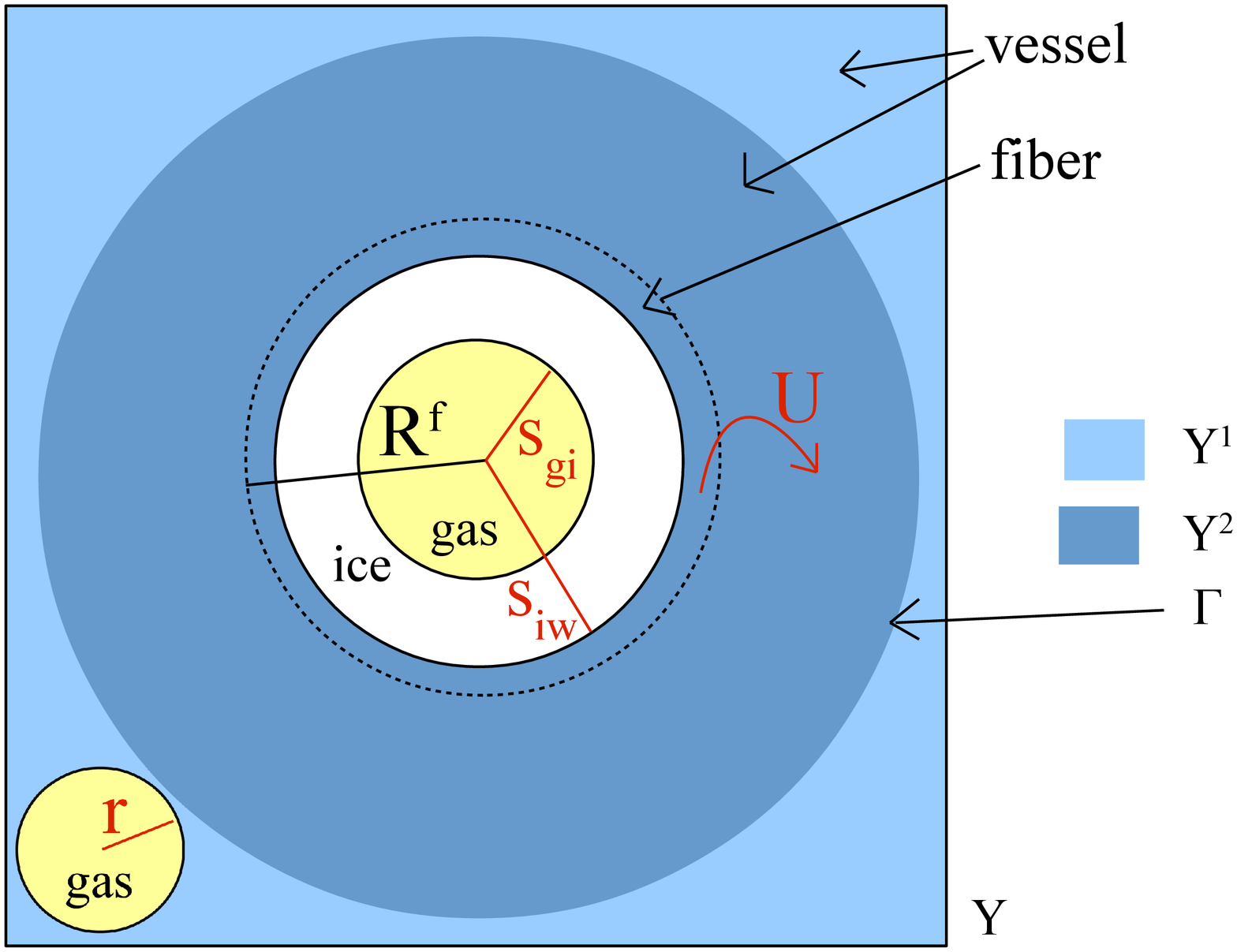}
    && 
    \includegraphics[height=0.25\textheight]{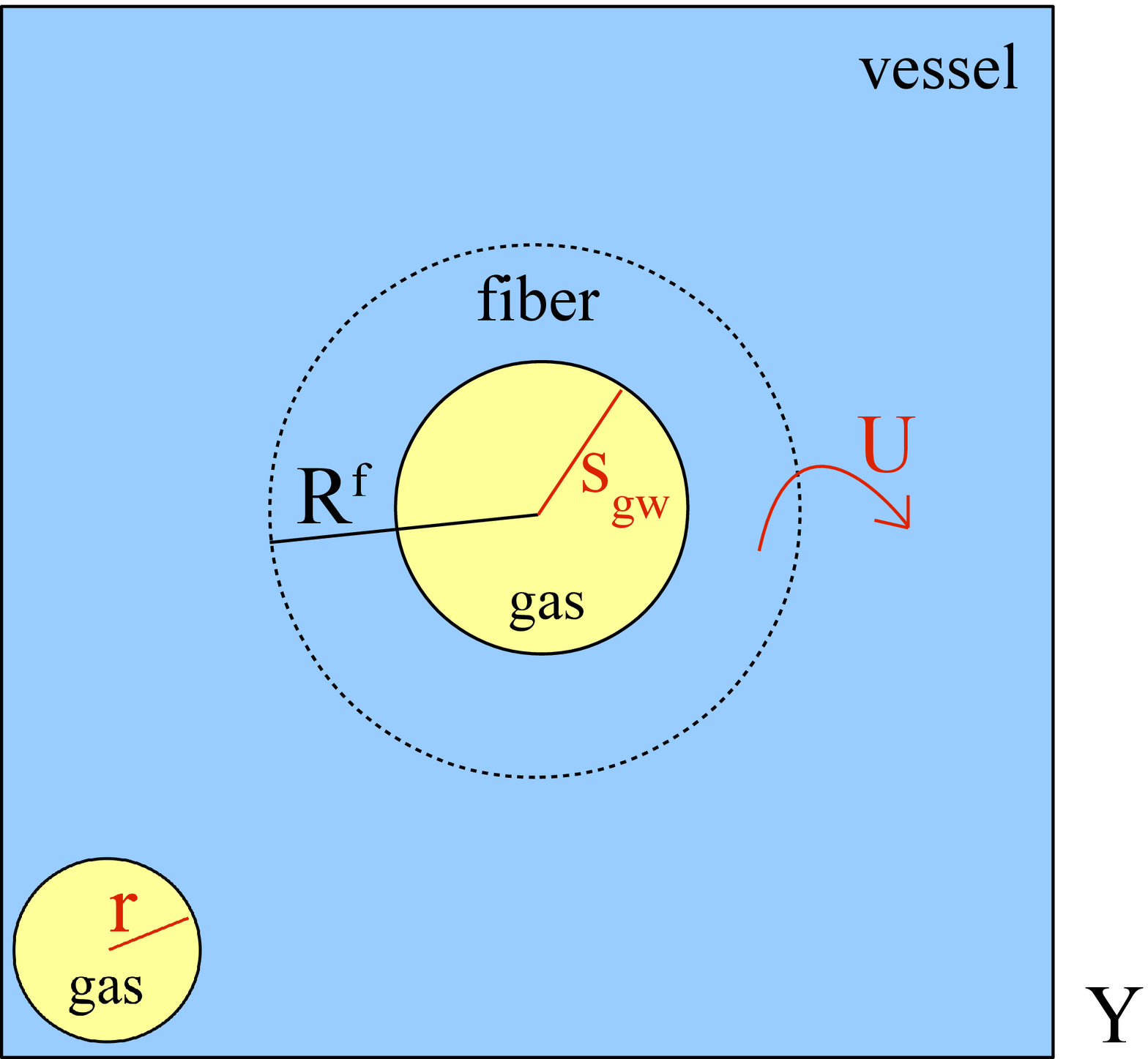}
  \end{tabular}
  \caption{Geometry of the reference cell, $Y$.  (a)~For the sap
    exudation problem, the vessel compartment contains a gas bubble with
    radius $r$, while the circular fiber (radius $\Rf$) contains a gas
    bubble (radius $s_\mathrm{gi}$), surrounded by an ice layer (with
    thickness $s_\mathrm{iw}-s_\mathrm{gi}$), and finally a layer of
    melt-water (with thickness $\Rf-s_\mathrm{iw}$).  The porous wall
    between fiber and vessel is denoted by a dotted line.  As ice melts,
    the melt-water is forced out by gas pressure through the porous
    fiber wall into the surrounding vessel compartment.  The total
    volume of melt-water transferred from fiber to vessel is denoted
    $U$.  An artificial boundary $\Gamma$ is introduced in the
    homogenization process to differentiate between a region $Y^1$
    (outside $\Gamma$) on which thermal diffusion is fast, and $Y^2$
    (inside $\Gamma$) on which diffusion is slow.  (b) After the ice has
    completely melted there remains only a gas/water interface,
    $s_\mathrm{gw}$, and a single temperature field can be used to
    describe the entire reference cell domain $Y \equiv Y^1$.}
  \label{fig:sap-refcells}
\end{figure}

With the above assumptions in mind, we approximate the
sapwood as a periodic array of square reference cells $Y$ pictured in
Figure~\ref{fig:sap-refcells}a.  Each reference cell contains a circular
fiber of radius $\Rf$ located at the centre, surrounded by a
vessel compartment that makes up the remainder of the cell.  Because the
vessels have considerably larger diameter, we assume that
on the scale of a fiber the cylindrical geometry of the vessel can be
neglected as long as we ensure that appropriate conservation principles
(for mass and energy) are maintained within the vessel.  This choice of
reference cell is obviously a mathematical idealization that may
influence fine details of vessel transport on the microscale but
ultimately has minimal impact on the homogenized solution.

The fiber compartment is sub-divided into nested annular regions
containing gas, ice and liquid, and the outer radii of the phase
interfaces are denoted $s_\mathrm{gi}$ (for gas/ice) and $s_\mathrm{iw}$ (for
ice/water).  The vessel contains a circular gas bubble of radius $r$
which has no specified location but rather is included simply to track
the amount of gas for mass-conservation purposes.  One additional
variable $U$ is introduced to measure the total volume of water
transferred from fiber to vessel.  The region lying outside the fiber
and inside the boundary of $Y$ represents the sugary sap-filled vessel.
Note that during a thawing cycle, we are only concerned with a vessel
containing liquid sap (no ice) because of the effect of freezing point
depression, which ensures that any given vessel thaws before the
adjacent fiber(s).  This reference cell geometry should be contrasted
with that depicted in~\citet[Fig.~3.1]{AHORN5}.

For the moment, we will consider the four solution variables
$s_\mathrm{iw}$, $s_\mathrm{gi}$, $r$, $U$ as depending on time $t$
only, with an additional dependence of temperature on the microscale
spatial variable; however, beginning in Section~\ref{sec:problem1} when
we derive macroscale equations for the homogenized problem, these
variables will also depend on the global spatial variable $x$ that
denotes the location of the reference cell within the tree stem.  Within
a reference cell, the dynamics for $s_\mathrm{iw}(t)$,
$s_\mathrm{gi}(t)$, $r(t)$ and $U(t)$ are governed by four differential
equations whose derivation can be found in~\citet{AHORN5}.  The
first is the Stefan condition for the ice/water interface in the fiber
\begin{subequations}\label{cell_problem_maple}
  \begin{gather}
    \partial_t s_\mathrm{iw} = 
    -\frac{D(\HE_2)}{\HE_\mathrm{w}-\HE_\mathrm{i}}\nabla_y T_2 \cdot
    \normvec + \frac{\partial_t U}{2\pi s_\mathrm{iw} \Lf} , 
    \label{cell_problem_maple_siw}
  \end{gather}
  where $T_2(y,t)$ denotes the microscale temperature variable that
  depends on both time and the local spatial variable $y\in Y^2$
  (which needs to be distinguished from the macroscale temperature
  variable $T_1(x,t)$ introduced later) and $\nabla_y T_2 \cdot
  \normvec$ is the normal derivative at the
  interface.  
  Here, the enthalpies of water and ice ($\HE_\mathrm{w}$ and
  $\HE_\mathrm{i}$, resp.) are evaluated at the freezing point, $T=\Tcrit$;
  consequently, the difference $\HE_\mathrm{w}-\HE_\mathrm{i}$
  represents the latent heat of fusion.  We describe heat transport
  using a mixed temperature--enthalpy formulation, in which the thermal
  diffusion coefficient $D(\HE_2)$ is written as a function of 
  enthalpy $\HE_2(y,t)$.  Following~\citet{STEFAN5}, we take
  $D$ to have the piecewise affine linear form
  \begin{align}
    D(\HE) = \left\{\begin{array}{cl}
        \frac{k_\mathrm{i}}{\rho_\mathrm{i}}, & \qquad \text{if } \HE<\HE_\mathrm{i},\\
        \frac{k_\mathrm{i}}{\rho_\mathrm{i}} + \frac{\HE-\HE_\mathrm{i}}{\HE_\mathrm{w}-\HE_\mathrm{i}}
        \left(\frac{k_\mathrm{w}}{\rho_\mathrm{w}} - \frac{k_\mathrm{i}}{\rho_\mathrm{i}}\right), 
        & \qquad \text{if } \HE_\mathrm{i}\leq \HE < \HE_\mathrm{w}, \\ 
        \frac{k_\mathrm{w}}{\rho_\mathrm{w}}, & \qquad \text{if } \HE_\mathrm{w} < \HE,
      \end{array}\right.
    \label{DH}
  \end{align}
  where $\rho_\mathrm{w}$, $\rho_\mathrm{i}$ are the densities of water
  and ice respectively, and $k_\mathrm{w}$, $k_\mathrm{i}$ are the
  thermal conductivities.  Note that $D$ in this temperature--enthalpy
  formulation has units of \myunit{W\, m^2/kg\,\degK} and is referred
  to as a thermal diffusion coefficient, to distinguish it from the more
  usual ``thermal diffusivity'' (which is defined as the ratio $k/\rho
  c$ and has units of \myunit{m^2/s}).  The governing equations for
  $T_2$ and $\HE_2$ are discussed later in
  Sections~\ref{sec:strong-form}--\ref{sec:sap-limit-problem} as a
  result of the two-scale convergence analysis and are the solutions of
  the system~(\ref{limit_problem_strong}a--e).  Note that the final term
  in the Stefan condition~\eqref{cell_problem_maple_siw} was neglected
  in~\citet{AHORN5} and serves to capture the effect on the phase
  interface of fiber--water volume changes due to porous flow through
  the fiber/vessel wall.
  
  The next two differential equations embody conservation of mass in
  the fiber
  \begin{gather}
    \partial_t s_\mathrm{gi} = 
    -\frac{(\rho_\mathrm{w}-\rho_\mathrm{i}) s_\mathrm{iw}\partial_t s_\mathrm{iw}}{s_\mathrm{gi} \rho_\mathrm{i}} +
    \frac{\rho_\mathrm{w} \partial_t U}{2\pi s_\mathrm{gi} \rho_\mathrm{i} \Lf} , 
    \label{cell_problem_maple_sgi}
  \end{gather}
  and the vessel 
  \begin{gather}
    \partial_t r = 
    -\frac{\Nf \partial_t U}{2\pi r \Lv}.
    \label{cell_problem_maple_r}
  \end{gather}
  Note that within the sapwood there are many more fibers than vessels
  (as depicted in Figure~\ref{fig:xylem}a), so that the effect of
  fiber--vessel flux terms should be increased to account for the
  multiplicity of fibers.  With this in mind, we have multiplied
  appropriate fluxes by the parameter $\Nf$ in
  \eqref{cell_problem_maple_r} that represents the average number of
  fibers per vessel and has a typical value of $\Nf=16$.  The final
  differential equation describes water transport through the porous
  fiber/vessel wall in response to both hydraulic and osmotic pressure
  \begin{gather}
    \partial_t U = 
    -\frac{\Lhyd A}{\Nf}(p_\mathrm{w}^\mathrm{v} -
    p_\mathrm{w}^\mathrm{f} - \Rgas C_\mathrm{s} T_1). 
    \label{cell_problem_maple_U}
  \end{gather}
\end{subequations}
Here, we denote the pressure variable by $p$, where superscripts
$\mathrm{f}$/$\mathrm{v}$ refer to fiber/vessel and subscript $\mathrm{w}$
denotes the liquid water phase.  The constant parameter $\Lhyd$ is the
fiber/vessel wall conductivity, $A$ is the wall surface area,
$C_\mathrm{s}$ is the vessel sugar concentration, and $\Rgas$ is the
universal gas constant.  Note that because $U$ is defined inside
$Y^2$, we should strictly be using the microscale temperature $T_2$ in
the osmotic term in \eqref{cell_problem_maple_U}, but this would lead
to a significant complication in any numerical algorithm due an
additional nonlinear coupling between scales.  Therefore, we have used
$T_1$ instead, which is a reasonable approximation because temperature
variations throughout the reference cell are small.

Several intermediate variables have been introduced into the above
equations. They are determined by the following algebraic relations:
\begin{subequations}\label{algebraic}
  \begin{alignat}{3}
    p_\mathrm{w}^\mathrm{f} &= p_\mathrm{g}^\mathrm{f}(0) \left(\frac{s_\mathrm{gi}(0)}{s_\mathrm{gi}}\right)^2 -
    \frac{2\sigma}{s_\mathrm{gi}},
    && \quad\text{(Young--Laplace equation for fiber)}
    \label{algebraic_pwf}
    \\
    p_\mathrm{w}^\mathrm{v} &= p_\mathrm{g}^\mathrm{v} - \frac{2\sigma}{r},
    && \quad\text{(Young--Laplace equation for vessel)}
    \label{algebraic_pwv}
    \\
    p_\mathrm{g}^\mathrm{v} &= \frac{\rho_\mathrm{g}^\mathrm{v} \Rgas T_1}{M_\mathrm{g}},
    && \quad\text{(ideal gas law for vessel)}
    \label{algebraic_pgv}
    \\
    \rho_\mathrm{g}^\mathrm{v} &= \rho_\mathrm{g}^\mathrm{v}(0) \left(\frac{r(0)}{r}\right)^2.
    \qquad\qquad\qquad
    && \quad\text{(vessel gas density)}
    \label{algebraic_rhogv}
  \end{alignat}
\end{subequations}
All constant parameters appearing in the above equations are listed in
Table~\ref{tab:params2} along with typical values.
\begin{table}[tbhp]
  \centering
  \caption{Constant parameter values appearing in the
    sap exudation model (taken from~\citet{AHORN5}).}  
  \label{tab:params2}
  \begin{tabular}{clcc}\hline
    {Symbol} & {Description} & {Value} & {Units}\\\hline
    \multicolumn{4}{l}{Geometric parameters:}\\
    $\delta$   & Side length of reference cell  & $4.33\times 10^{-5}$ & \myunit{m}\\
    $\gamma$   & $=\Rf+W$, Radius of $\Gamma$ & $7.88\times 10^{-6}$ & \myunit{m}\\
    $\Rf$      & Fiber radius            & $3.5\times 10^{-6}$  & \myunit{m}\\
    $\Rv$      & Vessel radius           & $2.0\times 10^{-5}$  & \myunit{m}\\
    $\Lf$      & Fiber length            & $1.0\times 10^{-3}$  & \myunit{m}\\  
    $\Lv$      & Vessel length           & $5.0\times 10^{-4}$  & \myunit{m}\\
    $\Vf$      & Fiber volume $=\pi({\Rf})^2 \Lf$    & $3.85\times 10^{-14}$ & \myunit{m^3}\\
    $\Vv$      & Vessel volume $=\pi({\Rv})^2 \Lv$   & $6.28\times 10^{-13}$ & \myunit{m^3}\\
    $A$        & Area of fiber/vessel wall $=2\pi\Rf\Lf$& $2.20\times 10^{-8}$  & \myunit{m^2}\\
    $W$        & Thickness of fiber/vessel wall & $4.38\times 10^{-6}$ & \myunit{m}\\
    $\Nf$      & Number of fibers per vessel& 16                & --\\
    $\Rtree$   & Tree stem radius           & 0.25              & \myunit{m} \\
    \hline
    \multicolumn{4}{l}{Thermal parameters:}\\
    $c_\mathrm{w}$    & Specific heat of water        & $4180$      & \myunit{J/kg \, \degK} \\ 
    $c_\mathrm{i}$    & Specific heat of ice          & $2100$      & \myunit{J/kg \, \degK} \\
    $\HE_\mathrm{w}$  & Enthalpy of water at $\Tcrit$ & $9.07\times 10^5$ & \myunit{J/kg} \\ 
    $\HE_\mathrm{i}$  & Enthalpy of ice at $\Tcrit$   & $5.74\times 10^5$ & \myunit{J/kg} \\ 
    $k_\mathrm{w}$    & Thermal conductivity of water & $0.556$     & \myunit{W/m\, \degK} \\
    $k_\mathrm{i}$    & Thermal conductivity of ice   & $2.22$      & \myunit{W/m\, \degK} \\
    $\rho_\mathrm{w}$ & Density of water              & $1000$      & \myunit{kg/m^3} \\
    $\rho_\mathrm{i}$ & Density of ice                & $917$       & \myunit{kg/m^3} \\
    $\Tcrit$   & Freezing temperature for water   & 273.15 & \myunit{\degK}\\
    $\Tout$    & Ambient temperature $=\Tcrit+10$ & 283.15 & \myunit{\degK}\\
    $\alpha$ & Heat transfer coefficient &  $10$ & \myunit{W/m^2\degK}\\
    \hline
    \multicolumn{4}{l}{Other parameters:}\\
    $M_\mathrm{g}$      & Molar mass of air           & 0.029           & \myunit{kg/mol}\\
    $\Rgas$    & Universal gas constant      & 8.314           & \myunit{J/mol\, \degK}\\
    $\sigma$   & Gas/liquid surface tension  & 0.076           & \myunit{kg/s^2}\\
    $C_\mathrm{s}$      & Vessel sugar concentration (2\%) & 58.4  & \myunit{mol/m^3}\\
    $\Lhyd$    & Hydraulic conductivity of fiber/vessel wall & $5.54\times 10^{-13}$& \myunit{m^2\,s/kg}\\
    \hline
  \end{tabular}
\end{table}

There is one special case to consider, namely when a fiber initially
containing an ice layer is above the freezing point for long enough time
that the ice melts completely.  In the moment the ice layer disappears,
the reference cell geometry appears as in Figure~\ref{fig:sap-refcells}b
and the cell-level equations must be modified as follows.  First of all,
$D(\HE_2)$ must change to account for the fact that there are two
possible values of thermal diffusivity, one in the region containing the
gas and another in the liquid.  Furthermore, the gas/ice and ice/water
interfaces merge so that Eq.~\eqref{cell_problem_maple_siw} drops out
and we identify a new fiber gas/water interface as $s_\mathrm{gw} \defeq
s_\mathrm{iw} \equiv s_\mathrm{gi}$.  This leads to the following
simplified version of~\eqref{cell_problem_maple_sgi}
\begin{gather*}
  \partial_t s_\mathrm{gw} = \frac{\partial_t U}{2\pi s_\mathrm{gw} \Lf} ,  
\end{gather*}
but otherwise the microscale equations
\eqref{cell_problem_maple}--\eqref{algebraic} remain the same.

The equations for the temperature and enthalpy variables appearing in
the microscale model above are derived in the next section in the
context of a simpler problem involving only melting ice.  Despite the
fact that this reduced model involves only a single microscale variable
for the dynamics of the ice--water interface $s_\mathrm{iw}$ (in
addition to the temperature), the equations for temperature and enthalpy
remain the same, and we will show that the microscale model above is
completed by Eqs.~\eqref{limit_problem_strong1a}--\eqref{limit_problem_strong2c}.

\section{Reduced model: Melting ice bars}
\label{sec:problem1}

We now shift our attention to the macroscale problem, which captures the
dynamics of thawing sap within a cylindrical tree stem having a circular
cross-section $\Omega$.  There is a clear separation
of scales in that the tree has radius on the order of tens of
centimetres whereas the cell-level processes occur over distances on the
order of microns.  Let $x\in\Omega$ represent the macroscale spatial
variable and $y\in Y$ the microscale variable on the reference cell.
Then, our main aim in this section is to determine equations for the
temperature and enthalpy variables not only in the reference cell, 
$T_2(y,t)$ and $\HE_2(y,t)$, but also on the macroscale,
$T_1(x,t)$ and $\HE_1(x,t)$.

\begin{figure}[tbhp]
  \centering\footnotesize
  \begin{tabular}{ccc}
    (a) Periodically-tiled domain $\Omega$ & \qquad\qquad & 
    (b) Reference cell $Y$ for reduced model\qquad\quad \\
    \includegraphics[width=0.40\textwidth]{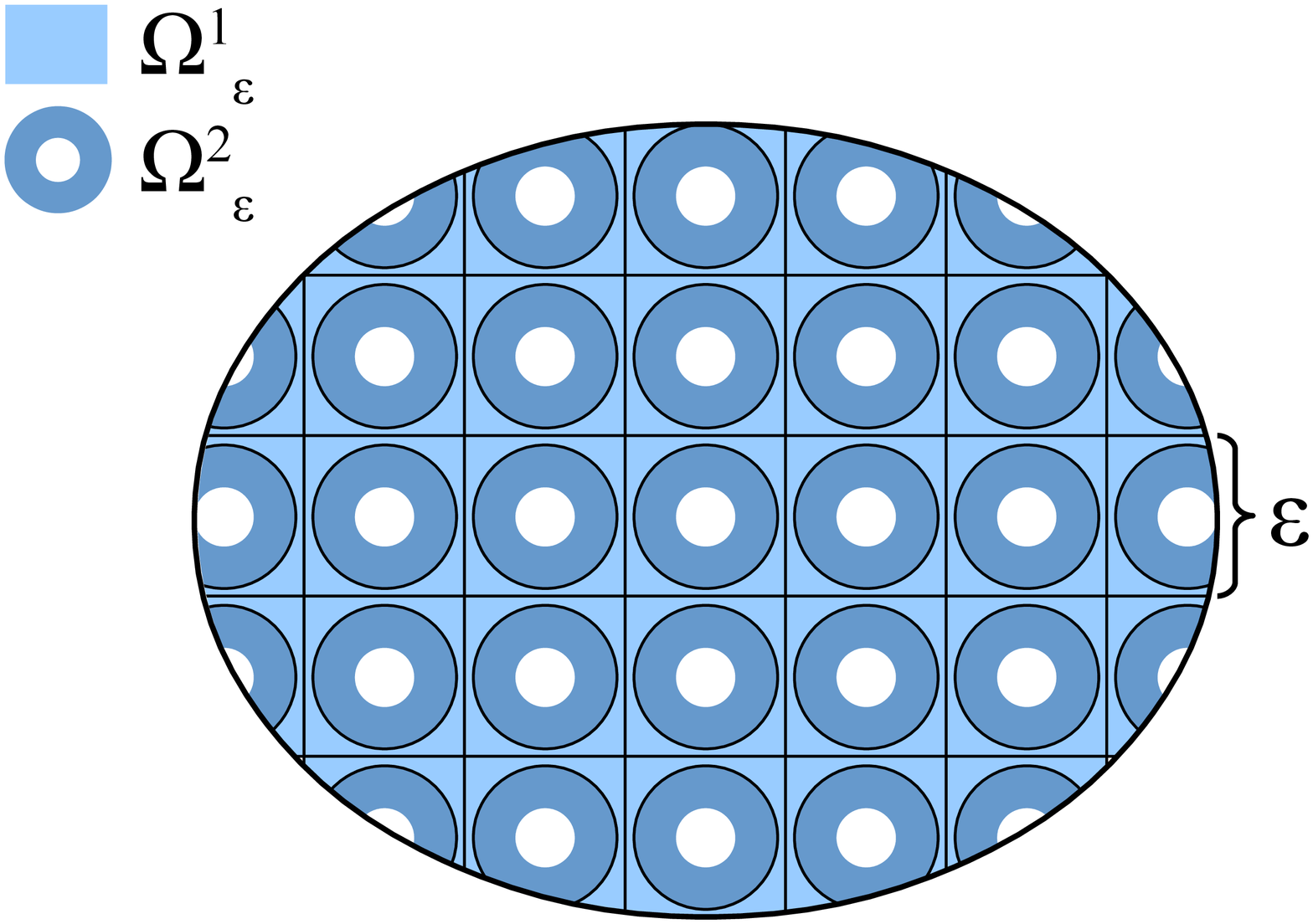}
    & & 
    \raisebox{0.7cm}{\includegraphics[width=0.280\textwidth]{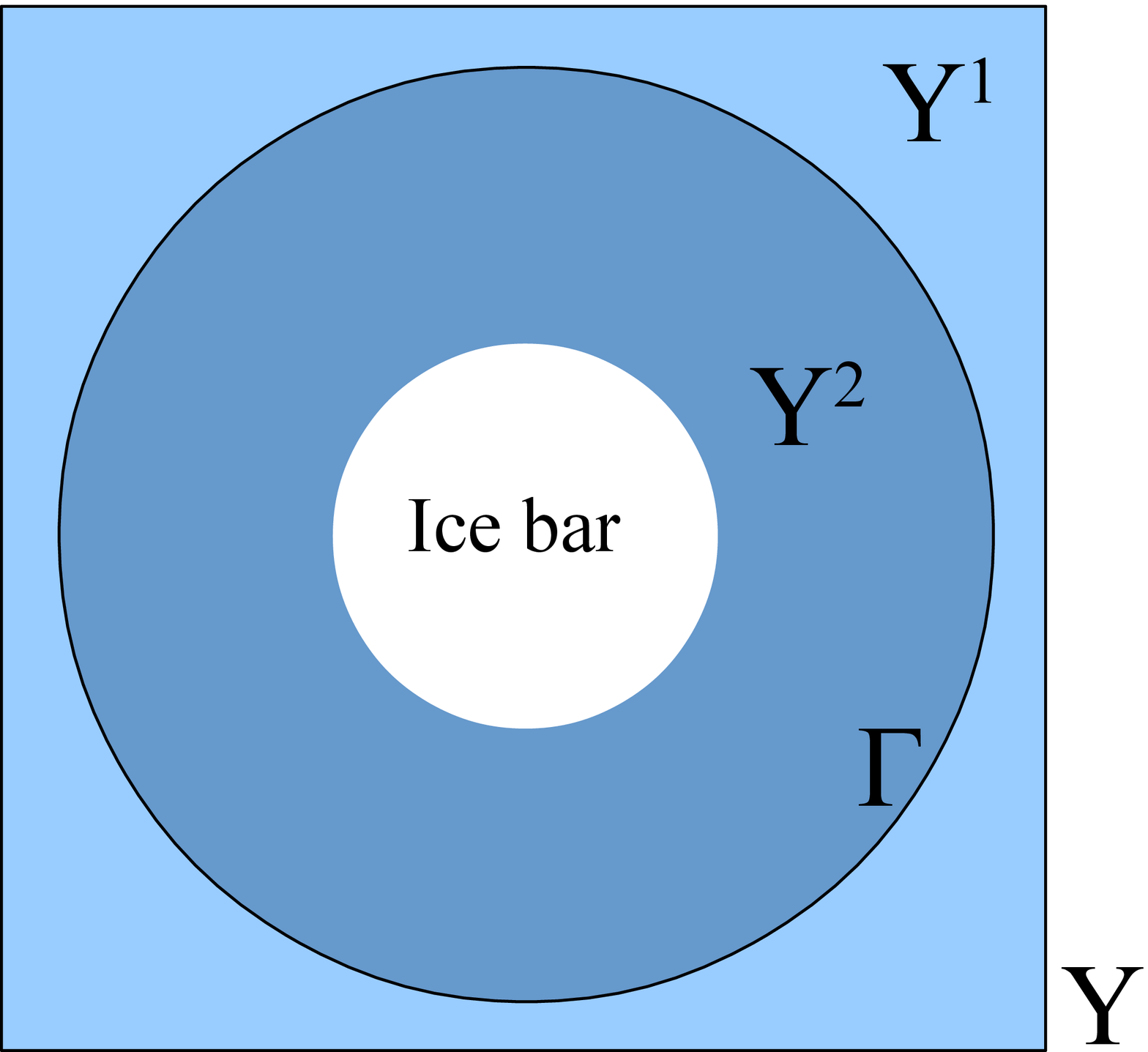}}
  \end{tabular}
  \caption{Periodic microstructure of the reduced model for melting ice
    bars immersed in water. (a)~The tree stem cross-section $\Omega$ is
    tiled periodically with copies of the reference cell $Y$, each of
    which is scaled to have side length $\eps$.  The homogenization
    process then takes the limit as $\eps\rightarrow 0$.  (b)~The
    reference cell $Y$ for the reduced model, illustrating the
    decomposition into fast ($Y^1$ and $\Omoneeps$) and slow ($Y^2$
    and $\Omtwoeps$) diffusing regions, with $Y=Y^1\cup Y^2 \cup \Gamma$.}
  \label{fig:ice-refcell}
\end{figure}

The derivation of these equations may be simplified significantly by
considering a reduced problem that involves only ice/water phase change
and leaves out all other physical processes (porous flow, gas bubbles,
surface tension, etc.).  To this end, we consider a periodic array of
melting ``ice bars'' as pictured in Figure~\ref{fig:ice-refcell}a,
situated inside a slightly more general domain $\Omega\subset \R^d$
having Lipschitz boundary that contains both water and ice in the form
of circular inclusions.  Let $Y=[0,\delta]^d$ be a \emph{reference cell}
that captures the configuration of the periodic microstructure, and for
which $\delta$ represents its actual physical size with $0\leq\delta\ll
1$ (although we focus on dimension $d=2$, the theoretical results proven
here apply to any dimension).  The reference cell is divided into
two sub-domains $Y^1$ and $Y^2$ that are separated by a Lipschitz
boundary $\Gamma=Y^1 \cap Y^2$ as shown in
Figure~\ref{fig:ice-refcell}b.  For simplicity, we take $\Gamma$ to be a
circle of radius $\gamma$ satisfying $0<\gamma<\frac{1}{2}\,\delta$.  The
primary feature that we exploit in our homogenization approach is that
within $Y^1$ heat must diffuse rapidly, whereas in $Y^2$ there is a
relatively slow diffusion of heat.

We next introduce a small parameter $0<\eps\ll 1$ that corresponds to
the size of the periodic microstructure (and must be distinguished
from the physical size $\delta$ because we will eventually take the
limit as $\eps\to 0$).  The domain $\Omega$ may then be decomposed
into three $\eps$-dependent sub-domains: $\Omoneeps\defeq
\text{int}\bigcup_{k\in\Z^d}\eps(k+\overline{Y^1})\cap \Omega$ (which
is connected), and two disconnected components consisting of the
region $\Omtwoeps \defeq \bigcup_{k\in\Z^d}\eps(k+Y^2)\cap\Omega$ and
the boundary curves $\Gameps \defeq \bigcup_{k\in\Z^d}
\eps(k+\Gamma)\cap\Omega$. This decomposition is illustrated in
Figure~\ref{fig:ice-refcell}. To avoid technical difficulties, we
assume that $\Gameps$ does not touch the outer boundary of $\Omega$,
so that $\Gameps\cap\partial\Omega=\emptyset$ and
$\Omtwoeps\cap\partial\Omega=\emptyset$.

The major advantage of this reduced model is that the reference cell
problem simplifies significantly, with the only unknowns being
$s_\mathrm{iw}$ and temperature.
We proceed with the temperature and enthalpy equations.

\subsection{Temperature and enthalpy equations}
\label{sec:reduced-temperature}

Throughout the analytical developments of this paper, we employ what is
known as the two-phase formulation of the Stefan problem, in which the
heat diffusion equation is posed in a mixed form involving both
temperature and enthalpy.  Assuming that material properties of
water and ice remain constant, the temperature $T$ can be written as a
piecewise linear function of enthalpy $\HE$ as follows \citep{STEFAN5}
\begin{gather*}
  T = \widetilde{\omega}(\HE) = 
  \left\{\begin{array}{cl}
      \ds\frac{1}{c_\mathrm{i}} \HE, & \qquad \text{if $\HE < \HE_\mathrm{i}$}, \\[0.3cm]
      \Tcrit,               & \qquad \text{if $\HE_\mathrm{i} \leq \HE < \HE_\mathrm{w}$},\\[0.1cm]
      \ds\Tcrit + \frac{1}{c_\mathrm{w}}(\HE-\HE_\mathrm{w}),& \qquad \text{if $\HE_\mathrm{w} \leq \HE$},
    \end{array}\right.
\end{gather*}
where $c_\mathrm{w}$ and $c_\mathrm{i}$ denote specific heats of water and ice
respectively, and $\Tcrit=273.15\myunit{\degK}$ is the freezing point of
water (parameter values are listed in Table~\ref{tab:ics1}).  A
distinguishing feature of this temperature--enthalpy relationship is
that when temperature is equal to the freezing point, the enthalpy
varies while temperature remains constant -- this behavior derives from
the fact that a certain amount of energy (called latent heat) is
required to effect a change in phase from solid to liquid at the phase
interface.

\begin{figure}[tbhp]
  \centering\footnotesize
  \begin{tabular}{ccc}
    \includegraphics[width=0.55\textwidth]{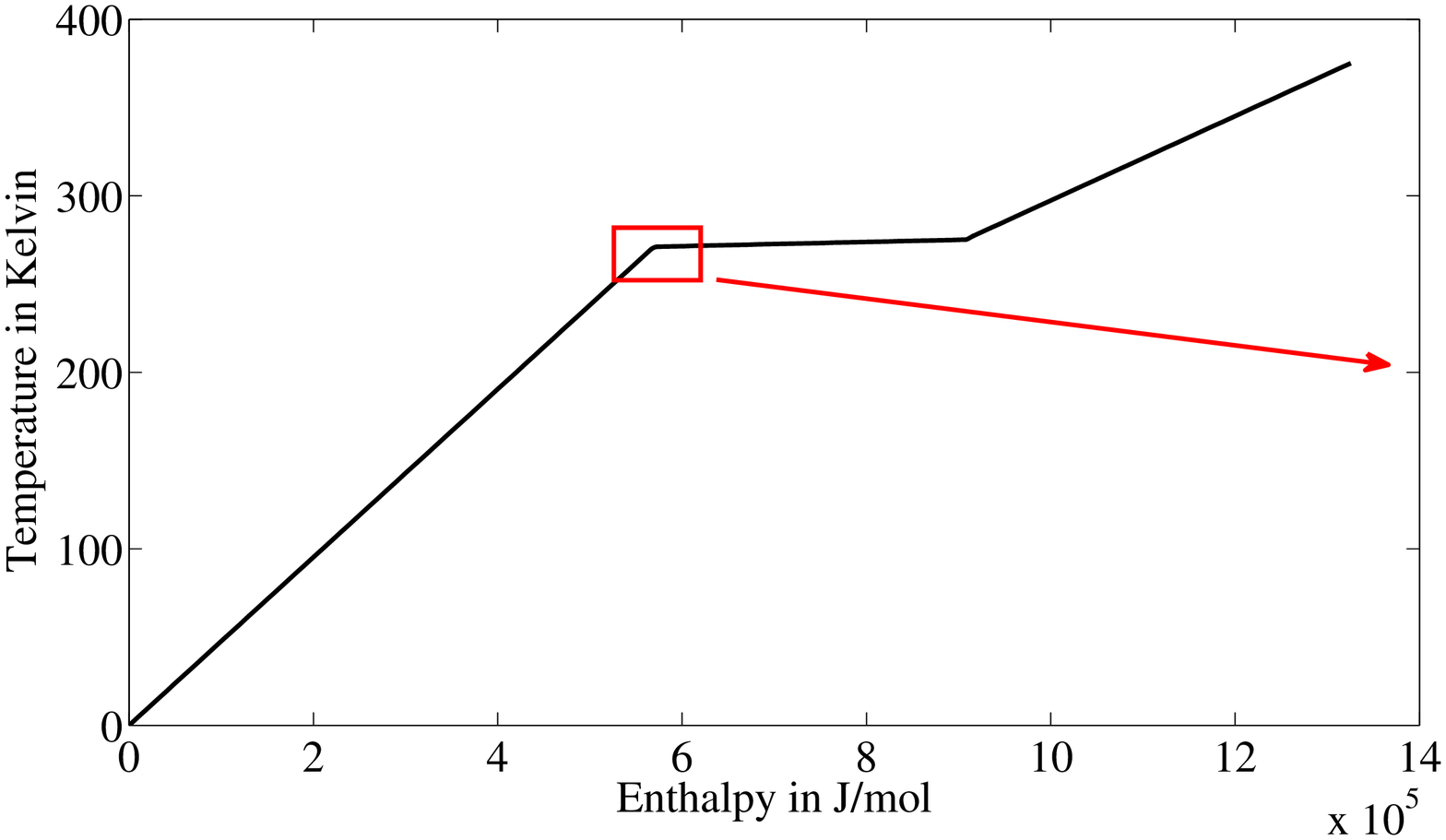}
    &\quad&
    \includegraphics[width=0.36\textwidth]{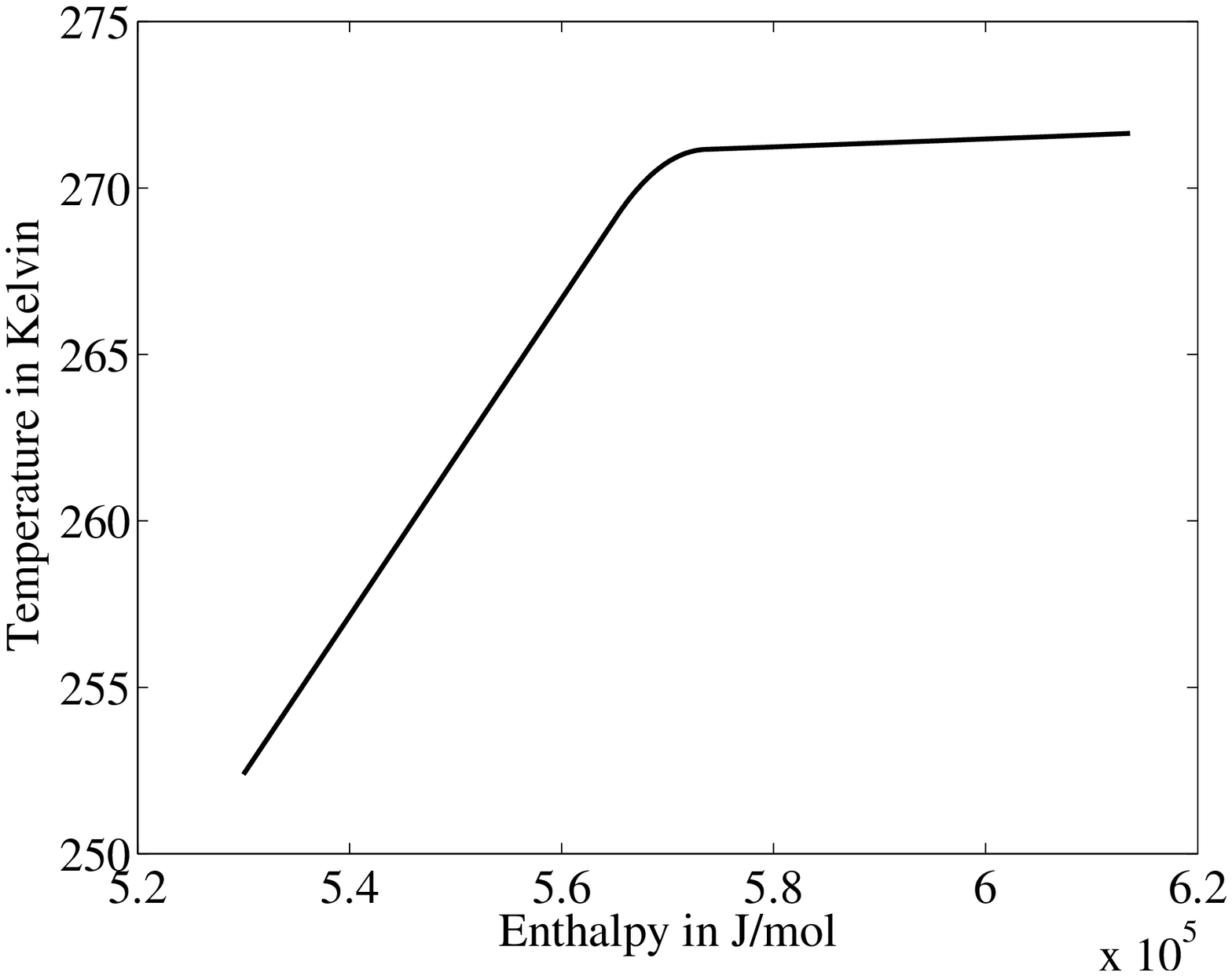}
  \end{tabular}
  \caption{Regularized temperature--enthalpy function, $T=\omega(\HE)$,
    for a generic ice/water phase change problem in which the critical
    (freezing) temperature is $\Tcrit=273.15\degK$.  The zoomed-in view
    on the right illustrates the smoothed corners in the
    regularization.}
  \label{fig:omega-smooth}
\end{figure}

Because the function $\widetilde{\omega}(\HE)$ is neither
differentiable nor invertible, we instead employ in our model a
regularized version $\omega(\HE)$ defined as
\begin{gather}\label{def_omega}
  T = \omega(\HE) = 
  \left\{\begin{array}{cl}
      \ds\frac{1}{c_\mathrm{i}} \HE,                   
      & \qquad \text{if $\HE < \HE_{\mathrm{i-}}$},\\[0.2cm]
      \text{\footnotesize[smooth connection]},
      & \qquad \text{if $\HE_{\mathrm{i-}} \leq \HE < \HE_{\mathrm{i+}}$},\\[0.2cm]
      \ds\Tcrit-\frac{2\HE - (\HE_{\mathrm{i+}} + \HE_{\mathrm{w-}})}{2c_\infty},
      & \qquad \text{if $\HE_{\mathrm{i+}} \leq \HE < \HE_{\mathrm{w-}}$},\\[0.2cm]
      \text{\footnotesize[smooth connection]},
      & \qquad \text{if $\HE_{\mathrm{w-}} \leq \HE < \HE_{\mathrm{w+}}$},\\[0.2cm]
      \ds\Tcrit + \frac{1}{c_\mathrm{w}}(\HE-\HE_{\mathrm{w+}}),
      & \qquad \text{if $\HE_{\mathrm{w+}} \leq \HE$},
    \end{array}\right.
\end{gather}
which has ``rounded corners'' that are smoothed over the short intervals
$\HE_{\mathrm{i-}} \lessapprox \HE_\mathrm{i} \lessapprox
\HE_{\mathrm{i+}}$ and $\HE_{\mathrm{w-}} \lessapprox \HE_\mathrm{w}
\lessapprox \HE_{\mathrm{w+}}$.  Note that we have also introduced a
small positive slope $c_\infty^{-1} \ll 1$ within the central
plateau region near $T\approx \Tcrit$ (refer to
Figure~\ref{fig:omega-smooth}).  These modifications ensure that
$\omega$ is a continuously differentiable, invertible and monotone
increasing function of enthalpy.  Incidentally, such a regularized
function is most likely a more accurate representation of what one would
actually observe in a real physical system.

We now describe the solution decomposition into slow and fast diffusing
variables.  Let functions $\Toneeps$ and $\HEoneeps$ denote the
fast-diffusing temperature and enthalpy components respectively, with
both defined on the sub-region $\Omoneeps$. Similarly, let $\Ttwoeps$
and $\HEtwoeps$ denote the slowly-diffusing temperature and enthalpy on
$\Omtwoeps$.  We may then state the strong formulation of the two-phase
Stefan problem as
\begin{subequations}
  \label{micro_eq_strong}
  \begin{align}
    \partial_t\HEoneeps - \nabla\cdot[D(\HEoneeps)\nabla \Toneeps] =
    0 & \qquad\text{in}\ \Omoneeps, \label{micro_eq_strong1}\\ 
    D(\HEoneeps)\nabla \Toneeps\cdot \normvec =
    -\eps^2D(\HEtwoeps)\nabla\Ttwoeps\cdot\normvec & \qquad\text{on}\
    \Gameps, \label{micro_eq_strong2}\\ 
    -D(\HEoneeps)\nabla T_{1,\eps}\cdot \normvec = \alpha(T_{1,\eps} - \Tout)
     & \qquad\text{on}\ 
    \partial\Omega\cap\partial\Omoneeps, \label{micro_eq_strong3}  
  \end{align}
  \begin{align}
    \partial_t\HEtwoeps - \eps^2\nabla\cdot[D(\HEtwoeps)\nabla \Ttwoeps]
    = 0 & \qquad\text{in}\ \Omtwoeps, \label{micro_eq_strong4}\\ 
    \HEtwoeps = \HEoneeps & \qquad\text{on}\ \Gameps, 
    \label{micro_eq_strong5}   
  \end{align}
\end{subequations}
where $D(\HE)$ is given by~\eqref{DH} and $\Tout$ is the ambient
temperature imposed on the outer domain boundary.  Note that these equations
capture the phase interface location implicitly through the relationship
$T=\omega(\HE)$ and no explicit Stefan-type phase interface
condition is imposed.

\section{Two-scale homogenization of the reduced model}
\label{sec:analysis}

It is not feasible to solve the system of equations
\eqref{micro_eq_strong} derived in the previous section using direct
numerical simulations owing to the presence of the microstruture with
$\eps$ being very small.  This section contains the primary theoretical
results pertaining to an upscaling of the reduced problem which is
achieved by characterising the limit as $\eps \to 0$. Lemmas and
theorems are stated here, and the proofs are relegated to the Appendix.

\subsection{Weak formulation}
In order to make the problem ameanable to periodic-homogenization
techniques, we begin by transforming Eqs.~\eqref{micro_eq_strong}
into a weak formulation. To avoid technical difficulties, we replace the
Robin boundary condition \eqref{micro_eq_strong3} by a Dirichlet
condition $\HEoneeps =\omega^{-1}(\Tout)$ (i.e.,~we consider the formal
limit $\alpha \to \infty$) and refer the reader
to~\citet{graf-peter-2014} for a detailed discussion.  This requires
first defining some appropriate solution spaces:
\begin{align*}
  \V^1_\eps &\defeq \left\{ u \in L^2([0,\tend], \Hil^1(\Omoneeps)) \cap
    \Hil^1([0,\tend], \Hil^1(\Omoneeps)')\; | \; u=0\text{ on
    }\partial\Omoneeps\cap\partial\Omega \right\},\\
  \V^2_\eps &\defeq \left\{ u \in L^2([0,\tend], \Hil^1(\Omtwoeps)) \cap
    \Hil^1([0,\tend], \Hil^1(\Omtwoeps)')\; | \; u=0\text{ on }\Gameps \right\},\\
  \V&\defeq L^2([0,\tend],\Hil^1_0(\Omega))\cap \Hil^1([0,\tend],\Hil^1(\Omega)'),
\end{align*}
where the ``primes'' denote dual spaces and $[0,\tend]$ represents the
time interval of interest for some fixed $\tend>0$.  The corresponding
test spaces are $V^1_\eps=\{u\in\Hil^1(\Omoneeps) \; | \; u=0 \text{ on
    }\partial\Omoneeps\cap\partial\Omega\}$,
$V^2_\eps=\Hil^1_0(\Omtwoeps)$ and $V = \Hil^1_0(\Omega)$.  We also need
to introduce notation for inner products, with $(u,v)_{\Omeps^\alpha} =
\int_{\Omeps^\alpha} uv\,\text{d}x$ representing the $L^2$-inner product
with respect to space of two functions in $\V^\alpha_\eps$ for
$\alpha=1,2$, whereas $(u,v)_{\Omeps^\alpha,\,t} =
\int_0^{t}\int_{\Omeps^\alpha} uv\,\text{d}x\, \text{d}\tau$ denotes
that an additional time integration is performed over the interval
$[0,t]$ with $0\leq t\leq \tend$.  Finally, we let $\langle
u,v\rangle_{\Gameps} = \langle u,v\rangle_{H^\frac12(\Gameps)'\times
  H^\frac12(\Gameps)}$ denote the dual pairing on $\Gameps$. Later, we
will show that $D(\HEtwoeps)\nabla \Ttwoeps\cdot \normvec \in
L^2(\Gameps)$, so that we can interpret $\langle u,v\rangle_{\Gameps}$
as $\int_{\Gameps}g_\eps uv\,\text{d}S$, where $g_\eps$ represents the
Riemann curvature tensor.

We are now prepared to state the weak form of the heat-diffusion
problem.  Assuming that initial values $T_{1,\eps,\mathrm{init}} =
\omega(\HE_{1,\eps,\mathrm{init}})$ and $T_{2,\eps,\mathrm{init}} =
\omega(\HE_{2,\eps,\mathrm{init}})$ are smooth, non-negative and bounded
functions, and that a Dirichlet condition $\Toneeps = \Tout$ is imposed
at the outer boundary $\partial\Omega\cap\partial\Omoneeps$, our goal is
to find $(\Toneeps, \Ttwoeps) \in (\V^1_\eps+\Tout) \times
(\V^2_\eps+\Toneeps)$ such that
\begin{subequations}\label{s_problem_weak}
  \begin{align}
    (\partial_t \HEoneeps,\; \phi)_{\Omoneeps} + (D(\HEoneeps)\nabla
    \Toneeps,\; \nabla\phi)_{\Omoneeps} 
    + \eps^2 \langle D(\HEtwoeps) \nabla \Ttwoeps\cdot \normvec,\;
    \phi\rangle_{\Gameps} &= 0, \label{s_problem_weak_1}\\  
    (\partial_t \HEtwoeps,\; \psi)_{\Omtwoeps} + (\eps^2
    D(\HEtwoeps)\nabla \Ttwoeps,\; \nabla\psi)_{\Omtwoeps} 
    &= 0,
    \label{s_problem_weak_2}
  \end{align}
\end{subequations}
for all $\phi, \psi \in V^1_\eps\times V^2_\eps$.  Note that
$\normvec$ represents the outward-pointing unit normal vector on
$\Gameps$, and that temperature and enthalpy are connected via $\Taeps
= \omega(\HEaeps)$, or equivalently $\HEaeps=\omega^{-1}(\Taeps)$.  We
assume that $\omega^{-1}(\Tout)$ is positive, bounded and smooth such
it can be extended to $[0,\tend]\times\Omega$ in $\Hil^1$.  We note
again in closing that slow diffusion is induced in the problem via the
factor $\eps^2$ multiplying terms in Eqs.~\eqref{s_problem_weak} that
involve the diffusion coefficient $D(\HEtwoeps)$.

\subsection{Transformation of the model}
\label{sec:transformation}

In this section, we apply a procedure developed by 
\citet{PU5} to transform the model \eqref{s_problem_weak} by
combining enthalpies $\HEoneeps$ and $\HEtwoeps$ into a single
function $\Theta_\eps$ defined on the whole $\eps$-independent
domain $\Omega$.  We use the fact that $\Gameps$ is the only boundary of
$\Omtwoeps$ to obtain
\begin{equation}\label{transition}
\begin{aligned}
  \langle\eps^2D(\HEtwoeps)\nabla\Ttwoeps\cdot
  \normvec,\psi\rangle_{\Gameps}  
  &= (\eps^2\nabla\cdot[D(\HEtwoeps)\nabla \Ttwoeps],\psi)_{\Omtwoeps} +
  (\eps^2D(\HEtwoeps)\nabla \Ttwoeps,\nabla\psi)_{\Omtwoeps}\\ 
  &= (\partial_t\HEtwoeps,\psi)_{\Omtwoeps} + (\eps^2D(\HEtwoeps)\nabla
  \Ttwoeps,\nabla\psi)_{\Omtwoeps}, 
\end{aligned}
\end{equation}
for all $\psi\in V$.  After substituting this expression into
\eqref{s_problem_weak_1} we obtain
\begin{align*}
  (\partial_t\HEoneeps,\phi)_{\Omoneeps} + (D(\HEoneeps)\nabla
  \Toneeps,\nabla\phi)_{\Omoneeps} +
  (\partial_t\HEtwoeps,\psi)_{\Omtwoeps} 
  + (\eps^2D(\HEtwoeps)\nabla\Ttwoeps,\nabla\psi)_{\Omtwoeps} = 0, 
\end{align*}
for all $\phi, \psi \in V^1_\eps\times V^2_\eps$. Hence,
Eqs.~\eqref{s_problem_weak} have been replaced with 
\begin{equation}\label{problem_weak_together1}
  \begin{aligned}
    (\partial_t\HEoneeps + \partial_t\HEtwoeps,\phi)_{\Omega} +
    (D(\HEoneeps)\nabla \Toneeps + \eps^2 D(\HEtwoeps)\nabla
    \Ttwoeps,\nabla\phi)_\Omega = 0, &\mbox{}\\ 
    \HEoneeps = \HEtwoeps \quad\text{on}\ \Gameps, &\mbox{}
  \end{aligned}
\end{equation}
for all $\phi\in V$. We then define the function $\Theta_\eps \in
L^2([0,\tend],\Hil^1(\Omega))$ by  
\begin{align*}
  \Theta_\eps = \left\{
    \begin{array}{ll}\HEoneeps &\ \text{in}\ \Omoneeps,\\ 
      \HEtwoeps & \ \text{in}\ \Omtwoeps,
    \end{array}\right. 
\end{align*}
so that with conditions \eqref{micro_eq_strong2} and
\eqref{micro_eq_strong5} the function $\Theta_\eps$ is guaranteed to be
continuous and weakly differentiable.  Furthermore, we define
$\kappa_\eps = \chi_{1,\eps} + \eps^2\chi_{2,\eps}$ where
$\chi_{i,\eps}$ for $i=1,2$ are indicator functions for $\Omoneeps$
and $\Omtwoeps$ respectively.  Solving \eqref{s_problem_weak} is then
equivalent to finding $\Theta_\eps \in (\V + \omega^{-1}(\Tout))$ such
that
\begin{equation}
  \label{problem_weak_together2}
  (\partial_t\Theta_\eps,\phi)_\Omega + (\kappa_\eps
  D(\Theta_\eps)\omega'(\Theta_\eps)\nabla\Theta_\eps,\nabla\phi)_\Omega
  = 0 ,
\end{equation}
for all $\phi\in V$, where we have used that $\Toneeps =
\omega(\HEoneeps)$ and $\Ttwoeps = \omega(\HEtwoeps)$.

We perform one further transformation of
\eqref{problem_weak_together2} that makes the Dirichlet boundary
condition homogeneous.  To this end, we define $\varrho_\eps = \Theta_\eps -
\omega^{-1}(\Tout)$ where $\omega^{-1}(\Tout)$ is extended continuously to $\Omega$.  
Then \eqref{problem_weak_together2} is
equivalent to finding $\varrho_\eps \in\V$ such that
\begin{equation}
  \label{problem_weak_together3}
  (\partial_t\varrho_\eps,\phi)_\Omega + (\kappa_\eps
  D\omega'(\varrho_\eps +
  \omega^{-1}(\Tout))\nabla\varrho_\eps,\nabla\phi)_\Omega =
  (-\partial_t\omega^{-1}(\Tout),\phi)_\Omega ,
\end{equation}
for all $\phi\in V$.


\subsection{Existence of a weak solution}
\label{sec:existence}

\subsubsection{Theorem of existence}
\label{sec:theorem_existence}

To prove the existence of a solution to \eqref{problem_weak_together2}
for every $\eps>0$, we formulate a theorem of existence, which is strongly
inspired by a proof for a related result found in a set of
unpublished lecture notes by \citet{wolff-lecturenotes}.  In
Theorem \ref{theorem_of_existence} we introduce a general existence
result for parabolic equations with non-monotone non-linearities in
the diffusion operator. The proof is based on the Rothe method.  We
state the theorem next and provide the proof in Appendix~\ref{app:existence}.

Consider the initial--boundary value problem
\begin{subequations}\label{ext_problem1}
  \begin{alignat}{3}
    \partial_tu + \sum_{j=1}^n \partial_{x_j}\left(a(x,t,u) \partial_{x_j}u
    \right) &= F(x,t)\qquad && \text{in}\ S\times \Omega, \\ 
    u&=0\qquad && \text{on}\ S\times\partial\Omega, \\
    u(0,x) &= u_0(x)\qquad &&\text{in}\ \Omega, 
  \end{alignat}
\end{subequations}
with 
\begin{subequations}\label{ext_cond:abcde}
  \begin{align}
    & \Omega\subset \R^n \quad\text{bounded Lipschitz domain}, \qquad
    S=[0,\tend], \label{ext_cond:a}\\ 
    & a(x,t,u):\Omega\times S \times \R \rightarrow\R
    \quad\text{Bochner-measurable in } x \ \text{and continuous in } t \text{ and }
    u, \label{ext_cond:b}\\ 
    & \exists\ 0 < \lambda \leq \Lambda < \infty\ \text{such that}\
    \lambda\leq a(x,t,u)\leq \Lambda 
    \;\; \forall s \in \R, \; \text{for a.e.}\ t\in\R \; \text{and} \;
    \text{a.e.}\ x\in\Omega, \label{ext_cond:c}\\ 
    &V = \Hil^1_0(\Omega), \qquad \V = L^2(S,V),\qquad  \V^* =
    L^2(S,V^*), \label{ext_cond:d}\\ 
    & u_0 = L^2(\Omega),\qquad F\in \V^*. \label{ext_cond:e}
  \end{align}
\end{subequations}

\begin{lemma}\label{lemma_ext_prep}
  Let the conditions \eqref{ext_cond:abcde} be satisfied. Then 
  \begin{align}
    \langle A(t)(u,v),w\rangle &= \sum_{j=1}^n\int_\Omega
    a(x,t,u)\, \partial_{x_j}v\, \partial_{x_j}w\, \text{\textup{d}}x, \\ 
    \langle f,v\rangle &= \int_0^T\int_\Omega F(x,t)\, v(x,t) \, 
    \text{\textup{d}}x\, \text{\textup{d}}t,  
  \end{align}
  define a family of operators $A(t):V\times V\rightarrow V^*$ and an
  element $f\in \V^*$ for which the following hold:
  \begin{subequations}\label{ext_cond2:abcd}
    \begin{alignat}{3}
      &\forall\ u,v \in V,\ \text{for a.e.}\ t\in \R:\quad
      && \norm{A(t)(u,v)}_{V^*}\leq \Lambda\norm{v}_V, \label{ext_cond2:a}\\
      &\forall\ u \in V,\ \text{for a.e.}\ t\in \R: \quad && A(t)(u,\cdot): V
      \rightarrow V^* \ \text{is linear and
        continuous}, \label{ext_cond2:b}\\
      &\forall\ u,v \in V,\ \text{for a.e.}\ t\in\R: \quad &&\langle
      A(t)(u,v),v\rangle \geq \lambda\norm{v}^2_V, \label{ext_cond2:c}\\
      & \A:\V\times \V\rightarrow\V^* \ &&\text{is the realization of $A$
        with $\A(u,u) = A(t)(u,u)$}. \label{ext_cond2:d}
    \end{alignat}
  \end{subequations}
\end{lemma}
This lemma can be proven using standard arguments as described by
\citet{PDE16} and \citet{wolff-lecturenotes}.  Lemma \ref{lemma_ext_prep}
implies that any initial value problem of the form
\begin{equation}\label{ext_problem2}
  \begin{aligned}
    u^\prime  + \A(u,u) &= f \qquad\text{in}\ V^*,\\
    u(0) &= u_0,
  \end{aligned}
\end{equation}
which includes the problem \eqref{ext_problem1}, is well-defined.

\begin{theorem}\label{theorem_of_existence}
  \theoremone
\end{theorem}

The proof of Theorem \ref{theorem_of_existence} is given in Appendix~\ref{app:existence}.
We apply this theorem to the reduced model
\eqref{problem_weak_together3} for which the space $V$ is
$\Hil^1_0(\Omega)$, so that both solution and test space correspond to
$\V = L^2([0,\tend],\Hil^1_0(\Omega))$.  We have $u_0 = \varrho_\eps(0)$
and $f\equiv \partial_t\omega^{-1}(\Tout)$, which is in $\V^*$ by
assumption.  The function $a$ corresponding to the problem
\eqref{problem_weak_together3} is
\begin{align*}
  a(x,t,u)  = \kappa_\eps(x) D\omega'(u + \omega^{-1}(\Tout(t))), 
\end{align*}
which satisfies the conditions \eqref{ext_cond:abcde} for $u\in V$.
With all conditions of Theorem~\ref{theorem_of_existence}
fulfilled, we obtain the following result

\begin{corollary}
  There exists a solution of Eq.~\eqref{problem_weak_together3}.
\end{corollary}

\subsection{A~priori estimates and limit functions}
\label{sec:apriori}

Our main results on a~priori estimates are stated in
Lemma~\ref{lemma_estimates} below, which bounds independently of $\eps$
the functions $\HEoneeps$ and $\Toneeps$ in
$L^2([0,\tend],\Hil^1(\Omoneeps))$, and similarly $\HEtwoeps$ and
$\Ttwoeps$ in $L^2([0,\tend],\Hil^1(\Omtwoeps))$. 
The proof of this Lemma is given in Appendix~\ref{app:estimates}.

\begin{lemma}\label{lemma_estimates}
  \lemmatwo
\end{lemma}

Using standard results of two-scale convergence
\citep{TWOSCALE1,TWOSCALE2,TWOSCALE16}, we immediately obtain the
following result.

  \begin{lemma}\label{lemma_limits}
    There exist functions $\HE_{1,0} \in
    L^2([0,\tend],\Hil^1(\Omega))$, $\widehat{E}_{1,0}\in
    L^2([0,\tend],L^2(\Omega,\Hil^1_\#(Y^1)))$ and $\HE_{2,0} \in
    L^2([0,\tend], L^2(\Omega,\Hil^1_\#(Y^2)))$ such that, up to
    subsequences,
    \begin{align*}
      \HE_{1,\eps} & \xrightarrow{\text{~2-scale~}} \HE_{1,0},\\
      \nabla \HE_{1,\eps} & \xrightarrow{\text{~2-scale~}} \nabla_x
      \HE_{1,0} + \nabla_y\widehat{E}_{1,0},\\ 
      \HE_{2,\eps} & \xrightarrow{\text{~2-scale~}} \HE_{2,0},\\
      \nabla \HE_{2,\eps} & \xrightarrow{\text{~2-scale~}} \nabla_y\HE_{2,0}.
    \end{align*}
  \end{lemma}

  Note that $\HE_{1,0}$ is independent of $y$, and we have also
  introduced $\widehat{E}_{1,0}\in
  L^2([0,\tend],L^2(\Omega,\Hil^1_\#(Y^1)))$ and $\HE_{2,0}\in
  L^2([0,\tend],L^2(\Omega,\Hil^1_\#(Y^2)))$, where the subscript $\#$
  denotes $Y$-periodicity in space. The limit of $\nabla E_{1,\eps}$ has
  a special form obtained in~\citet{TWOSCALE1} that consists of
  two terms: one involving a gradient with respect to the slow variable,
  and a second term with respect to the fast variable.

\subsection{Identification of the two-scale limit}
\label{sec:two-scale}

Owing to the nonlinear dependence of the diffusion coefficient on
enthalpy in \eqref{problem_weak_together3}, we have not yet been able to
identify the system of equations satisfied by the limit functions of
Lemma~\ref{lemma_limits} without further assumptions. In order to allow
the limit passage without difficulties, we \emph{assume} in the
remainder of this section that the function $D\omega'$
in~\eqref{problem_weak_together3} is independent of $\Theta_\eps$, which
makes the model linear.  Note that having strong convergence of the
function $\Theta_\eps$ in $L^2([0,\tend],\Omega)$ would lead to the same
results. A homogenization proof for the fully nonlinear problem is left
for future work.

In order to characterize the limit functions from Lemma
\ref{lemma_limits}, we define test functions that vary on length
scales of size $O(1)$ and $O(\eps)$ according to
\begin{gather*}
  \phi_\eps\left(x,\textstyle\frac{x}{\eps}\right) =
  \chi_1 \left(\textstyle\frac{x}{\eps}\right) 
  \left(\phi_0(x) + \eps\phi_1 \left(x,\textstyle\frac{x}{\eps}\right)\right) +
  \chi_2 \left(\textstyle\frac{x}{\eps}\right) 
  \phi_2(x,\textstyle\frac{x}\eps),
\end{gather*}
where 
$(\phi_0,\phi_1,\phi_2) \in C^\infty_0(\Omega)\times
C^\infty(\Omega,C^\infty_\#(Y))\times C^\infty(\Omega,C^\infty_\#(Y))$.

By substituting $\phi_\eps$ into \eqref{s_problem_weak_1} and using
$\chi_i$ to write the resulting integrals over the entire domain
$\Omega$, we obtain
\begin{multline*}
  \int_{\Omega}\chi_1
  \left(\textstyle\frac{x}{\eps}\right) \partial_t\HEoneeps(x,t)
  \phi_\eps \left(x,\textstyle\frac{x}{\eps}\right) \, \text{d}x +
  \int_{\Omega}\chi_1 \left(\textstyle\frac{x}{\eps}\right)D\omega'
  \nabla \HEoneeps(x,t) \nabla\phi_\eps
  \left(x,\textstyle\frac{x}{\eps}\right) \,\text{d}x\\
  +\int_{\Omega}\chi_2 \left(\textstyle\frac{x}{\eps}\right)\partial_t
  \HEtwoeps(x,t)\phi_\eps
  \left(x,\textstyle\frac{x}{\eps}\right)\,\text{d}x +
  \int_\Omega\chi_2\left(\textstyle\frac{x}{\eps}\right) D\omega'
  \eps^2\nabla
  \HEtwoeps(x,t)\phi_\eps\left(x,\textstyle\frac{x}{\eps}\right)\,\text{d}x=
  0.
\end{multline*}
Then, taking the limit as $\eps\to 0$ yields
\begin{multline}
  \label{eq:limit1}
  \int_{\Omega\times Y^1}\partial_t \HE_{1,0}(x,t)\phi_0(x)\,
  \text{d}y\, \text{d}x 
  + \int_{\Omega\times Y^1}D\omega'[\nabla_xE_{1,0}(x,t) + \nabla_y
  \widehat{E}_{1,0}(x,y,t)]
  [\nabla_x\phi_0(x) + \nabla_y\phi_1(x,y)]\, \text{d}y\, \text{d}x\\
  + \int_{\Omega\times
    Y^2}\partial_t\HE_{2,0}(x,y,t)\phi_2(x,y)\;\text{d}y\, \text{d}x 
  + \int_{\Omega\times
    Y^2}D\omega'\nabla_y\HE_{2,0}(x,y,t)\nabla_y\phi_2(x,y)\, \text{d}y\,
  \text{d}x=
  0,
\end{multline}
where $y$ denotes the spatial variable on the reference cell $Y$. 

We are free at this point to choose any test function and so we take
$\phi_0=0$ and $\phi_2=0$ in Eq.~\eqref{eq:limit1}. To start with, we introduce
functions $\mu_k\in \Hil^1_\#(Y^1)$ in order to express
$\widehat{E}_{1,0}(x,y,t) = \sum_{k=1}^d\partial_{x_k}E_{1,0}(x,t)
\mu_k(y)$ in separable form.  The weak formulation of the cell problem
for $k=1,\dots,d$ may then be expressed in the simpler form
\begin{align}
  \label{cell_problem1}
  (e_k + \nabla_y\mu_k, \; \nabla_y\phi_1)_{Y^1} = 0,
\end{align}
where the $\mu_k$ are $Y$-periodic.  Alternatively, we may take $\phi_1=0$
in Eq.~\eqref{eq:limit1} to obtain
\begin{multline*}
  \int_{\Omega\times Y^1}\partial_t \HE_{1,0}(x,t)\phi_0(x)\,
  \text{d}y\, \text{d}x 
  + \int_{\Omega\times Y^1}D\omega' \sum_{k=1}^d \partial_{x_k}
  E_{1,0}(x,t) [e_k + \nabla_y \mu_k(y)] \nabla_x\phi_0(x)\, \text{d}y\,
  \text{d}x\\ 
  + \int_{\Omega\times Y^2} 
  \partial_t\HE_{2,0}(x,y,t)\phi_2(x,y)\;\text{d}y\, \text{d}x 
  + \int_{\Omega\times Y^2}D\omega'\nabla_y\HE_{2,0}(x,y,t)
  \nabla_y\phi_2(x,y)\;\text{d}y\, \text{d}x= 0,
\end{multline*}
which can be rewritten in the more suggestive form
\begin{multline}\label{eq:limit3}
  \int_{\Omega\times Y^1}\partial_t \HE_{1,0}(x,t)\phi_0(x)\,
  \text{d}y\, \text{d}x 
  + \int_{\Omega}D\omega'
  \sum_{k,\ell=1}^d\partial_{x_k}E_{1,0}(x,t) \int_{Y^1}[\delta_{k\ell}
  + \partial_{y_\ell}
  \mu_k(y)]\, \text{d}y\, \partial_{x_\ell}\phi_0(x)\, \text{d}x\\
  + \int_{\Omega\times Y^2} 
  \partial_t\HE_{2,0}(x,y,t) \phi_2(x,y)\;\text{d}y\, \text{d}x
  +  \int_{\Omega\times Y^2} 
  D\omega'\nabla_y\HE_{2,0}(x,y,t)
  \nabla_y\phi_2(x,y)\;\text{d}y\, \text{d}x = 0.
\end{multline}
The diffusion term involves the factors
\begin{align}
  \label{cell_problem2}
  \Pi_{k\ell} = \int_{Y^1}(\delta_{k\ell} + \partial_{y_\ell}\mu_k)
  \,\text{d}y 
\end{align}
for $k,\ell=1,\dots,d$, which can be represented as a matrix $\Pi$ that
multiplies the diffusion coefficient $D\omega'$.  We then obtain from
\eqref{eq:limit1} and \eqref{eq:limit3} the equation
\begin{multline}
  \label{s_problem_limit_inter}
  \abs{Y^1}(\partial_t\HE_{1,0},\; \phi_0)_{\Omega} + (\Pi
  D\omega'\nabla_x \HE_{1,0},\; \nabla_x\phi_0)_{\Omega} + (\partial_t
  \HE_{2,0},\; \phi_2)_{\Omega\times Y^2} + (D\omega' \nabla_y \HE_{2,0},\;
  \nabla_y\phi_2)_{\Omega\times Y^2} = 0.
\end{multline}

As a final step, we obtain the limit equation for $\HE_{2,0}$ by setting
$\phi_0=0$ in Eq.~\eqref{s_problem_limit_inter}, and similarly for
$\HE_{1,0}$ by setting $\phi_2=0$ in Eq.~\eqref{s_problem_limit_inter},
which is found using a similar transition as in \eqref{transition}.  The
resulting limit equations are
\begin{subequations}\label{s_problem_limit}
  \begin{align}
    \abs{Y^1}(\partial_t\HE_{1,0},\; \phi_0)_{\Omega} + (\Pi
    D\omega'\nabla_x \HE_{1,0},\; \nabla_x\phi_0)_{\Omega}
    + \langle D\omega'\nabla_y\HE_{2,0},\;
    \phi_0\rangle_{\Gamma\times\Omega} & = 0,\\ 
    (\partial_t \HE_{2,0},\; \phi_2)_{\Omega\times Y^2} + (D\omega'
    \nabla_y \HE_{2,0},\; \nabla_y\phi_2)_{\Omega\times Y^2} & = 0,
  \end{align}
\end{subequations}
for all $\phi_0\in \Hil^1_0(\Omega)$ and $\phi_2\in
L^2(\Omega,\Hil^1_\#(Y^2))$, where $\Pi$ is the $d\times d$ matrix of
scaling factors defined in \eqref{cell_problem2}, $\HE_{1,0}\in
L^2([0,\tend],\Hil^1(\Omega))$ with $\omega(\HE_{1,0})=\Tout$ on
$\partial\Omega$, and $\HE_{2,0}\in L^2([0,\tend],
L^2(\Omega,\Hil^1_\#(Y^2)))$ with $\HE_{2,0}=\HE_{1,0}$ on
$\Omega\times\Gamma$.

To simplify notation in the remainder of the paper, we drop the zero
subscripts in $\{T_{1,0},\HE_{1,0}, T_{2,0},\HE_{2,0}\}$ and denote them
instead by $\{T_1,\HE_1,T_2,\HE_2\}$. Note again that we have only
rigorously derived the limit problem in the linear case and so we would
need to prove strong convergence of the function $\Theta_\eps$ in
$L^2([0,\tend], \Omega)$ for the analysis to hold for
\eqref{problem_weak_together2}; we will nevertheless transition back to
the nonlinear problem with an enthalpy-dependent diffusion coefficient
$D\omega'(\HE)$, for which the corresponding limit equations are
\begin{subequations}
  \label{s_problem_limit_nonlinear}
  \begin{align}
    \abs{Y^1}(\partial_t\HE_{1},\; \phi_0)_{\Omega} + (\Pi
    D(\HE_{1})\omega'(\HE_{1})\nabla_x \HE_{1},\;
    \nabla_x\phi_0)_{\Omega} 
    + \left\langle D(\HE_{2})\omega'(\HE_{2})\nabla_y\HE_{2},\;
      \phi_0 \right\rangle_{\Gamma\times\Omega} &= 0,\\ 
    (\partial_t \HE_{2},\; \phi_2)_{\Omega\times Y^2} +
    \left( D(\HE_{2})\omega'(\HE_{2}) \nabla_y \HE_{2},\;
      \nabla_y\phi_2 \right)_{\Omega\times Y^2} &= 0.
  \end{align}
\end{subequations}

\subsection{Uniqueness}
\label{sec:uniqueness}

The uniqueness of the solution to the nonlinear problem
\eqref{s_problem_limit_nonlinear} subject to suitable boundary and
initial conditions may be formulated compactly in terms of the following
theorem, which is proven in Appendix~\ref{app:uniqueness}.
\begin{theorem}\label{uniqueness}
  \theoremtwo
\end{theorem}

We note that the uniqueness of the limit problem implies that already
the whole sequences of solutions converge to the functions satisfying
\eqref{s_problem_limit_nonlinear}.

\subsection{Strong formulation of the limit problem}
\label{sec:strong-form}

We now state an equivalent strong formulation of the limit problem
corresponding to the weak form in~\eqref{s_problem_limit_nonlinear},
but with the Dirichlet condition at the outer boundary switched back to a
Robin condition again. This consists of a PDE for $T_1$ and $\HE_1$ on
the macroscale domain $\Omega$
\begin{subequations}
  \label{limit_problem_weak}
  \begin{alignat}{3}
    \abs{Y^1}\partial_t\HE_1 - \nabla_x\cdot(\Pi D(\HE_1)\nabla_xT_1) =
    \int_\Gamma D(\HE_2)\nabla_y T_2\cdot \normvec \, \text{d}S
    & &\text{in}\ \Omega, 
    \label{limit_problem_weak1a}\\
    -D(\HE_1)\nabla_x T_1\cdot \normvec = \alpha(T_1 - \Tout) & &
    \text{on}\ \partial\Omega, 
    \label{limit_problem_weak1b}\\
    \intertext{along with a second PDE for $T_2$ and $\HE_2$ on the
      microscale}
    \partial_t\HE_2 - \nabla_y\cdot (D(\HE_2)\nabla_y T_2) = 0 
    & & \text{  on}\ \Omega\times Y^2, 
    \label{limit_problem_weak2a}\\
    T_2 = T_1 
    & & \text{on}\ \Omega\times \Gamma, 
    \label{limit_problem_weak2b}\\
    \intertext{and initial values for enthalpy that we denote
      $\HE_{1,\mathrm{init}}$ and $\HE_{2,\mathrm{init}}$.  These two
      problems are coupled through the heat flux integral term in
      \eqref{limit_problem_weak1a} and the matching condition
      \eqref{limit_problem_weak2b}, both of which are enforced on
      $\Gamma$.  This is again the \emph{two-phase formulation of the
        Stefan problem}, which contains no explicit equation for the
      motion of the phase interface; instead, the interface location is
      captured implicitly through the temperature--enthalpy relation}
    T_1 = \omega(\HE_1) \ \text{in}\ \Omega \qquad\text{and}\qquad T_2 =
    \omega(\HE_2) & & \text{in}\ \Omega\times Y^2.
    \label{limit_problem_weak3}
  \end{alignat}
\end{subequations}
Under the assumption that temperature within the ice phase is constant
in space during a thawing event \citep{STEFAN5}, the problem
\eqref{limit_problem_weak} may be rewritten in an equivalent
\emph{one-phase formulation} that obeys the same macroscale problem
\begin{subequations}
  \label{limit_problem_strong}
  \begin{align}
    \abs{Y^1}\partial_t\HE_1 - \nabla_x\cdot(\Pi D(\HE_1)\nabla_xT_1) =
    \int_\Gamma D(\HE_2)\nabla_y T_2\cdot \normvec\, \text{d}S
    \quad \text{in}\ \Omega,
    \label{limit_problem_strong1a}
  \end{align}
  \begin{alignat}{3}
    -D(\HE_1)\nabla_x T_1\cdot \normvec = \alpha(T_1 - \Tout) & &
    \text{on}\ \partial\Omega.
   \label{limit_problem_strong1b}\\
    \intertext{On the reference cell, however, the ice temperature is
      taken equal to $\Tcrit$ and the water temperature obeys the
      following microscale equations}
    c_\mathrm{w}\partial_t T_2 - \nabla_y\cdot
    (D(\HE_2)\nabla_y T_2) = 0 
    & & \text{on}\ \Omega\times \widetilde{Y}^2(x,t),
    \label{limit_problem_strong2a}\\ 
    T_2 = T_1 & & \text{on}\ \Omega\times \Gamma,
    \label{limit_problem_strong2b}\\ 
    T_2 = \Tcrit & & \text{on}\ \Omega\times\partial \widetilde{Y}^2(x,t),
    \label{limit_problem_strong2c}\\
    \intertext{which are solved only on the water-filled annular region
      $\widetilde{Y}^2(x,t)\subseteq {Y}^2$ lying between $\Gamma$ and
      the moving phase boundary $s_\mathrm{iw}(x,t)$.  Consequently, in this
      one-phase formulation both the domain $\widetilde{Y}^2$ and its
      boundary $\partial\widetilde{Y}^2$ depend on $x$ and $t$ through
      $s_\mathrm{iw}$.  In the case of a freezing event, the ordering of the
      ice/water layers is reversed in which case the water temperature
      is held constant at $\Tcrit$ instead and $\widetilde{Y}^2$
      corresponds to the sub-region containing ice.  Finally, rather than
      imposing a temperature--enthalpy relation, the one-phase
      formulation fixes the temperature on the phase boundary via
      \eqref{limit_problem_strong2c} and provides an explicit Stefan
      condition governing the dynamics of the phase interface}
    \partial_t s_\mathrm{iw} =
    -\frac{D(\HE_2)}{(\HE_\mathrm{w}-\HE_\mathrm{i})}\nabla_y T_2\cdot \normvec
    &\quad &\text{on}\ \Omega\times\partial \widetilde{Y}^2(x,t).
    \label{limit_problem_strong3}
  \end{alignat}
\end{subequations}
The primary reason that we employ the one-phase formulation of the
Stefan problem is that it makes numerical simulations of the limit
problem much more convenient.  A detailed derivation of this one-phase
formulation from the corresponding two-phase formulation can be found
in~\citet{STEFAN5}.

\subsection{Limit problem for the sap exudation model}
\label{sec:sap-limit-problem}

Based on the limit problem we just derived for the reduced model using
homogenization techniques, it is now straightforward to pose the
analogous limit problem for the sap exudation model.  The two-scale heat
transport
equations~\eqref{limit_problem_strong1a}--\eqref{limit_problem_strong2c}
remain identical, but the Stefan condition~\eqref{limit_problem_strong3}
in the reduced model is replaced by the full set of
differential--algebraic equations (DAEs)
\eqref{cell_problem_maple}--\eqref{algebraic} for the microscale sap
exudation problem.

Although we have only performed the periodic homogenization procedure on
the reduced model, there are several features of the sap exudation
problem that can be exploited to extend our analytical results:
\begin{enumerate}
\item \emph{Presence of the gas phase:} which takes the form of gas
  bubbles in both fiber and vessel and introduces a spatial dependence
  in the thermal diffusion coefficient, $D(\HE,x)$.  Extending our
  analytical results to the case when $D$ also depends on $x$ would be a
  straightforward generalization.

\item \emph{Dissolved sugar in the vessel sap:} which gives rise to an
  osmotic potential between fiber and vessel that is essential for
  generating realistic exudation pressures.  Sugar within the vessel sap
  also depresses the freezing point so that the function $\omega$ in
  \eqref{def_omega} differs between vessel and fiber.  Although we do not
  need to consider freezing point depression explicitly in this paper
  (since we treat only a single thawing cycle) this effect could still
  be incorporated into the analysis, for example by adding an extra
  spatial dependence in $\omega$. Alternatively, the fiber could be
  defined as separate domain that is connected to the vessel via
  appropriate boundary conditions, thereby ensuring that the
  homogenization results carry through for the sap exudation model.  We
  have chosen not to incorporate this effect into the analysis, although
  periodic homogenization has previously been applied to Stefan problems
  having various functional forms for $\omega$
  in~\citet{bossavit-damlamian-1981}.

\item \emph{Extension to a freezing cycle:} which requires modifications
  only to the microscale equations in the reference cell as outlined
  in~\citet{graf-ceseri-stockie-2015}.  Consequently, this extension has
  no effect on the homogenization procedure.
  %
\end{enumerate}

\section{Multiscale numerical simulations}
\label{sec:simulations}

\subsection{Solution algorithm}
\label{sec:algo}

We now propose a multiscale numerical solution algorithm that computes
approximate solutions to both the reduced and sap exudation models.  The method is
based on a time-splitting approach that alternates in each time step
between solving the microscale (reference cell) and macroscale
equations, and exploits three main approximations:
\begin{itemize}
\item Because of the simple form of coupling between microscale and
  macroscale problems that involves only interfacial solution values, we
  propose a ``frozen coefficient'' splitting approach in
  which variables on the microscale are advanced to the next time step
  by holding all macroscale variables constant at their
  previous values, and vice versa.

\item The multiplier matrix $\Pi$ defined in~\eqref{cell_problem1} for
  the thermal diffusion coefficient in the macroscale heat equation is
  independent of the local temperature state and phase interface
  configuration.  Consequently, the entries $\Pi_{k\ell}$ are constants
  that only need to be computed once at the beginning of a simulation.

\item Both models have an inherent radial symmetry on the microscopic
  scale, and we restrict ourselves here to problems that have an
  analogous symmetry on the macroscale.  This is a natural choice for
  the tree sap exudation problem since a tree stem is well-approximated
  by a circular cylinder with cross-section $\Omega$ having radius
  $\Rtree$.  Consequently, all variables and governing equations are
  cast in terms of a single radial coordinate labelled $x$ or $y$ on the
  macro- or microscale respectively, so that only 1D problems need to be
  solved on both scales.
\end{itemize}
The spatial discretization of the governing equations is performed
separately on each spatial scale:
\begin{description}
\item[\emph{Macroscale heat problem:}] The circular domain $\Omega$ is
  discretized on an equally-spaced radial mesh of $\Mmacro=40$
  points, denoted $x_i = i\Rtree/\Mmacro$ for $i =
  0,1,\dots,\Mmacro-1$.  Discrete values of the unknowns $T_1(x,t)$ and
  $\HE_1(x,t)$ are defined at each mesh point $x_i$.

\item[\emph{Microscale heat problem:}] Within each reference cell $Y$,
  the portion of the domain ${Y}^2(x,t)$ consisting of ice will grow or
  shrink according to the location of the local phase boundary
  $s_\mathrm{iw}$. We therefore employ a \emph{moving mesh}
  discretization wherein the annular-shaped water region in the fiber is
  discretized at $\Mmicro$ equally-spaced radial points that move in
  time according to $y_j(t) = s_\mathrm{iw}(t) + j (\gamma -
  s_\mathrm{iw}(t))/\Mmicro$ for $j=0,1, \dots, \Mmicro$, where we
  recall that $\gamma$ is the radius of the artificial boundary
  $\Gamma$.  In practice, it suffices to use a coarse grid in the
  reference cell with $\Mmicro = 4 \ll \Mmacro$.  Discrete values of the
  solution variables $T_2(x,y,t)$ and $\HE_2(x,y,t)$ are then defined
  at each location $x_i$ and $y_j$.
\end{description}
Recall that the temperature $T_2$ is treated as the primary solution
variable in the microscale problem, whereas enthalpy $\HE_1$ is the
primary variable in the macroscale problem.  We employ a method-of-lines
approach in which spatial derivatives of solution quantities in both $x$
and $y$ are approximated using finite differences.  The resulting
coupled system of time-dependent DAEs is then integrated in time using
the stiff ODE solver {\tt ode15s} \citep{matlab-r2015a}.\
This solver requires absolute and relative error tolerances, which we
choose as {\tt AbsTol~=~7e-8} and {\tt RelTol~=~2e-14}.

We may then summarize the multiscale numerical algorithm
as follows:
\begin{description}
\item[\bf Step 1:] For a single canonical reference cell having the
  shape of a square with a circular hole, we use the package {COMSOL
    Multiphysics} \citep{comsolv5-2015}\ 
  to discretize the domain, approximate the functions $\mu_i(y)$ in
  \eqref{cell_problem1}, and then to calculate the corresponding
  integrals in \eqref{cell_problem2}.  This yields precomputed constant
  values of the four entries in matrix $\Pi$ that are used in the
  remainder of the computation (in Step 3c).

\item[\bf Step 2:] At each macroscale point $x_i$, set the initial value
  of $T_2=T_{2,\mathrm{init}}$.  Then within the $i$th reference cell,
  set $\HE_1=\HE_{1,\mathrm{init}}$ at each point $y_j$, and initialize
  either $s_\mathrm{iw}$ for the reduced model or $\left\{
    s_\mathrm{iw}, s_\mathrm{gi}, r, U \right\}$ for the sap exudation
  model.  Initial values are listed in Table~\ref{tab:ics1}.

\item[\bf Step 3:] At each time step, advance the solution variables as
  follows:
  \begin{description}
  \item[\bf 3a.] Set $T_1=\omega(\HE_1)$ and $\HE_2=\omega^{-1}(T_2)$.

  \item[\bf 3b.] Update $T_2$ by integrating the microscale heat
    diffusion problem
    \eqref{limit_problem_strong2a}--\eqref{limit_problem_strong2c} one
    time step within each reference cell ${Y}^2(x_i,t)$.  The values of
    $T_1$, $\HE_1$, $\HE_2$ and $s_\mathrm{iw}$ are frozen at the
    previous time step.

  \item[\bf 3c.] Update $\HE_1$ by integrating the macroscale heat
    diffusion problem
    \eqref{limit_problem_strong1a}--\eqref{limit_problem_strong1b} at
    all grid points $x_i$. Due to radial symmetry of the reference cell,
    the integral in the right hand side of~\eqref{limit_problem_strong1a}
    reduces to $2\pi R_2 D(\HE_2)\nabla T_2 \cdot \normvec$ where $R_2$
    is the radius of $Y^2$. The values of $T_2$ and $\HE_2$ are frozen
    at the values computed in step 3b.

  \item[\bf 3d.] Update the microscale variables within each reference
    cell $\widetilde{Y}^2(x_i,t)$ by integrating the governing
    differential(--algebraic) equations in time, and freezing values of
    $T_1$ and $T_2$.  Here, the equations being solved depend on the
    model problem:
    \begin{itemize}
    \item For the reduced problem, include the reduced Stefan condition
      \eqref{limit_problem_strong3} only. 
    \item For the sap exudation problem, use the system
      of DAEs \eqref{cell_problem_maple}--\eqref{algebraic}. 
    \end{itemize}

  \item[\bf 3e.] Increment the time variable and return to Step 3a. 
  \end{description}
\end{description}

\begin{table}[tbhp]
  \centering
  \caption{Initial values for the reduced and sap exudation models,
    taken from~\citet{AHORN5}.}   
  \label{tab:ics1}
  \begin{tabular}{clcc}\hline
    {Symbol} & {Description}   & {Initial Values} & {Units}\\\hline
    $T_\mathrm{init}$ & Initial temperature     & $\Tcrit$            & \myunit{\degK} \\
    $\Tout$    & Ambient temperature     & $\Tcrit+10$         & \myunit{\degK} \\
    $s_\mathrm{iw}(0)$& $= \Rf$                 & $3.5\times 10^{-6}$ & \myunit{m}\\
    $s_\mathrm{gi}(0)$& $= {\Rf}/{\sqrt2}$      & $2.5\times 10^{-6}$ & \myunit{m}\\
    $r(0)$     &                         & $6.0\times 10^{-6}$ & \myunit{m}\\
    $U(0)$     &                         & 0                   & \myunit{m^3}\\
    $p^\mathrm{f}_\mathrm{g}(0)$ &                         & $2.0\times 10^5$    & \myunit{N/m^2}\\
    $p^\mathrm{v}_\mathrm{g}(0)$ &                         & $1.0\times 10^5$    & \myunit{N/m^2}\\
    $p^\mathrm{f}_\mathrm{w}(0)$ &                         & $9.89\times 10^4$   & \myunit{N/m^2}\\
    $p^\mathrm{v}_\mathrm{w}(0)$ &                         & $9.95\times 10^4$   & \myunit{N/m^2}\\
    \hline
  \end{tabular}
\end{table}

The above algorithm must be modified slightly whenever the ice 
completely melts, since the loss of the Stefan condition
\eqref{limit_problem_strong3} induces a change in the governing
equations.  At the same time, the separation of the reference cell into
two sub-domains $Y^1$ and $Y^2$ (which was required to handle the
Dirichlet condition on temperature at the phase interface) is no longer
necessary and hence the temperature can be described by the single field
$T_1$ that obeys
\eqref{limit_problem_strong1a}--\eqref{limit_problem_strong1b} with
zero right hand side, constant $D$, and $\Pi\equiv 1$.  This alteration
to the governing equations can be triggered easily within the numerical
algorithm above by exploiting the ``event detection'' feature in
Matlab's {\tt ode15s} solver, signalling an event based on a
zero-crossing of the ice layer thickness,
$b=s_\mathrm{iw}-s_\mathrm{gw}$: when $b>0$, ice is still present and
the original equations are solved; when $b=0$, ice is totally melted and
the modified equations just described are solved instead (and Steps 3b
and 3d are omitted).  We mention in closing that although only 2D
simulations are performed in this paper, our algorithm extends in a
straightforward manner to 3D stem geometries by stacking a number of 2D
stem slices in series and enforcing suitable flux continuity conditions.

\subsection{Simulations of the reduced model}
\label{sec:sims-reduced}

We begin by presenting numerical simulations of the reduced model
wherein a periodic array of melting ice bars fills a circular domain
$\Omega$ with radius $\Rtree=0.25\;\myunit{m}$.  The periodic reference
cell $Y=[0,\delta]^2$ depicted in Figure~\ref{fig:ice-refcell}b is given
a side length of $\delta=4.33\times 10^{-5}\;\myunit{m}$.  Each
reference cell is initialized with an ice bar of radius
$s_\mathrm{iw}(0)= \Rf / \sqrt{2}$ surrounded by water, such that the
initial volume of ice in the reduced model and the sap exudation model
is equal.  The initial temperature throughout the domain is set to
$T_1(x,0)=T_2(x,y,0)=T_\mathrm{init}=\Tcrit$.  On the outer boundary of
the domain, a Robin boundary condition $-D(E_1(\Rtree,t))\nabla_x
T_1(\Rtree,t)\cdot \normvec = \alpha(T_1(\Rtree,t) - \Tout)$ is imposed
with $\Tout = \Tcrit + 10$, while a symmetry condition $\partial_x
T_1(0,t)=0$ is imposed at the center of the domain.  We take the size of
the artificial boundary $\Gamma$ in each reference cell to be larger
than the fiber radius $\Rf$ by an amount equal to the typical thickness
$W=4.38\times 10^{-6}\;\myunit{m}$ of the vessel wall; in other words,
$\gamma=\Rf+W$ which is well-separated from the phase interface.  Note
that the system is solved in dimensional variables so that there is no
need to non-dimensionalize and hence the size of the reference cell
corresponds simply to the physical dimension $\delta$.  All physical
parameter values and initial conditions are listed in
Tables~\ref{tab:params2} and~\ref{tab:ics1}.

\begin{figure}[tbhp]
  \centering\footnotesize
  \begin{tabular}{ccc}
    (a) Macroscale temperature $T_1(x,t)$ && 
    (b) Local ice-bar radius $s_\mathrm{iw}(x,t)$\\
    \includegraphics[width=0.40\textwidth]{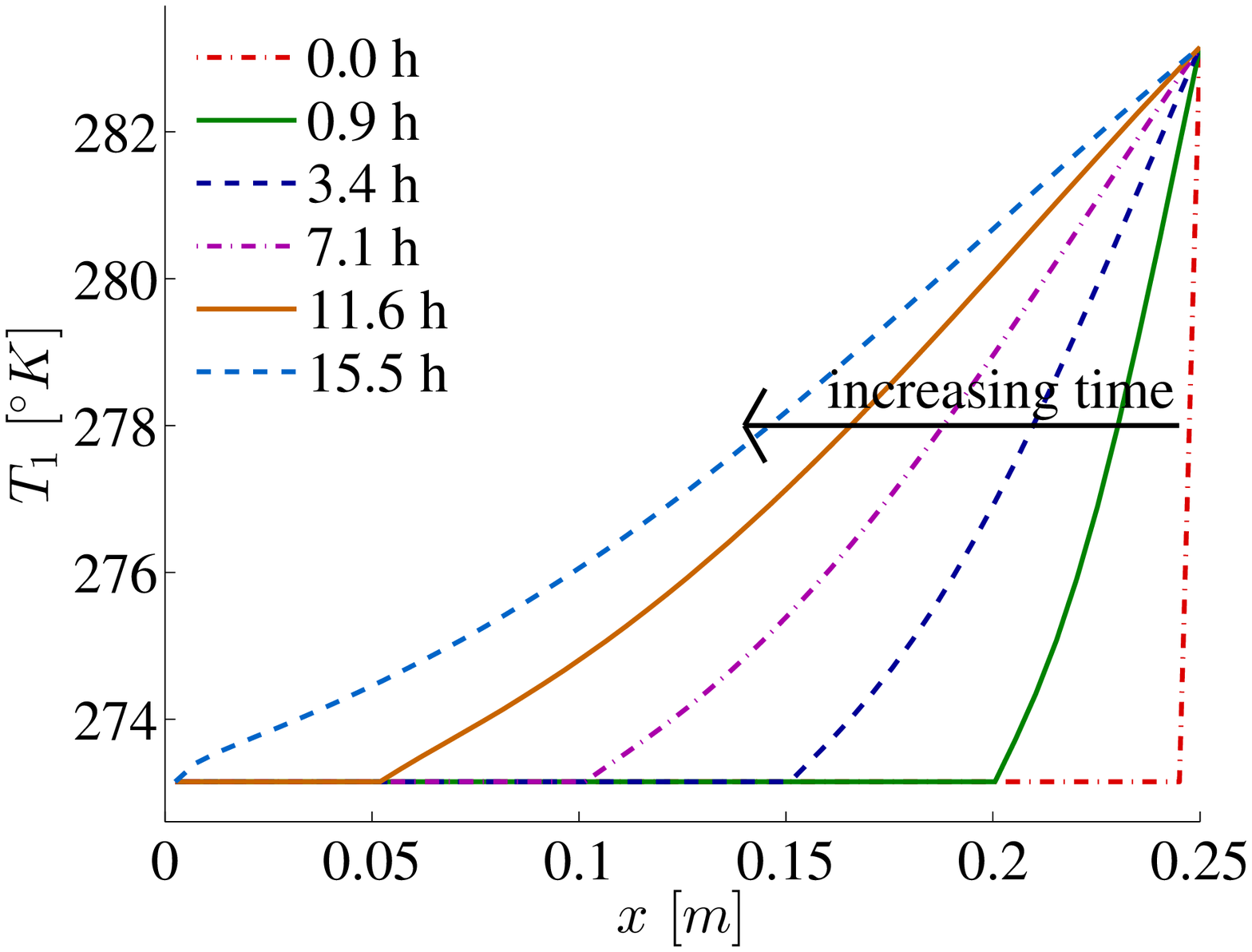}
    &&
    \includegraphics[width=0.40\textwidth]{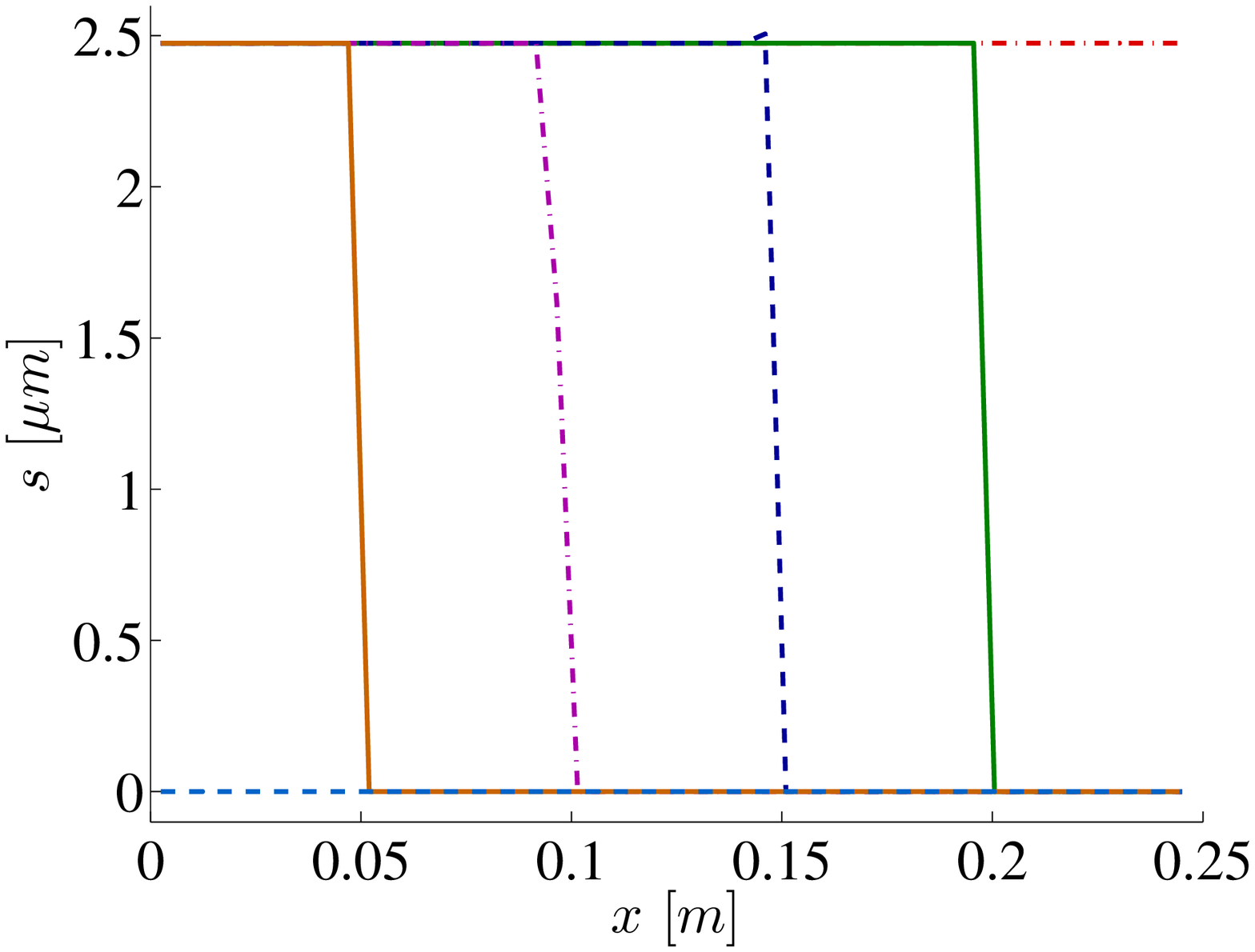}
  \end{tabular}
  \caption{Simulations of the reduced model, showing $T_1$ and $s_\mathrm{iw}$
    as functions of the global variable $x$, pictured at selected times
    between 0 and 16~\myunit{h}.}
  \label{fig:sims-reduced}
\end{figure}

Figure~\ref{fig:sims-reduced} displays a sequence of solution snapshots
at selected times between 0 and 16~h that illustrate the spatial and
temporal variations in the macroscale temperature $T_1$ and ice-bar
radius $s_\mathrm{iw}$. In each plot, the horizontal ($x$) axis corresponds to
the radial distance measured from the center of the circular domain
$\Omega$.  As time progresses, the temperature gradually increases and
penetrates the domain interior as heat from the outer boundary diffuses
inwards.  In response to this rise in temperature, the ice melts and the
ice bar within each local reference cell shrinks in size.  The ice bars
in the outermost region melt first, and by time $t\approx 16~\myunit{h}$
the entire domain is completely melted (i.e., $s_\mathrm{iw}=0$ throughout
$\Omega$).  The formation of a steep thawing front that progresses from
the outer boundary to the center of the domain is clearly visible in
Figure~\ref{fig:sims-reduced}b.  These results should be contrasted with
the study in~\citet{AHORN5} that investigated only the local behaviour of
the solution to the thawing model (at some fixed location on the
microscale); on the other hand, our homogenized model results illustrate
the progress of the thawing front on the macroscale, while at the same
time incorporating physical processes taking place on the microscale.

\subsection{Simulations of the sap exudation model}
\label{sec:sims-sap}

Next we perform simulations of the sap exudation model for a tree stem
with the same parameters and boundary conditions as for the reduced
model.  Recall that we chose the length of the reference cell to be
$\delta=4.33\times 10^{-5}\; \myunit{m}$, consistent with the size of
fibers and vessels in actual sapwood.  The initial temperature is again
taken to be $T_\mathrm{init}=\Tcrit$ throughout, with the water in the fiber
initially frozen and the vessel sap in liquid form.  Recall that this
initial state captures the effect of freezing point depression
due to the presence of sugar within the vessel sap; and besides setting
these initial conditions, there is no need to incorporate any
concentration dependence in the freezing point for this thawing-only
model.  To initiate a thawing cycle, we apply the Robin boundary
condition at the outer boundary of the tree stem as in the reduced
model, $-D(E_1(\Rtree,t))\nabla_x T_1(\Rtree,t)\cdot \normvec =
\alpha(T_1(\Rtree,t) - \Tout)$, where $\Tout = \Tcrit + 10$.  All other
parameters and initial values specific to the sap exudation model can be
found in Tables~\ref{tab:params2} and~\ref{tab:ics1}.

\begin{figure}[tbhp]
  \centering\footnotesize
  \begin{tabular}{ccc}
    (a) Temperature $T_1$ & \qquad & 
    (b) Phase interfaces $s_\mathrm{iw}$, $s_\mathrm{gi}$\\
    \includegraphics[width=0.40\textwidth]{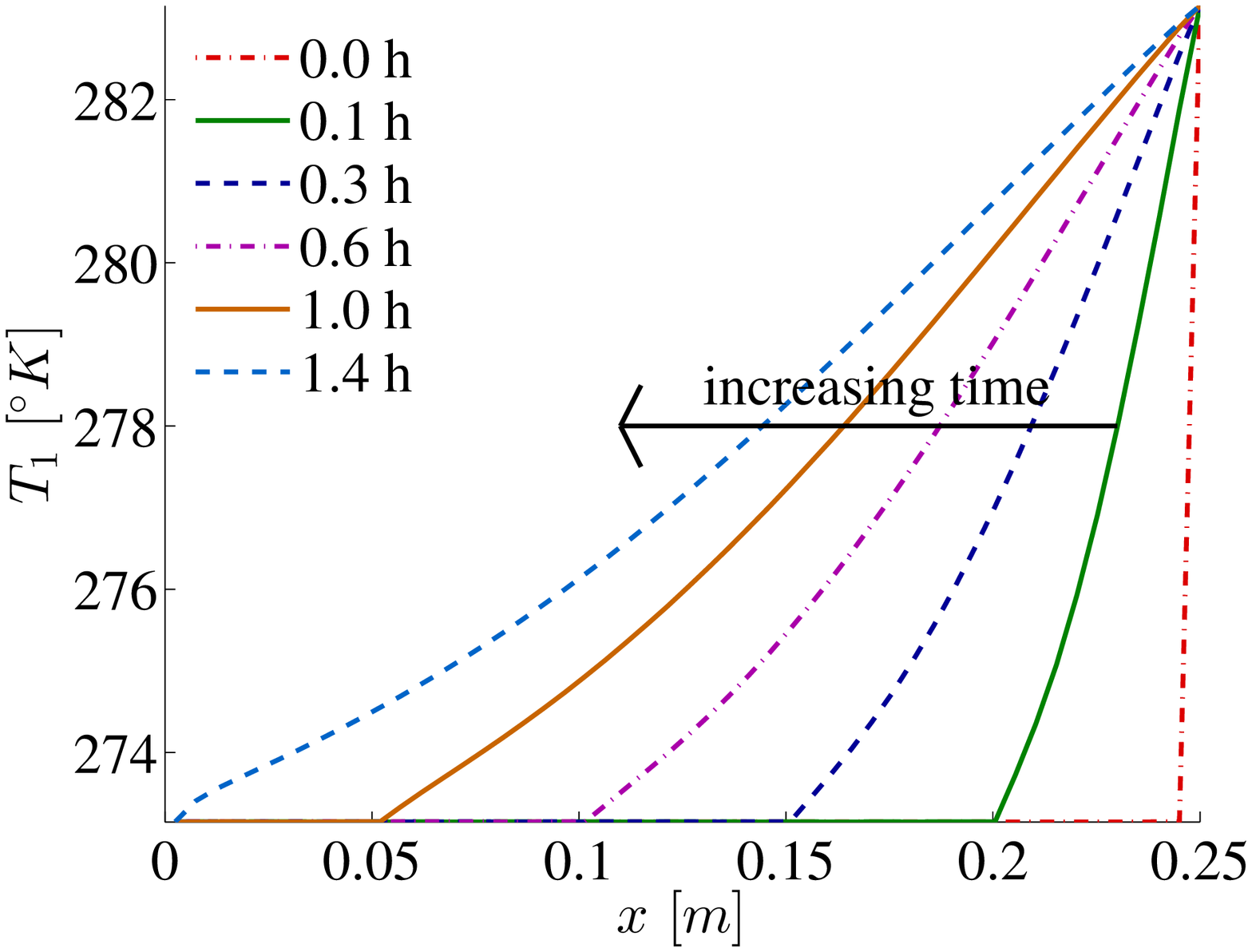}
    & &
    \includegraphics[width=0.40\textwidth]{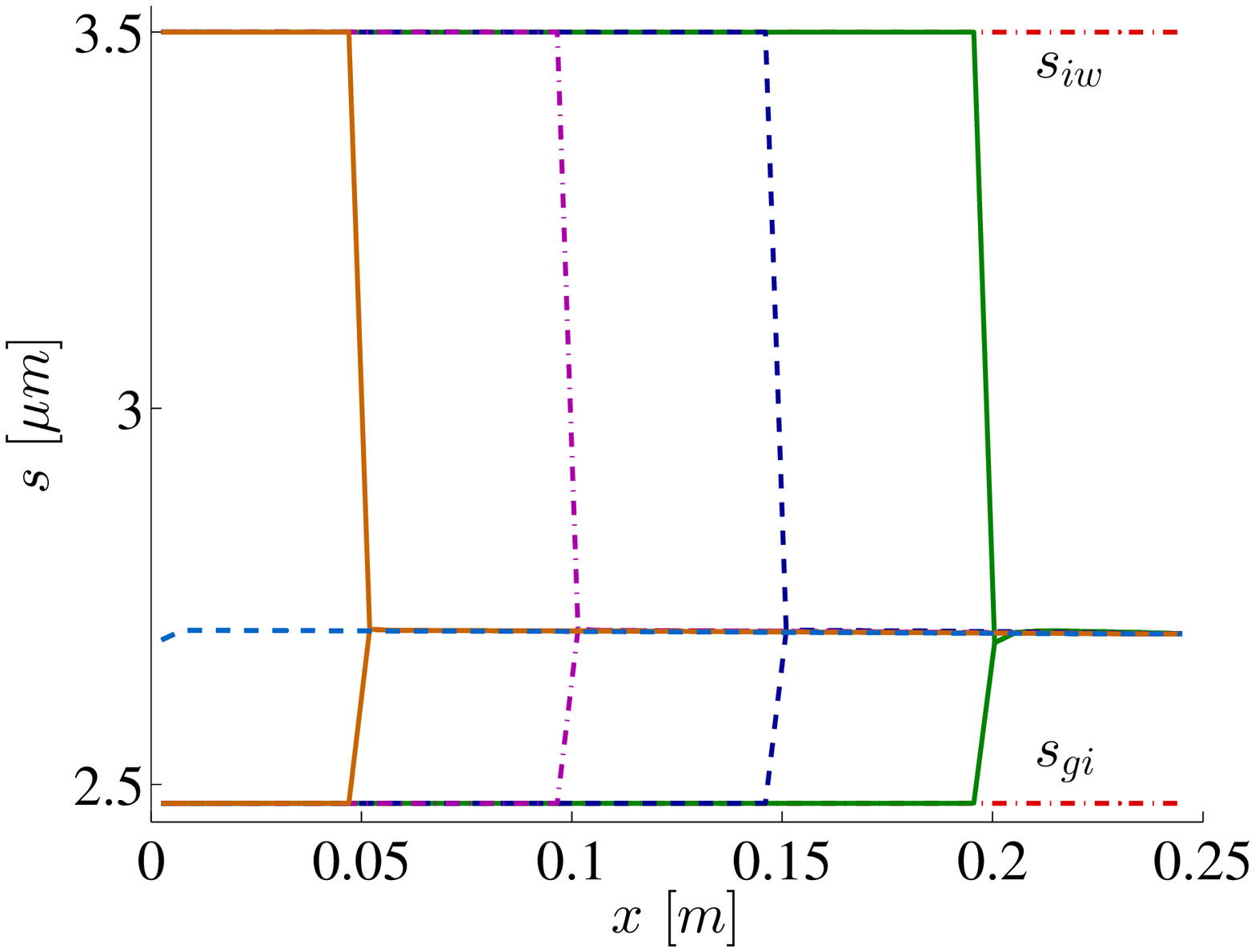}
    \\
    (c) Melt-water volume $U$ & & 
    (d) Vessel bubble radius $r$\\
    \includegraphics[width=0.40\textwidth]{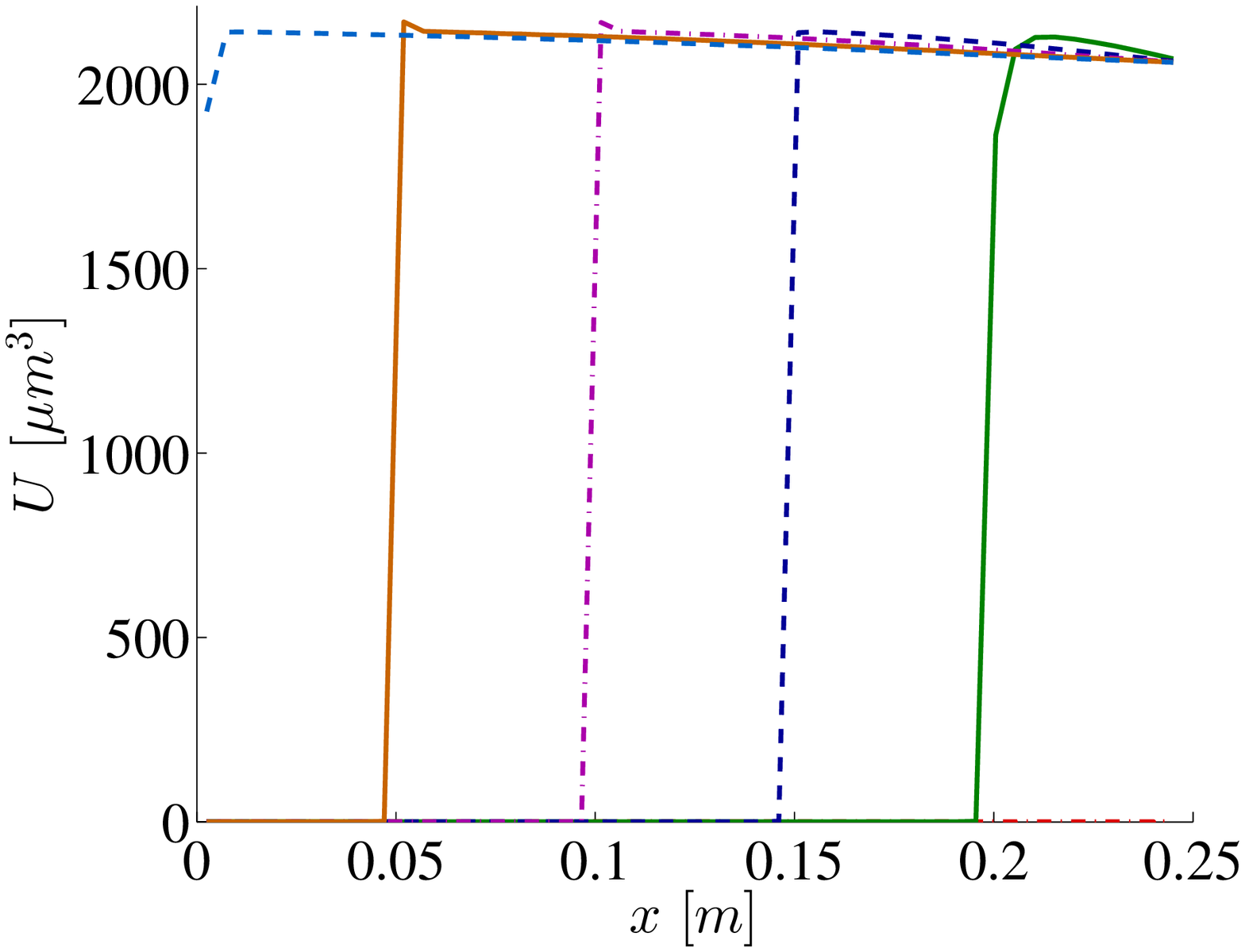}
    & &
    \includegraphics[width=0.40\textwidth]{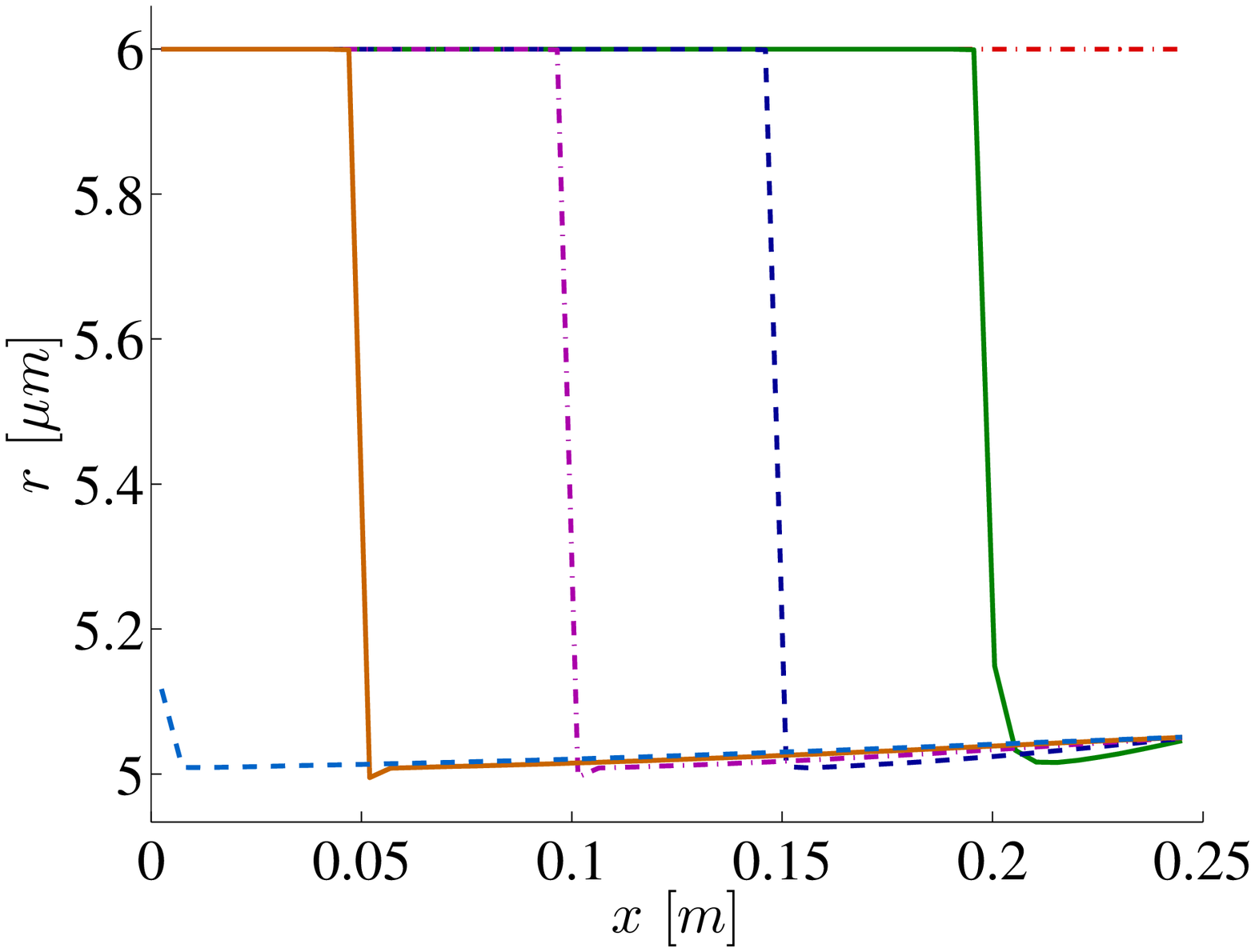}
    \\
    (e) Fiber water pressure $p^\mathrm{f}_\mathrm{w}$ & & 
    (f) Vessel sap pressure $p^\mathrm{v}_\mathrm{w}$\\
    \includegraphics[width=0.40\textwidth]{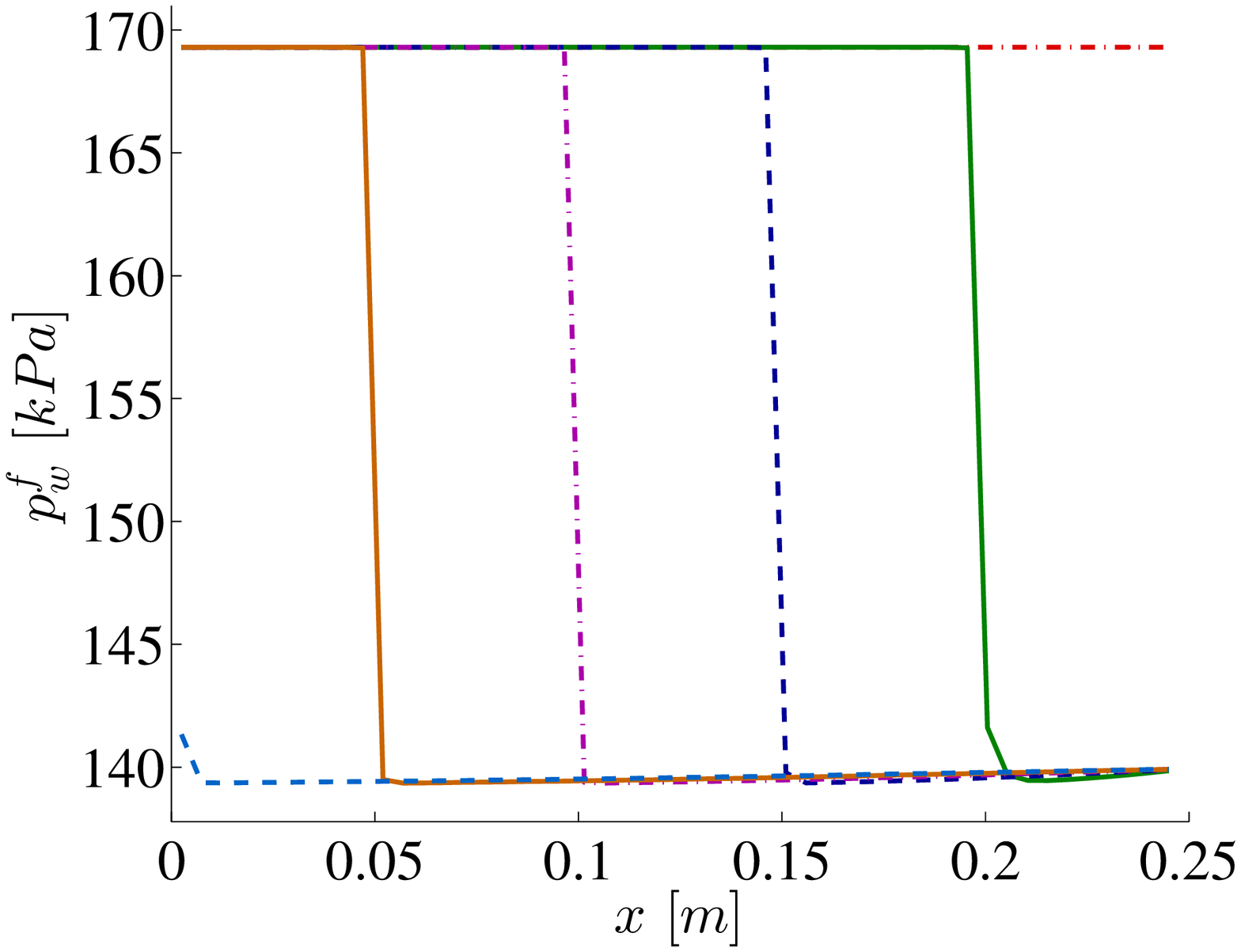}
    & &
    \includegraphics[width=0.40\textwidth]{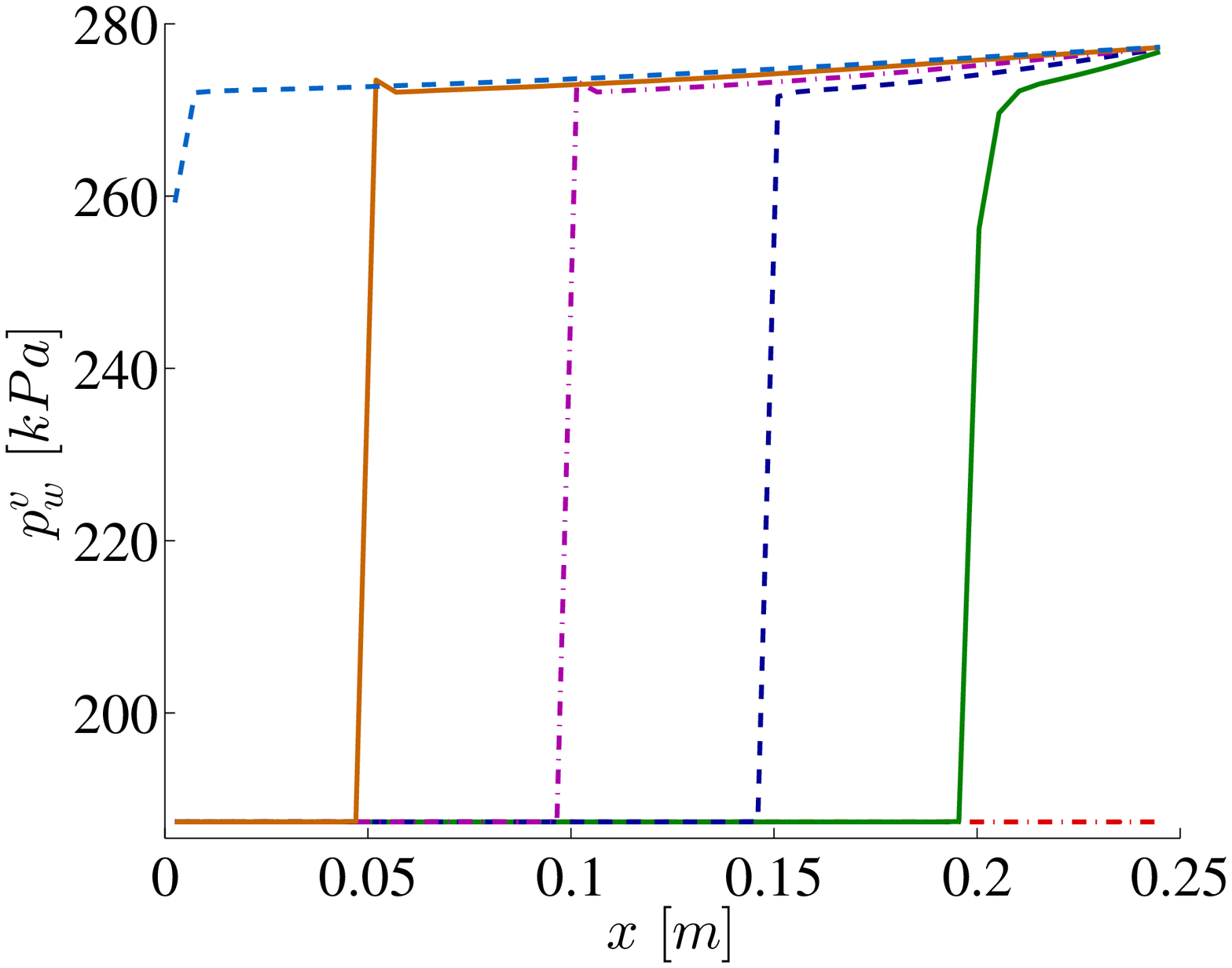}
  \end{tabular}
  \caption{Simulations of the sap exudation model showing
    $(x,t)$-dependent solution profiles at a sequence of time points.
    In all cases, the profiles evolve from right to left (toward the
    center of the tree) as indicated by the arrow in (a).}
  \label{fig:sap1}
\end{figure}

After applying the multiscale algorithm described in
Section~\ref{sec:algo}, the solutions for $T_1$, $s_\mathrm{iw}$ and
$s_\mathrm{gi}$, $r$, $U$, $p_\mathrm{w}^\mathrm{f}$ and
$p_\mathrm{w}^\mathrm{v}$ are illustrated in Figure~\ref{fig:sap1} at a
sequence of six times between 0 and $1.4\;\myunit{h}$.
From these plots, it is evident that the solution dynamics for all
variables are characterized by a \emph{melting front} that progresses
through the tree from the outer boundary toward the center (from right
to left in the plots) as the warm ambient air gradually heats up the
interior.  Furthermore, the time required for complete melting of the ice
contained in the fibers is just under 1.4\;\myunit{h}.

The temperature profiles vary smoothly in space as one expects from a
diffusion problem, while the other solution quantities are characterized
by a steep front that propagates toward the centre of the tree with
a speed that decreases with time.  The steepness of the melting front
derives from the thawing of ice and subsequent adjustment of liquid
between vessels and fibers on the microscale, all of which occur very
rapidly in the instant after the temperature exceeds the freezing point
$\Tcrit$ at any given location $x$.  The reason for the gradual slowing
of the melting front with time is that the heat flux naturally decreases
as the front approaches the center of the tree, which in turn leads
to a speed decrease owing to the Stefan condition.

\begin{figure}[tbhp]
  \centering\footnotesize
  \begin{tabular}{ccc}
    (a) Temperature $T_1$ & 
    (b) Phase interfaces $s_\mathrm{iw}$, $s_\mathrm{gi}$ &
    (c) Liquid pressures $p^\mathrm{f}_\mathrm{w}$, $p^\mathrm{v}_\mathrm{w}$ \\
    \includegraphics[width=0.30\textwidth]{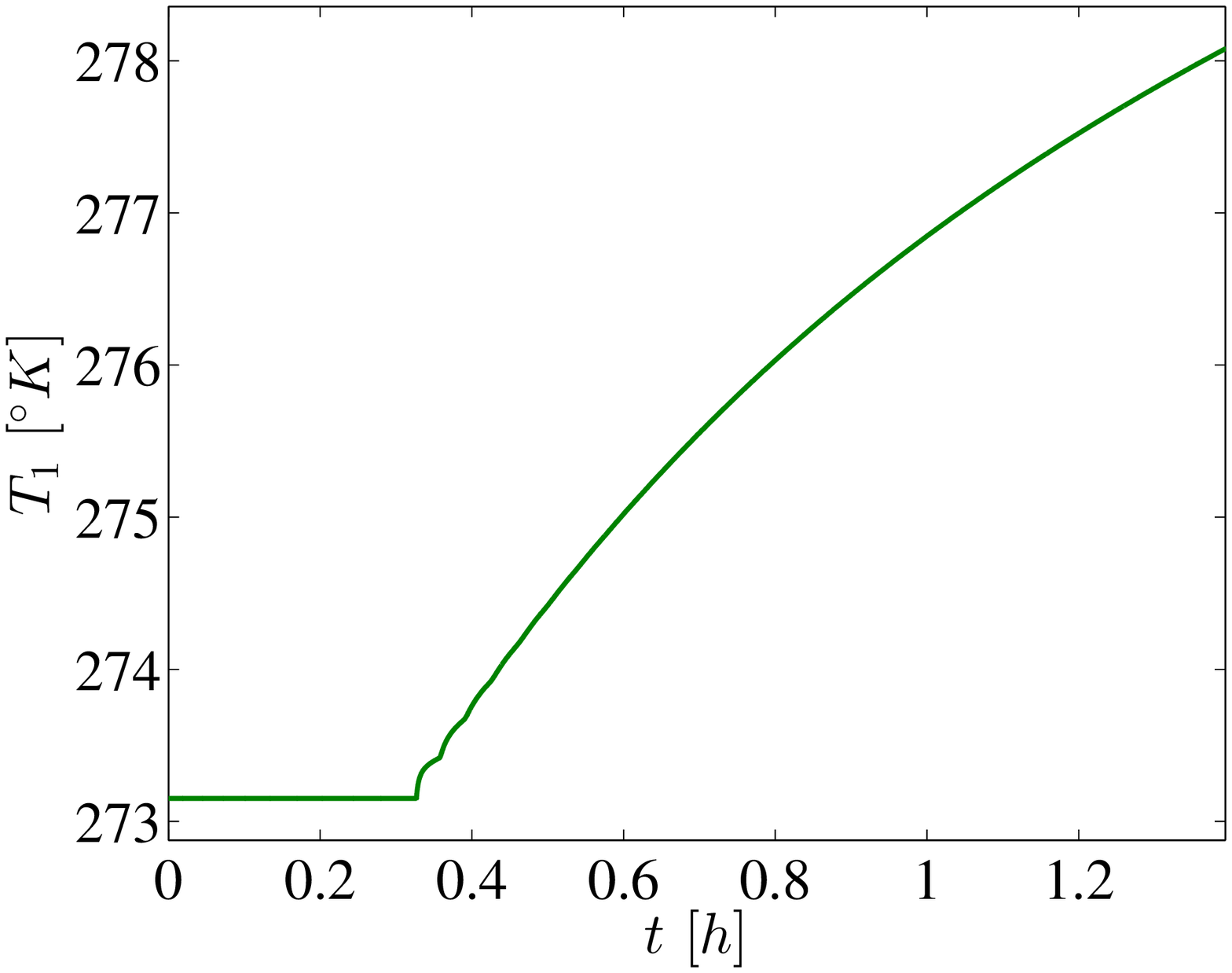} & 
    \includegraphics[width=0.30\textwidth]{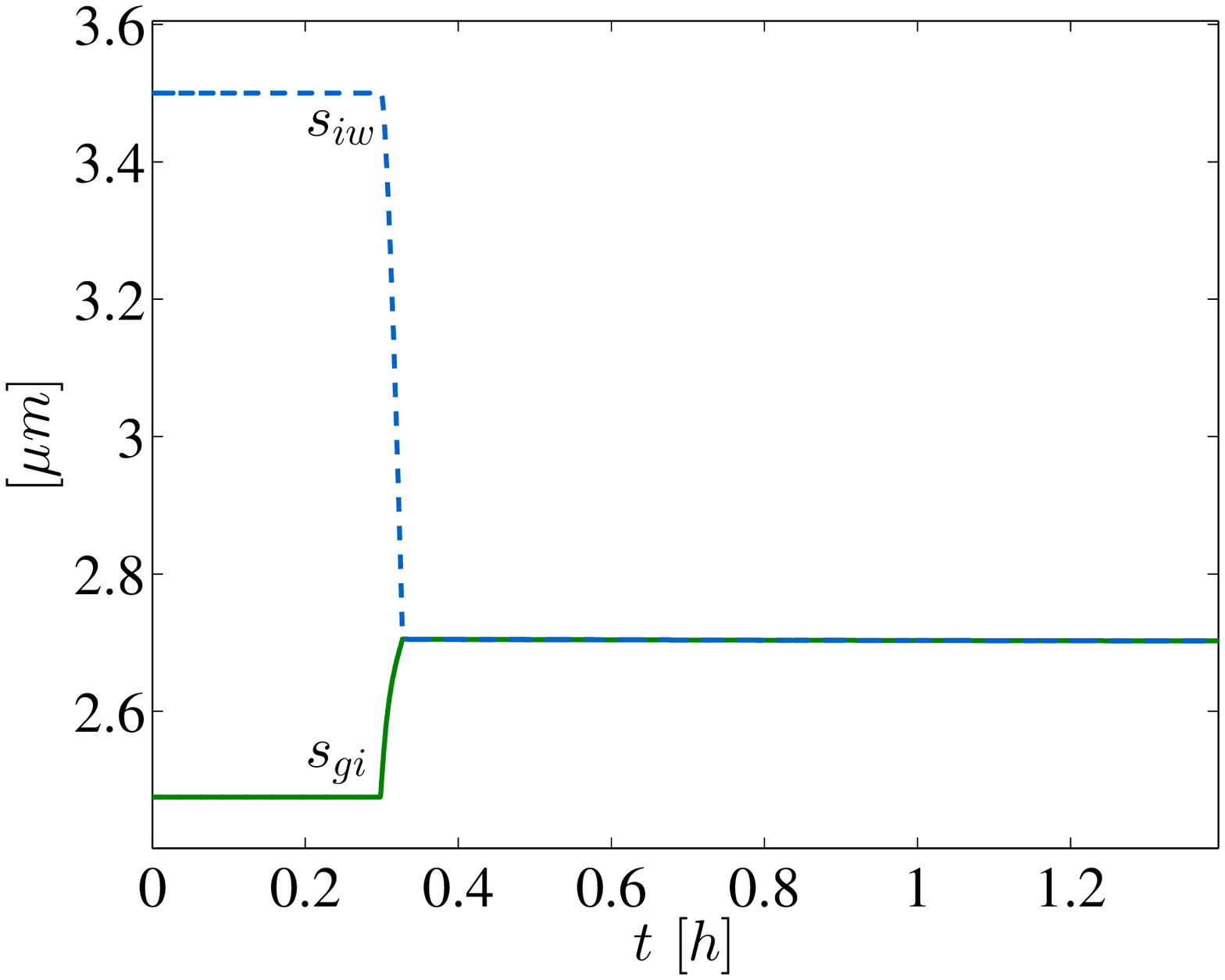} &
    \includegraphics[width=0.30\textwidth]{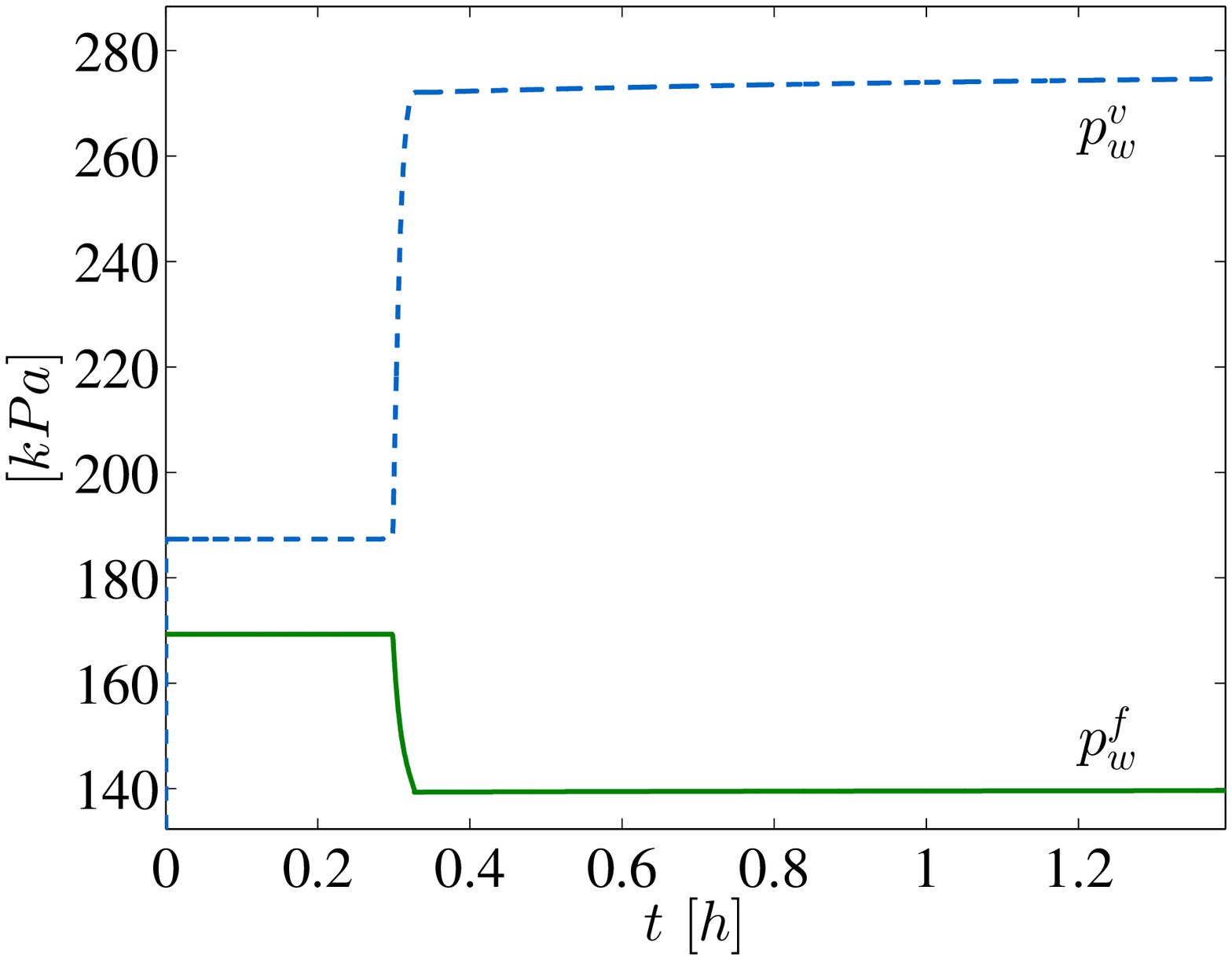} 
  \end{tabular}
  \caption{Simulations of the sap exudation model, showing the time
    evolution of various solution components at the fixed radial
    position $x=0.15\;\mathrm{m}$.} 
  \label{fig:sap2}
\end{figure}

We also observe a clear separation in time scales between the slow
evolution of temperature on the macroscale and the relatively rapid
phase change and sap redistribution within fibers and vessels on the
microscale.  This scale separation is easily seen by comparing
Figure~\ref{fig:sap1} with plots of the time evolution of local solution
variables at a fixed radial location $x=0.15\;\myunit{m}$ shown in
Figure~\ref{fig:sap2}.  The thickness of the fiber--ice layer can be
determined as the vertical distance between the $s_\mathrm{iw}$ and
$s_\mathrm{gi}$ curves in Figure~\ref{fig:sap2}b, which rapidly drops to
zero as the ice melts.  At the same time, melt-water is driven from
fiber to vessel by the pressure stored in the fiber gas bubble, and the
pressure plot in Figure~\ref{fig:sap2}c clearly illustrates the
subsequent increase in $p^\mathrm{v}_\mathrm{w}$ that we attribute to
{exudation pressure}.  After the melting process is complete, the
vessel--liquid pressure continues to increase (although at a slow rate
that is not easily visible to the naked eye) owing to a slight expansion
of gas in the fiber and vessel in response to further temperature
increase as heat continues to diffuse through the liquid phase from the
outer tree surface.

Upon more careful inspection of the gas/ice and ice/water interfaces in
Figure~\ref{fig:sap2}b, we observe that there is a slight time delay in
the motion of $s_\mathrm{iw}$ relative to $s_\mathrm{gi}$.  Indeed, the ice begins to
melt at the gas/ice interface (leading to an increase in $s_\mathrm{gi}$) at a
time that is roughly $25\;\myunit{s}$ in advance of when $s_{wi}$ starts
to drop, which is when a water layer appears between the ice and fiber
wall.  This phenomenon can be explained as follows.  When melt-water
first appears in a particular fiber, the gas bubble pressure is so high
that water is immediately forced out into the vessel, leaving the ice
layer in contact with the fiber wall.  The gas pressure then declines
until approximately $25\;\myunit{s}$ elapses, at which time the rate of
water melting exceeds that of the porous outflow and a water layer
begins to accumulate along the fiber wall.  By this time, roughly half
of the water volume contained in the fiber has been transferred into the
vessel.

One of the most significant results from our sap exudation model is the
prediction that vessel liquid pressure increases by roughly
$120\,\myunit{kPa}$, which is within the range of exudation pressures
actually observed in sugar maple trees \citep{AHORN4} and closely-related
species such as black walnut \citep{AHORN3}.  Furthermore, simulations of
multiple freeze and thaw cycles with this homogenized model yield
results that are consistent with controlled experiments on
walnut \citep{graf-ceseri-stockie-2015}.  Work is currently underway on
comparing our model results to experiment measurements of sugar maple
saplings \citep{brown-2015}.

Finally, we draw a comparison between the solution of the sap exudation
model for temperature $T_1$ and ice layer thickness
$s_\mathrm{iw}-s_\mathrm{gi}$, and the corresponding solution variables
($T_1$, $s_\mathrm{iw}$) from the reduced problem in
Section~\ref{sec:problem1}.  Although the shape of the temperature and
ice interface profiles are similar, there is a significant difference in
that the melting process for the reduced problem takes over 10 times
longer than for the sap exudation problem even though the macroscopic
domain and outer temperature are the same.  This discrepancy may seem at
first glance to be inconsistent, but this result can be easily explained
as follows. The diffusion coefficient for the sap exudation problem is
roughly 10 times higher because of the much larger value of thermal
diffusivity ($k/\rho c$) in the gas phase (with an upper bound of
$2\times 10^{-5}\;\myunit{m^2/s}$ based on atmospheric conditions)
compared with the corresponding values for ice and water
($1.2\times10^{-6}$ and $1.3\times 10^{-7}\;\myunit{m^2/s}$
respectively) which are the only phases appearing in the reduced model.

\section{Conclusions}
\label{sec:conclusions}

The aim of this paper was to apply techniques from periodic
homogenization to derive a multiscale model for a multiphase flow
problem arising in the context of maple sap exudation.  Because of the
complexity of the physics underlying the sap exudation problem
(involving liquid/ice phase change, dissolving gas, flow through
porous cell membranes, osmosis and other effects), we started by
deriving a simpler reduced model that focuses on the melting of a
periodic array of ice bars.  This reduced model belongs to the class
of Stefan problems, which have been well-studied in the context of
homogenization in case of fast diffusion.  We prove results on
existence, uniqueness and a~priori estimates for the weak form of the
reduced governing equations involving fast \emph{and} slow diffusion,
which we then use to derive a strong form of the homogenized limit
problem in which there is a clear separation between equations for the
cellular level processes on the microscale, and heat transport on the
macroscale.  Our approach has the advantage that it applies
homogenization techniques in a straightforward manner in order to
obtain an uncomplicated limit model.  The primary novelty of the
analytical results, relative to other work on homogenization of
Stefan-type problems, derives from our directly imposing a Dirichlet
condition on temperature at the phase interface, which gives rise to a
decomposition into fast and slow variables on the sub-regions $Y^1$
and $Y^2$ of the reference cell.  A major advantage of this
decomposition is that it leads immediately to a simple and efficient
numerical method based on a time-splitting approach that exploits the
scale separation in the limit equations.  In particular, we are able
to encapsulate all microscale processes specific to the phase change
within the reference cell domain $Y^2$, wherein the temperature
diffuses slowly.  Consequently, this homogenized limit structure (and
the corresponding numerical algorithm) can be easily adapted to the
sap exudation problem by simply ``plugging in'' the corresponding
microscale equations governing the cell-level processes. In passing,
we proved a general existence result for quasi-linear parabolic
differential equations having a non-monotone nonlinearity in the
diffusion operator. Numerical simulations are performed for both the
sap exudation problem and the reduced model, and the results are shown
to be consistent, although there are significant differences that we
attribute to the absence of a gas phase in the reduced model.  The
homogenized limit equations derived here have been extended
elsewhere \citep{graf-ceseri-stockie-2015} to handle the freezing case,
and then applied to simulate multiple daily cycles of freeze and thaw
cycles; these results show an excellent match with sap exudation
experiments.

There are several natural avenues for future work that arise from this
study.  Most notably, we would like to fill the gap in our analytical
results by extend the proofs to handle the nonlinearity arising from an
enthalpy-dependent thermal diffusion coefficient.  This effort will be
guided by results on other related nonlinear
problems \citep{kanschatkrebs-mscthesis-2015,STEFAN2}.  We also plan to
extend our model to handle the three-dimensional geometry of a
cylindrical tree stem and include radial flow of sap and the effect of
gravitational pressure head on vertical transport.


\appendix
\setcounter{theorem}{0} 
\setcounter{lemma}{1}   

\section{Proofs of three main results}
\label{appendix}

This appendix contains proofs of the lemmas and theorems introduced in
Section~\ref{sec:analysis}.  Throughout, we use $\bigC$ or $C_i$ to
denote a generic, real, positive constant whose value may change from
line to line.

\subsection{Proof of existence}
\label{app:existence}

\newcounter{tmp}
\begingroup
\setcounter{tmp}{\value{theorem}}
\renewcommand\thetheorem{\ref{sec:analysis}.\arabic{theorem}}
\begin{theorem}
  \theoremone
\end{theorem}
\endgroup

\begin{proof}
  \begin{enumerate}[{(i)}]
  \item To handle the nonlinearities of $A$ we perform the
    semi-discretization
    \begin{align}\label{ext_problem3}
      \frac{u^m - u^{m-1}}{k} + A(mk)(u^{m-1},u^m) = f^m
      \qquad\text{in}\ V^*, 
    \end{align}   
    for $m=1,\dots, N$, with $N\in\N$, $N>2$, $k=\frac{\tend}{N}$,
    $u^0=u_0$ and
    \begin{gather*}
      f^m = \frac{1}{k} \int_{J_m} f(t)\, \text{d}t,
    \end{gather*}
    for $m=1,\dots, N$ and $J_m = [(m-1)k, mk]$.  Then $f_k(t) = f^m$ on
    $J_m$ for $m=1,\dots, N$. With these conditions the existence of the
    semi-discrete problem \eqref{ext_problem3} holds.
    
  \item (A~priori estimates) We define the functions
    $u_k:[0,T]\rightarrow V$ and $w_k:[0,T]\rightarrow L^2(\Omega)$ with
    \begin{subequations}\label{ext_def1}
      \begin{align}\label{ext_def2}
        u_k(t) &= u^m,\\
        w_k(t) &= u^m + \frac{t-mk}{k}(u^{m+1} - u^m).
      \end{align}
    \end{subequations}
    Then it holds that
    \begin{align}\label{ext_eq1}
      \norm{u_k - w_k}^2_{L^2(S,L^2(\Omega))} \leq \frac{k}3
      \sum_{m=1}^N\norm{u^m - u^{m-1}}^2_{L^2(\Omega)}. 
    \end{align}
    For the next estimate we start with the fact that
    \begin{gather*}
      2(a-b\, , \,a) = \norm{a}^2_\Omega - \norm{b}^2_\Omega +
      \norm{a-b}^2_\Omega \quad \forall\ a,b \in L^2(\Omega),
    \end{gather*}
    which implies
    \begin{multline}
      \norm{u^m}^2_\Omega - \norm{u^{m-1}}^2_\Omega 
      + \norm{u^m - u^{m-1}}^2_\Omega 
      + 2k\langle A(mk)(u^{m-1},u^{m}), u^m\rangle_{V^*V} 
      = 2k\langle f^m,u^m\rangle_{V^*V}.
    \end{multline}
    Using H\"older's inequality we obtain that
    \begin{align}\label{ext_eq2}
      \norm{u^m}^2_\Omega - \norm{u^{m-1}}^2_\Omega + \norm{u^m -
        u^{m-1}}^2_\Omega + k\lambda\norm{u^m}^2_V \leq
      \frac{k}{\lambda}\norm{f^m}^2_{V^*}. 
    \end{align}
    Making use of the inequality
    \begin{align}\label{ext_eq3}
      k\sum_{m=1}^N\norm{f^m}^2_{V^*}\leq \int_0^T\norm{f(t)}^2_{V^*}\,
      \text{d}t , 
    \end{align}
    and then summing \eqref{ext_eq2} over $m$ from 1 to $r\leq N$ yields
    \begin{align}\label{ext_eq4}
      \norm{u^r}^2_{\Omega} + \sum_{m=1}^r \norm{u^m - u^{m-1}}^2_\Omega
      + k\lambda\sum_{m=1}^r\norm{u^m}^2_V \leq \norm{u_0}^2+
      \frac1\lambda\int_0^T\norm{f(t)}^2_{V^*}\, \text{d}t.
    \end{align}
    It then follows that 
    \begin{subequations}\label{ext_eq5:abc}
      \begin{align}
        \norm{u^r}^2_\Omega &\leq \norm{u_0}^2_\Omega +
        \frac1\lambda\int_0^T\norm{f(t)}^2_{V^*}\, \text{d}t
        \quad\text{for } 1\leq r \leq N, \label{ext_eq5:a} \\
        \sum_{m=1}^r\norm{u^m - u^{m-1}}^2_\Omega &\leq
        \norm{u_0}^2_\Omega + \frac1\lambda\int_0^T\norm{f(t)}^2_{V^*}\,
        \text{d}t,  \label{ext_eq5:b}\\
        k\sum_{m=1}^r \norm{u^m}^2_V &\leq \frac{1}{\lambda}
        \norm{u_0}^2_\Omega + \frac{1}{\lambda^2}
        \int_0^T\norm{f(t)}^2_{V^*}\, \text{d}t. \label{ext_eq5:c}
      \end{align}
    \end{subequations}
    Using \eqref{ext_eq5:a} we obtain the estimates
    \begin{subequations}\label{ext_eq6:ab}
      \begin{align}
        \norm{u_k}_{L^\infty(S,L^2(\Omega))} \leq c, \label{ext_eq6:a}\\
        \norm{w_k}_{L^\infty(S,L^2(\Omega))} \leq c, \label{ext_eq6:b}
      \end{align}
    \end{subequations}
    where the constant $c$ depends only on the right-hand side of
    \eqref{ext_eq5:a}. Because of
    \begin{align*}
      \int_0^T\norm{u_k(t)}^2_V\, \text{d}t =
      k\sum_{m=1}^N\norm{u^m}^2_V,  
    \end{align*}
    it follows that 
    \begin{align}\label{ext_eq7}
      \norm{u_k}_{L^2(S,V)} \leq c.
    \end{align}
    The Banach space $L^\infty(S,L^2(\Omega))$ is the dual space of the
    separable space $L^1(S,L^2(\Omega))$; hence, when taken together
    with the estimates \eqref{ext_eq6:ab} and \eqref{ext_eq7} and the
    theorems of Eberlein--Shmuljan and Banach--Alaoglu, we are
    guaranteed the existence of subsequences
    \begin{subequations}\label{ext_eq8:abc}
      \begin{align}
        u_k &\rightharpoonup u\quad \text{in}\ L^2(S,V), \label{ext_eq8:a}\\
        u_k &\rightharpoonup^* u \quad \text{in}\ L^\infty(S, L^2(\Omega)), \label{ext_eq8:b}\\
        w_k &\rightharpoonup^* w \quad \text{in}\ L^\infty(S, L^2(\Omega)). \label{ext_eq8:c} 
      \end{align}
    \end{subequations}
    We next want to show that 
    \begin{align}\label{ext_eq9}
      u=w.
    \end{align}
    Using \eqref{ext_eq1} and \eqref{ext_eq5:a} we conclude that 
    \begin{align}\label{ext_eq9a}
      u_k - w_k \rightarrow 0 \quad\text{in}\ L^2(S,L^2(\Omega)), 
    \end{align}
    so that for a subsequence 
    \begin{gather*}
      u_k(t) - w_k(t) \rightarrow 0.
    \end{gather*}
    With Eqs.~\eqref{ext_eq6:ab} we obtain that
    \begin{gather*}
      u_k - w_k \rightarrow 0 \quad\text{in}\ L^\infty(S,L^2(\Omega)), 
    \end{gather*}
    and hence we can rewrite \eqref{ext_eq8:c} as
    \begin{gather*}
      w_k \rightharpoonup^* u \quad\text{in}\ L^\infty(S, L^2(\Omega)).
    \end{gather*}
    The nonlinearity of $A$ requires another a~priori estimate to
    perform the limit  $k\rightarrow 0$.  From \eqref{ext_eq1},
    \eqref{ext_cond:e} and \eqref{ext_cond2:a} we deduce that
    \begin{align}
      \Norm{\frac{u^m - u^{m-1}}{k}}_{V^*} \leq \norm{f^m}_{V^*} +
      \Lambda\norm{u^m}_V. 
    \end{align}
    After that, we apply \eqref{ext_eq3} and \eqref{ext_eq5:c} to obtain
    \begin{align}\label{ext_eq10}
      k\sum_{m=1}^N \Norm{\frac{u^m - u^{m-1}}{k}}^2_{V^*} \leq d <
      \infty, 
    \end{align}
    where $d$ only depends on the data in \eqref{ext_cond:e} and
    constants $\lambda$ and $\Lambda$. This estimate implies that
    \begin{align}\label{ext_eq11}
      \norm{w'_k}_{L^2(S,V^*)} \leq d.
    \end{align}
    Eq.~\eqref{ext_eq8:a} and the construction of $w_k$ in
    \eqref{ext_def2} yield 
    \begin{align}\label{ext_eq12}
      \forall \delta>0\quad \forall k\leq \delta:\quad
      \norm{w_k}_{L^2(]\delta,T[, V)} \leq c, 
    \end{align}
    for $c$ independent of $\delta$.  When taken together with
    \eqref{ext_eq11}, \eqref{ext_eq9a}, and the theorem of
    Lions--Aubin, we obtain subsequences
    \begin{subequations}\label{ext_eq13:abc}
      \begin{align}
        w_k &\rightarrow u \quad \text{in}\ L^2(]\delta,T[,
        L^2(\Omega)),  \label{ext_eq13:a}\\
        u_k - w_k &\rightarrow0 \quad\text{in}\
        L^2(S,L^2(\Omega)),  \label{ext_eq13:b}\\ 
        u_k &\rightarrow u \quad\text{in}\ L^2(]\delta, T[,
        L^2(\Omega)). \label{ext_eq13:c} 
      \end{align}
    \end{subequations}
    Then there exist a subsequence $(u_k)$, converging pointwise a.e. on
    $S$ to $u$.  Using \eqref{ext_eq8:a} and Lebesgue's theorem we
    obtain
    \begin{align}\label{ext_eq14}
      u_k\rightarrow u \quad\text{in}\ L^2(S,L^2(\Omega)).
    \end{align}
    
  \item (Limit) We define a translation of the function
    $u_k:[0,T]\rightarrow V$ by
    \begin{align}\label{ext_eq15}
      u_k(t-k) := u^{m-1} \quad\text{on}\ J_m \quad\text{for}\
      m=1,\dots, N,  
    \end{align}
    for $J_m = [(m-1)k,mk]$. Then the semi-discretization yields 
    \begin{align}\label{ext_eq16}
      w'_k(t) + A(t)(u_k(t-k),u_k(t)) = f_k(t) \quad \text{for a.e.}\ t\in
      S\quad \text{in}\ V^*. 
    \end{align}
    Applying \eqref{ext_eq16} to $v\in\V$ with $\V^*\subset
    L^2(S,L^2(\Omega))$ and $v(T) = 0$, and integrating over $S$ using
    integration by parts in the first term yields
    \begin{multline}\label{ext_eq17}
      -\int_0^T(v^\prime(t)\; , \; w_k(t))\, \text{d}t + \int_0^T\langle
      A(t)(u_k(t-k),u_k(t)), v(t)\rangle \, \text{d}t 
      = \int_0^T\langle f_k(t), v(t)\rangle \, \text{d}t + (u_0\; , \;v(0)).
    \end{multline}
    With \eqref{ext_cond2:a} and \eqref{ext_eq8:a} it holds for a
    subsequence that 
    \begin{align}\label{ext_eq18}
      A(\cdot )(u_k(\cdot - k), u_k(\cdot)) \rightharpoonup \zeta \quad \text{in}\ V^*.
    \end{align}
    Taking the limit in \eqref{ext_eq17} we obtain
    \begin{align}\label{ext_eq19}
      -\int_0^T(v^\prime(t)\; , \;u(t))\, \text{d}t +
      \int_0^T\langle\zeta(t), v(t)\rangle\, \text{d}t 
      = \int_0^T\langle f(t),\ v(t)\rangle\, \text{d}t  + (u_0\; , \;
      v(0)). 
    \end{align}
    With $u\in\V$, Eq.~\eqref{ext_eq19} and $u^\prime\in\V^*$ yield
    \begin{align}\label{ext_eq20}
      u^\prime(t) + \zeta(t) = f(t)\quad\text{for a.e.}\ t\in S,\quad
      \text{in}\ V^*. 
    \end{align}
    It is left to show that 
    \begin{align}\label{ext_eq21}
      A(t)(u,u) = \zeta.
    \end{align}
    We use the monotonicity of $A$ in the second argument and
    compactness from the a~priori estimates to obtain
    \begin{align}\label{ext_eq22}
      X_k := \int_0^T\langle A(t)(u_k(t-k),u_k(t)) - A(t)(u_k(t-k), v(t)),
      u_k(t) - v(t)\rangle\, \text{d}t 
      \geq 0
    \end{align}
    for all $v\in \V$. Eq.~\eqref{ext_eq20} then implies
    \begin{align}\label{ext_eq23}
      \int_0^T\langle f(t), u(t)\rangle\, \text{d}t + \frac12\norm{u_0}^2_\Omega
      - \frac12\norm{u(T)}^2_\Omega = \int_0^T\langle\zeta(t),u(t)\rangle\,
      \text{d}t,  
    \end{align}
    and from Eq.~\eqref{ext_eq16} we obtain
    \begin{multline}\label{ext_eq24}
      \int_0^T\langle w'_k(t), u_k(t)\rangle\, \text{d}t + \int_0^T\langle
      A(t)(u_k(t-k),u_k(t), u_k(t)\rangle \, \text{d}t
      = \int_0^T\langle f_k(t), u_k(t)\rangle\, \text{d}t.
    \end{multline}
    We use the following transformation
    \begin{align}
      \int_0^T\langle w'_k(t), u_k(t)\rangle\, \text{d}t &=
      \int_0^T\langle w'_k(t), w_k(t)\rangle\, \text{d}t
      + \int_0^T\langle w'_k(t),u_k(t) - w_k(t)\rangle\, \text{d}t \nonumber\\
      &=\frac12\norm{u^N}^2_\Omega - \frac12\norm{u_0}^2_\Omega -
      \frac12\sum_{m=1}^N\norm{u^m - u^{m-1}}^2_\Omega.
      \label{ext_eq25}
    \end{align}
    Eqs.~\eqref{ext_eq24} and \eqref{ext_eq25} together lead to
    \begin{multline}\label{ext_eq26}
      \int_0^T\langle A(t)(u_k(t-k),u_k(t)), u_k(t)\rangle\, \text{d}t \\ 
      = \int_0^T\langle f_k(t), u_k(t)\rangle\, \text{d}t - \frac12\norm{u^N}^2_\Omega
      + \frac12\norm{u_0}^2_\Omega + \frac12\sum_{m=1}^N\norm{u^m - u^{m-1}}^2_\Omega.
    \end{multline}
    Using Eq.~\eqref{ext_eq5:b}, we conclude that the sum in
    \eqref{ext_eq25} and \eqref{ext_eq26} is convergent. From
    Eq.~\eqref{ext_eq22} we deduce that
    \begin{multline}\label{ext_eq27}
      0\leq \int_0^T\langle f_k(t), u_k(t)\rangle \, \text{d}t + \frac12\norm{u_0}^2_\Omega 
      - \frac12\norm{u^N}^2_\Omega\\
      + \frac12\sum_{m=1}^N\norm{u^m - u^{m-1}}^2_\Omega
      - \int_0^T\langle A(t)(u_k(t-k), u_k(t)), v(t)\rangle\, \text{d}t \\
      - \int_0^T\langle A(t)(u_k(t-k), v(t)), u_k(t) - v(t)\rangle\, \text{d}t.
    \end{multline}
    The limit superior in \eqref{ext_eq27} leads to
    \begin{multline}\label{ext_eq28}
      0 \leq \int_0^T\langle f(t), u(t)\rangle\, \text{d}t 
      + \frac12\norm{u_0}^2_\Omega 
      - \frac12\norm{u(T)}^2_\Omega
      +\gamma\\
      - \int_0^T\langle \zeta(t),v(t)\rangle\, \text{d}t 
      - \int_0^T\langle A(t)(u(t),v(t)), u(t)-v(t)\rangle\, \text{d}t 
    \end{multline}   
    where we used
    \begin{align}\label{ext_eq29}
      \lim_{N\rightarrow\infty} \inf\norm{u^N}^2_\Omega \geq \norm{u(T)}^2_\Omega
    \end{align}
    and
    \begin{align}\label{ext_eq30}
      u_k(\cdot - k) \rightarrow u \quad \text{in}\ L^2(S,L^2(\Omega)) =
      L^2(S\times\Omega) . 
    \end{align}
    This last result follows from Eq.~\eqref{ext_eq14}, the
    Lebesgue integration theory, and an application of the Nemyzki
    operator. Eqs.~\eqref{ext_eq28} and \eqref{ext_eq23} lead to
    \begin{align}\label{ext_eq31}
      -\gamma \leq \int_0^T\langle \zeta(t) - A(t)(u(t),v(t)),
      u(t)-v(t)\rangle\, \text{d}t \quad\forall\ v\in\V.
    \end{align}

    Now we consider two cases: 
    \begin{enumerate}[(a)]
    \item If the integral on the right-hand side of \eqref{ext_eq31}
      is always greater or equal to 0, then let $v=u-\alpha w$ with
      $\alpha>0$ and $w\in\V$. It follows that 
      \begin{align}\label{ext_eq32}
        0\leq \int_0^T\langle\zeta(t) - A(t)(u(t),u(t)-\alpha w(t)), w(t)\,
        \text{d}t \quad \forall w\in\V. 
      \end{align}
      Using condition \eqref{ext_cond2:b}, the limit passage for
      $\alpha\rightarrow0$ is admissible and we obtain
      \begin{align}\label{ext_eq33}
        0\leq \int_0^T\langle \zeta(t)-A(t)(u(t),u(t)),w(t)\rangle\,
        \text{d}t\quad\forall w\in\V,  
      \end{align}
      and with the standard linearity argument the proof is complete.
      
    \item If on the other hand the right-hand side of
      Eq.~\eqref{ext_eq31} becomes negative, then there exists a
      $v\in\V$ and (because of continuity) a whole ball
      $B_r(v)\subset\V$ such that
      \begin{align}\label{ext_eq34}
        \int_0^T\langle \zeta(t) - A(t)(u(t),v(t)),u(t)-v(t)\rangle
        \;\text{d}t \leq 0 \quad\forall v\in B_r(v). 
      \end{align}
      We set $v=u-w$ with $w\in B_r(u-v)$, and use Eq.~\eqref{ext_eq34}
      and the linearity condition \eqref{ext_cond2:b} to obtain
      \begin{align}
        \int_0^T\langle\zeta(t)-A(t)(u(t),u(t)-w(t)),w(t)\rangle\, \text{d}t
        =&\int_0^T\langle \zeta(t)-A(t)(u(t),u(t)),w(t)\rangle \, \text{d}t\nonumber\\
        &\;\; +\;\int_0^T\langle A(t)(u(t),w(t)),w(t)\rangle \,
        \text{d}t \leq 0. \label{ext_eq35}
      \end{align}
      With Eq.~\eqref{ext_cond:c} we deduce
      \begin{align}\label{ext_eq36}
        \int_0^T\langle \zeta(t) - A(t)(u(t),u(t)),w(t)\rangle \,
        \text{d}t\leq 0 \quad\forall w\in B_r(u-v),  
      \end{align}
      which yields
      \begin{align}\label{ext_eq37}
        \int_0^T\langle \zeta(t) - A(t)(u(t),u(t)), \alpha w(t)\rangle
        \, \text{d}t \leq 0 \quad\forall \alpha\geq 0 
        \quad \forall w\in B_r(u-v).
      \end{align}
    \end{enumerate}
    As a result, Eq.~\eqref{ext_eq36} holds for all $w\in\V$ and
    using the standard trick of linearity, statement \eqref{ext_eq21} is
    proven.
  \end{enumerate}
\end{proof}

\subsection{Proof of a~priori estimates}
\label{app:estimates}

\stepcounter{lemma}
\begingroup
\setcounter{tmp}{\value{lemma}}
\renewcommand\thelemma{\ref{sec:analysis}.\arabic{lemma}}
\begin{lemma}
  \lemmatwo
\end{lemma}
\endgroup

\begin{proof}
  Begin by testing Eq.~\eqref{problem_weak_together3} with
  $\varrho_\eps$ to obtain 
  \begin{align*}
    (\partial_t\varrho_\eps,\varrho_\eps)_\Omega + (\kappa_\eps
    D\omega'(\varrho_\eps +
    \omega^{-1}(\Tout))\nabla\varrho_\eps,\nabla\varrho_\eps)_\Omega =
    (-\partial_t\omega^{-1}(\Tout),\varrho_\eps)_\Omega.
  \end{align*}
  Because $D\omega'$ is bounded from below by a positive constant, 
  we can apply the definition of $\kappa_\eps$ to get
  \begin{align*}
    (\partial_t\varrho_\eps,\varrho_\eps)_\Omega +
    \min\{D\omega'\}\norm{\nabla \varrho_\eps}^2_{\Omoneeps}
    +\min\{D\omega'\}\norm{\eps\nabla\varrho_\eps}^2_{\Omtwoeps} \leq
    \norm{\partial_t\omega^{-1}(\Tout)}^2_\Omega +
    \norm{\varrho_\eps}^2_\Omega.
  \end{align*}
  Then, integrating with respect to time and using the boundedness of
  $\norm{\omega^{-1}(\Tout)}^2_\Omega$, we conclude using Gronwall's
  Lemma that 
  \begin{align*}
    \frac12\norm{\varrho_\eps(t)}^2_\Omega +
    \min\{D\omega'\}\norm{\nabla\varrho_\eps}^2_{\Omoneeps,\,t}
    +\min\{D\omega'\}\norm{\eps\nabla\varrho_\eps}^2_{\Omtwoeps,\,t}\leq C
    + \frac12\norm{\varrho_\eps(0)}^2_{\Omega} 
  \end{align*}
  for every $t\in[0,\tend]$, where we use that the initial
  conditions are bounded.  This yields for $\Theta_\eps$ that
  \begin{multline*}
    \frac12\norm{\Theta_\eps(t) + \omega^{-1}(\Tout(t))}^2_{\Omega}
    +\min\{D\omega'\}\norm{\nabla\Theta_\eps}^2_{\Omoneeps,\,t}
    +\min\{D\omega'\}\norm{\eps\nabla\Theta_\eps}^2_{\Omtwoeps,\,t}\\
    \leq C + \frac12\norm{\Theta_\eps(0) +
      \omega^{-1}(\Tout(0))}^2_\Omega, 
  \end{multline*}
  after which we obtain from the reverse triangle inequality that
  \begin{multline*}
    \frac12\norm{\Theta_\eps(t)}^2_{\Omega}
    +\min\{D\omega'\}\norm{\nabla\Theta_\eps}^2_{\Omoneeps,\,t}
    +\min\{D\omega'\}\norm{\eps\nabla\Theta_\eps}^2_{\Omtwoeps,\,t}\\
    \leq C + \frac12\norm{\Theta_\eps(0) + \omega^{-1}(\Tout(0))}^2_\Omega
    + \frac12\norm{\omega^{-1}(\Tout(t))}^2_\Omega.
  \end{multline*}
  This implies for $\HEoneeps$ and $\HEtwoeps$ that
  \begin{align*}
    \norm{\HEoneeps(t)}^2_{\Omoneeps} + \norm{\HEtwoeps(t)}^2_{\Omtwoeps}
    +\min\{D\omega'\}\norm{\nabla\HEoneeps}^2_{\Omoneeps,\,t}
    +\min\{D\omega'\}\norm{\eps\nabla\HEtwoeps}^2_{\Omtwoeps,\,t} \leq C_1, 
  \end{align*}
  where $C_1$ is a constant independent of $\eps$.
\end{proof}

\subsection{Proof of uniqueness theorem}\label{app:uniqueness}

\begingroup
\setcounter{tmp}{\value{theorem}}
\renewcommand\thetheorem{\ref{sec:analysis}.\arabic{theorem}}
\begin{theorem}
  \theoremtwo
\end{theorem}
\endgroup

\begin{proof}
  First we note that the cell problem \eqref{cell_problem1} has a
  unique solution, which was proven in~\citet{TWOSCALE26}. Hence, we
  will only prove uniqueness of the macroscopic problem by assuming
  that there are two solutions $(\HE_{1,\,a},\HE_{2,\,a})$ and
  $(\HE_{1,\,b},\HE_{2,\,b})$, and then showing that they are equal. To
  show uniqueness of solutions to \eqref{s_problem_limit_nonlinear},
  we use the equivalent version \eqref{s_problem_limit_inter} with
  nonlinear diffusion coefficient.  We start by substituting each of
  our solutions into \eqref{s_problem_limit_inter}, subtract the
  two equations, and then test with the functions
  $\HE_{1,\,a}-\HE_{1,\,b}$ and $\HE_{2,\,a}-\HE_{2,\,b}$:
  \begin{multline*}
    \abs{Y^1}(\partial_t\HE_{1,\,a}-\partial_t\HE_{1,\,b},\;
    \HE_{1,\,a}-\HE_{1,\,b})_\Omega \\ 
    + (\Pi (D\omega'(\HE_{1,\,a})\nabla \HE_{1,\,a} -
    D\omega'(\HE_{1,\,b})\nabla \HE_{1,\,b}),\; \nabla \HE_{1,\,a}-\nabla
    \HE_{1,\,b})_\Omega \\ 
    +  (\partial_t\HE_{2,\,a} - \partial_t\HE_{2,\,b},\; \HE_{2,\,a} -
    \HE_{2,\,b})_{\Omega\times Y^2}\\ 
    + (D\omega'(\HE_{2,\,a})\nabla_y \HE_{2,\,a} -
    D\omega'(\HE_{2,\,b})\nabla_y \HE_{2,\,b},\; \nabla_y\HE_{2,\,a} -
    \nabla_y\HE_{2,\,b})_{\Omega\times Y^2} = 0 .
  \end{multline*}
  By adding and subtracting an extra term we obtain 
  \begin{multline*}
    0=\abs{Y^1}(\partial_t\HE_{1,\,a}-\partial_t\HE_{1,\,b},\; \HE_{1,\,a}-\HE_{1,\,b})_\Omega\\ 
    + (\Pi (D\omega'(\HE_{1,\,a})\nabla \HE_{1,\,a} -  D\omega'(\HE_{1,\,a})\nabla \HE_{1,\,b} + D\omega'(\HE_{1,\,a})\nabla
    \HE_{1,\,b} - D\omega'(\HE_{1,\,b})\nabla \HE_{1,\,b}),\\
    \nabla  \HE_{1,\,a}-\nabla \HE_{1,\,b})_\Omega   + (\partial_t\HE_{2,\,a} - \partial_t\HE_{2,\,b},\; \HE_{2,\,a} -  \HE_{2,\,b})_{\Omega\times Y^2}\\ 
    + (D\omega'(\HE_{2,\,a})\nabla_y \HE_{2,\,a}  -   D\omega'(\HE_{2,\,a})\nabla_y \HE_{2,\,b}  + 
    D\omega'(\HE_{2,\,a})\nabla_y  \HE_{2,\,b} - D\omega'(\HE_{2,\,b})\nabla_y \HE_{2,\,b},\\
    \nabla_y\HE_{2,\,a} - \nabla_y\HE_{2,\,b})_{\Omega\times Y^2}  ,
  \end{multline*} 
  which yields the following estimates
  \begin{align*}
    \abs{Y^1}&(\partial_t\HE_{1,\,a} - \partial_t\HE_{1,\,b},\; 
    \HE_{1,\,a}-\HE_{1,\,b})_\Omega +  \min\{\norm{\Pi D\omega'}\}\norm{\nabla
      \HE_{1,\,a} - \nabla \HE_{1,\,b}}^2_\Omega \\ 
    &\quad +  (\partial_t\HE_{2,\,a} - \partial_t\HE_{2,\,b},\;
    \HE_{2,\,a} - \HE_{2,\,b})_{\Omega\times Y^2}  
    + \min\{D\omega'\}\norm{\nabla_y \HE_{2,\,a}  -\nabla
      \HE_{2,\,b}}^2_{\Omega\times Y^2} \\ 
    &\leq -  ((D\omega'(\HE_{1,\,a}) - D\omega'(\HE_{1,\,b}))\nabla \HE_{1,\,b},\; 
    \nabla \HE_{1,\,a}-\nabla \HE_{1,\,b})_\Omega  \\ 
    &\quad - ((D\omega'(\HE_{2,\,a}) - D\omega'(\HE_{2,\,b}))\nabla_y \HE_{2,\,b},\;
    \nabla_y\HE_{2,\,a} - \nabla_y\HE_{2,\,b})_{\Omega\times Y^2}\\ 
    &\leq \norm{D\omega'(\HE_{1,\,a})-D\omega'(\HE_{1,\,b})}_{L^\infty(\Omega)}
    \norm{\nabla \HE_{1,\,b}}_\Omega
    \norm{\nabla \HE_{1,\,a}-\nabla \HE_{1,\,b}}_\Omega\\ 
    &\quad + \norm{D\omega'(\HE_{2,\,a})-D\omega'(\HE_{2,\,b})}_{L^\infty(\Omega\times Y^2)}
    \norm{\nabla_y \HE_{2,\,b}}_{\Omega\times Y^2}
    \norm{\nabla_y\HE_{2,\,a}-\nabla_y\HE_{2,\,b}}_{\Omega \times  Y^2},  \\
    &\leq C_L\norm{D\omega'(\HE_{1,\,a})-D\omega'(\HE_{1,\,b})}_\Omega
    \norm{\nabla \HE_{1,\,b}}_\Omega
    \norm{\nabla \HE_{1,\,a}-\nabla \HE_{1,\,b}}_\Omega\\ 
    &\quad + C_L\norm{D\omega'(\HE_{2,\,a})-D\omega'(\HE_{2,\,b})}_{\Omega\times Y^2}
    \norm{\nabla_y \HE_{2,\,b}}_{\Omega\times Y^2}
    \norm{\nabla_y\HE_{2,\,a}-\nabla_y\HE_{2,\,b}}_{\Omega \times  Y^2},  \\
    &\leq \bigC_DC_L\norm{\HE_{1,\,a}-\HE_{1,\,b}}_\Omega
    \norm{\nabla \HE_{1,\,b}}_\Omega\norm{\nabla \HE_{1,\,a}-\nabla \HE_{1,\,b}}_\Omega\\ 
    &\quad + \bigC_DC_L\norm{\HE_{2,\,a}-\HE_{2,\,b}}_{\Omega\times Y^2}
    \norm{\nabla_y \HE_{2,\,b}}_{\Omega\times
      Y^2}\norm{\nabla_y\HE_{2,\,a}-\nabla_y\HE_{2,\,b}}_{\Omega \times
      Y^2},  
  \end{align*}
  Here, we first use H\"older's inequality; secondly that $D\omega'$ is
  bounded and greater than zero, and $\Omega$ is bounded which implies
  that $\norm{D\omega(\cdot)}_{L^\infty(\Omega)}\leq
  C_L\norm{D\omega(\cdot)}_\Omega$ for a constant $C_L>0$; and thirdly
  we use that $D\omega'$ is Lipschitz continuous with constant
  $\bigC_D$.  Next, we apply the quadratic formula and integrate with
  respect to time to get
  \begin{align*}
    \half\abs{Y^1}&\norm{\HE_{1,\,a}-\HE_{1,\,b}}^2_\Omega + 
    \min\{\norm{\Pi D\omega'}\}\norm{\nabla \HE_{1,\,a} - \nabla
      \HE_{1,\,b}}^2_{\Omega,\,t} \\
    &\qquad\qquad + \half\norm{\HE_{2,\,a} - \HE_{2,\,b}}^2_{\Omega\times
      Y^2} + \min\{D\omega'\}\norm{\nabla_y \HE_{2,\,a} -\nabla
      \HE_{2,\,b}}^2_{\Omega\times Y^2,\,t} \\
    &\leq \frac{\bigC_DC_L \lambda}{2}\int_0^{\tend}\norm{\HE_{1,\,a}-\HE_{1,\,b}}^2_\Omega
    \norm{\nabla \HE_{1,\,b}}^2_{\Omega}\;\text{d}t + \frac{\bigC_DC_L}{2\lambda}\norm{\nabla
      \HE_{1,\,a}-\nabla \HE_{1,\,b}}^2_{\Omega,\,t}\\
    &\quad +
    \frac{\bigC_DC_L\lambda}{2}\int_0^{\tend}\norm{\HE_{2,\,a}-\HE_{2,\,b}}^2_{\Omega\times Y^2}
    \norm{\nabla_y\HE_{2,\,b}}^2_{\Omega\times Y^2}\;\text{d}t + \frac{\bigC_DC_L}{2\lambda} 
      \norm{\nabla_y\HE_{2,\,a}-\nabla_y\HE_{2,\,b}}^2_{\Omega \times
      Y^2,\,t},
  \end{align*}
  for any $\lambda>0$, where we have taken advantage of the fact that
  terms containing the initial conditions are zero. Rearranging terms
  yields
  \begin{align*}
    \half\abs{Y^1}&\norm{\HE_{1,\,a}-\HE_{1,\,b}}^2_\Omega +
    \left( \min\{\norm{\Pi D\omega'}\} -
      \frac{\bigC_DC_L}{2\lambda}\right)\norm{\nabla \HE_{1,\,a} - \nabla
      \HE_{1,\,b}}^2_{\Omega,\,t} \\ 
    &\quad + \half\norm{\HE_{2,\,a} - \HE_{2,\,b}}^2_{\Omega\times 
      Y^2} + \left(\min\{D\omega'\}-
      \frac{\bigC_DC_L}{2\lambda}\right)\norm{\nabla_y \HE_{2,\,a}  -\nabla
      \HE_{2,\,b}}^2_{\Omega\times Y^2,\,t} \\ 
    &\leq \frac{\bigC_DC_L \lambda}{2}\int_0^{\tend}\norm{\HE_{1,\,a}-\HE_{1,\,b}}^2_\Omega
    \norm{\nabla \HE_{1,\,b}}^2_{\Omega}\;\text{d}t + 
    \frac{\bigC_DC_L\lambda}{2}\int_0^{\tend}\norm{\HE_{2,\,a}-\HE_{2,\,b}}^2_{\Omega\times Y^2}
    \norm{\nabla_y\HE_{2,\,b}}^2_{\Omega\times Y^2}\;\text{d}t.
  \end{align*}
  Finally, we choose $\lambda$ large enough such that all terms on the left-hand side
  are positive and exploit that $\norm{\nabla\HE_{1}}^2_\Omega$ and 
  $\norm{\nabla_y\HE_2}^2_\Omega$ are bounded, after which we can apply
  Gronwall's Lemma to obtain  
  \begin{multline*}
    \norm{\HE_{1,\,a} - \HE_{1,\,b}}^2_\Omega + \norm{\nabla \HE_{1,\,a} -
      \nabla \HE_{1,\,b}}^2_{\Omega,\,t} 
    + \norm{\HE_{2,\,a} -
      \HE_{2,\,b}}^2_{\Omega\times Y^2} + \norm{\nabla_y \HE_{2,\,a} -\nabla
      \HE_{2,\,b}}^2_{\Omega\times Y^2,\,t} \leq 0.
  \end{multline*}
  Consequently, $\HE_{1,\,a} = \HE_{1,\,b}$ and $\nabla \HE_{1,\,a} = \nabla
  \HE_{1,\,b}$ almost everywhere on $\Omega\times[0,\tend]$, and similarly
  $\HE_{2,\,a} = \HE_{2,\,b}$ and $\nabla \HE_{2,\,a} = \nabla \HE_{2,\,b}$
  almost everywhere on $\Omega\times Y^2\times[0,\tend]$.
\end{proof}

\section*{Acknowledgements}
IK and JMS were supported by Fellowships from the Alexander von Humboldt
Foundation.  MAP received support from the DFG Priority Pogram 1506. JMS was
funded partially by research grants from the Natural Sciences and
Engineering Research Council of Canada and the North American Maple
Syrup Council.

\bibliographystyle{imamat}
\bibliography{bibliothek}

\end{document}